%% file: smc-convergence-arxiv.tex
\documentclass{amsart}
\input{preamble}

\begin{document}
\title[Time-complexity of sampling form a multimodal distribution]{Time-complexity of sampling from a multimodal distribution using sequential Monte Carlo}
\begin{abstract}
  We study a sequential Monte Carlo algorithm to sample from
  the Gibbs measure with a non-convex energy function at a low temperature.
We use the practical and popular geometric annealing schedule, and use a Langevin diffusion at each temperature level.
  The Langevin diffusion only needs to run for a time that is long enough to ensure local mixing within energy valleys, which is much shorter than the time required for global mixing.
  Our main result shows convergence of Monte Carlo estimators with time complexity that, approximately, scales like the fourth power of the inverse temperature, and the square of the inverse allowed error. 
  We also study this algorithm in an illustrative model scenario where more explicit estimates can be given.
\end{abstract}
\author[Han]{Ruiyu Han}
\address{%
  Department of Mathematical Sciences, Carnegie Mellon University, Pittsburgh, PA 15213.
}
\email{ruiyuh@andrew.cmu.edu}
\author[Iyer]{Gautam Iyer}
\address{%
  Department of Mathematical Sciences, Carnegie Mellon University, Pittsburgh, PA 15213.
}
\email{gautam@math.cmu.edu}
\author[Slep\v cev]{Dejan Slep\v cev}
\address{%
  Department of Mathematical Sciences, Carnegie Mellon University, Pittsburgh, PA 15213.
}
\email{slepcev@andrew.cmu.edu}
\thanks{This work has been partially supported by the National Science Foundation under grants
  DMS-220606,
  DMS-2342349,
  DMS-2406853,
  and the Center for Nonlinear Analysis.}

\maketitle

\tableofcontents

\section{Introduction}

We show that under general non-degeneracy conditions, the Annealed Sequential Monte Carlo algorithm (detailed in Algorithm \ref{a:ASMC}) produces samples from multimodal distributions with time complexity that is a polynomial in the inverse temperature, with a precise dimension independent degree.
We begin (Section~\ref{s:informal}) with an informal description of the algorithm and our results.
Following this, we survey (Section~\ref{s:litReview}) the literature, provide a gentle introduction to the area, and place our work in the context of existing results.
Our main results are stated precisely (Section~\ref{s:results}) below, and the remainder of this paper is devoted to the proofs.

\subsection{Informal statement of main results}\label{s:informal}

Let~$U \colon \mathcal X \to \R$ be an energy function defined on a configuration space~$\mathcal X$.
Consider the Gibbs distribution~$\pi_\epsilon$ whose density is given by
\begin{equation}\label{e:piNu}
  \pi_{\epsilon}(x) = \frac{1}{Z_\epsilon} \tilde \pi_\epsilon(x)
  ,
  \quad\text{where}\quad
  \tilde \pi_\epsilon(x) \defeq e^{-U(x) / \epsilon}
  \text{ and }
  Z_\epsilon \defeq \int_{\mathcal X} \tilde \pi_\epsilon(y) \, d y
  .
\end{equation}
where~$dy$ denotes some fixed measure on the configuration space~$\mathcal X$.
In many applications arising in physics, the parameter~$\epsilon > 0$ is proportional to the absolute temperature.
We adopt (and abbreviate) this terminology and subsequently refer to the parameter $\epsilon$ as the temperature.
In this paper, the space~$\mathcal X$ will typically be the~$d$-dimensional Euclidean space~$\R^d$, or the torus~$\T^d$.

Our aim is to study convergence of an \emph{Annealed Sequential Monte Carlo (ASMC)} algorithm.
This is a \emph{Sequential Monte Carlo (SMC)} algorithm (see for instance~\cite[Chapters 3.3, 17]{ChopinPapaspiliopoulos20}, or~\cite[Chapter 3.4]{Liu08}), where particles are moved through a sequence of interpolating measures obtained by gradually reducing the temperature according to a specified annealing schedule.
We use the practical and popular geometric annealing schedule where the inverse temperatures are linearly spaced~\cite{SyedBouchardCoteEA24}.
  Our main result shows convergence of Monte Carlo estimators using ASMC with time complexity that, approximately, scales like the fourth power of the inverse temperature, and the square of the inverse allowed error.

Before stating our main result, we briefly recall the ASMC algorithm.
\begin{enumerate}[\indent 1.]
  \item
    Choose a finite sequence of temperatures~$\eta_1 > \eta_2 \cdots > \eta_M$ (called an \emph{annealing schedule}) so that~$\pi_{\eta_1}$ is easy to sample from and~$\eta_M = \eta$ is the desired final temperature.

  \item\label{i:mcmc}
    Choose a family of Markov processes~$\set{Y_{\epsilon, \cdot}}_{\epsilon > 0}$ so that for every~$\epsilon > 0$ the stationary distribution of~$Y_{\epsilon, \cdot}$ is~$\pi_\epsilon$, and fix a running time~$T > 0$.

  \item
    Choose arbitrary initial points~$y^1_1$, \dots, $y^1_N$.

  \item
    For each~$i \in \set{1, \dots, N}$, run (independent) realizations of~$Y_{\eta_1, \cdot}$ for time~$T$, starting from~$y^i_1$, to obtain~$x^i_1$.

  \item
    Assign each point~$x^i_1$ the weight~$\tilde \pi_{\eta_{2}}(x^i_1) / \tilde \pi_{\eta_1}(x^i_1)$.
    Choose~$(y^1_2, \dots, y^N_2)$ to be a \emph{resampling} of the points $\paren{x^1_1, \dots, x^N_1}$ from the multinomial distribution with probabilities proportional to the assigned weights.
   
  \item
    Repeat the previous two steps, reducing the temperature until the final temperature is reached.
\end{enumerate}

This is stated more precisely as Algorithm~\ref{a:ASMC} in Section~\ref{s:ASMC}, below.
Clearly if we choose~$T$ larger than the \emph{mixing time} of~$Y_{\eta, \cdot}$ at the final temperature~$\eta = \eta_M$, then the above procedure will produce good samples from~$\pi_{\eta}$.
This, however, is not practical -- when~$U$ has multiple wells the mixing time of~$Y_{\eta, \cdot}$ grows exponentially with~$1/\eta$.
So, when~$\eta$ is small waiting for the mixing time of~$Y_{\eta, \cdot}$ is computationally infeasible.
We will instead show that we only need to choose~$T$ to be larger than the mixing time of~$Y_{\eta_1, \cdot}$ at the \emph{initial temperature} $\eta_1$.
The initial temperature~$\eta_1$ can be chosen to be large enough to ensure the mixing time is computationally tractable.
Moreover, this algorithm requires \emph{only polynomially many} temperature levels~$M$, provided we use the 
\emph{geometric annealing schedule}, \cite{VacherChehabStromme25},
where the inverse temperatures are linearly spaced (and hence the densities form a geometric sequence).

Roughly speaking, our main result is as follows.
\begin{theorem}\label{t:mainIntro}
  Suppose~$U \colon \T^d \to \R$ is a non-degenerate double-well function with wells of equal depth (but not necessarily the same shape).
  For~$\epsilon > 0$ let~$Y_{\epsilon, \cdot}$ be a solution to the overdamped Langevin equation
  \begin{equation}\label{e:Langevin}
    dY_{\epsilon,t} = - \grad U(Y_{\epsilon,t}) \, dt + \sqrt{2 \epsilon} \, dW_t
    ,
  \end{equation}
  where~$W$ is a standard~$d$-dimensional Brownian motion on the torus.
  There exists constants~$C_N, C_T$, depending on~$U$ and~$d$, such that  the following holds.
  For any~$\delta > 0$, $\eta > 0$, choose~$M, N, T$ according to
  \begin{equation}
    M = \ceil[\Big]{\frac{1}{\eta}},
    \quad
    N = \frac{C_N M^2}{\delta^2},
    \quad\text{and}\quad
    T\geq C_T\paren[\Big]{ M^{2}+\log\paren[\Big]{\frac{1}{\delta}}+\frac{1}{\eta} }
  \end{equation}
  and a suitable geometric annealing schedule~$\set{ 1/\eta_k}_{k = 1, \dots, M}$ so that~$\eta_1$ is sufficiently large, and~$\eta_M = \eta$.
  Then the points~$x^1$, \dots, $x^N$ obtained from ASMC (with the parameters above) are such that for any bounded test function~$h$ we have
  \begin{equation}
    \E \paren[\Big]{
	\frac{1}{N}\sum_{i=1}^{N}h(x^i)-\int_{\mathbb T^d}h(x)\pi_{\eta}(x)\,d x}^2
	  < \|h\|_{\osc}^2 \delta^2
	.
  \end{equation}
\end{theorem}

Theorem~\ref{t:mainIntro} shows that the time complexity of obtaining good samples from~$\pi_{\eta}$ using ASMC is polynomial in~$1/\eta$, with a degree independent of dimension (but dimensional coefficients that are not explicit).
We note that the drift in~\eqref{e:Langevin} is independent of temperature $\epsilon$, and so computational complexity of ASMC is proportional to the time complexity.
In contrast, the time complexity of obtaining good samples by directly simulating the process~$Y_{\eta, \cdot}$ is~$e^{O(1/\eta)}$.

The assumption that~$U$ has a double-well structure is mainly to simplify the technical presentation.
Our proof will generalize without difficulty to the situation where~$U$ has more than two wells, at the expense of several technicalities that further obscure the heart of the matter.
As a result we present a detailed proof of Theorem~\ref{t:mainIntro} in the double-well setting (in Section~\ref{s:mainProof}, below).
We sketch the necessary modifications required to address the case when there are more than two wells in Section~\ref{sec:multi}, below.

We assumed that the wells have equal depth above only for simplicity. Our main result (Theorem~\ref{thm: main}) will generalize Theorem~\ref{t:mainIntro} so that it applies to a large class of double-well energy functions, where the low temperature sampling problem is nondegenerate in the sense that each well has a non-negligible fraction of the total mass.
The precise assumptions required are laid out in Section~\ref{s:assumptions}.
In particular, Lemma~\ref{lem: upper bound on mass ratio of two wells} shows that if the wells have \emph{nearly equal depth} (a condition that is necessary for nondegeneracy) then Theorem~\ref{thm: main} applies.

We also remark that Theorem~\ref{t:mainIntro} requires no prior knowledge of the location or the depth of the wells.
In particular, if the target distribution is a mixture, we require no knowledge of the decomposition of the domain into components of a mixture, and only require access to the energy and its gradient.

The main tool used in the proof is a spectral decomposition.
This decomposes any initial distribution into components corresponding to the (target) stationary distribution, a mass imbalance between wells, and higher order terms.
The higher order terms decay exponentially and do not present a problem.
The term corresponding to the mass imbalance decays extremely slowly (at a rate that is exponentially small in the inverse temperature), and is the bottleneck.
This, precisely, is the term that can be eliminated using ASMC.
At high temperatures, all terms converge rapidly, and it is easy to obtain samples with a small mass imbalance.
The resampling step used to move between temperature levels does not disturb this much, resulting in a distribution that has a small mass imbalance at a lower temperature.
Iterating this should, in principle, yield samples from distributions that have a small mass imbalance at every temperature level.
This is the main idea behind of the proof of Theorem~\ref{t:mainIntro}, and is presented in Section~\ref{s:mainProof}, below.

The details of the proof, however, are somewhat involved.
To precisely quantify the error at each level, we require precise bounds on the shape of the second eigenfunction, and how it changes with temperature.
In particular, the proof relies on bounding the inner-product between the normalized second eigenfunction at successive temperature levels, and this involves dimensional constants that are not explicit.
As a result, the constants~$C_N$, $C_T$ in Theorem~\ref{t:mainIntro} are not explicit. The assumption that the state space is the compact torus~$\T^d$, while inconsequential to issues arising from of multimodality, is required to ensure the validity of some of our spectral estimates.

To obtain a better understanding of the dynamics, and a more explicit dimensional dependence of the constants, we also study ASMC in an idealized scenario.
In this idealized scenario, we assume the domain is divided into~$J$ energy valleys, and we have access to a Markov process that mixes quickly in each valley but very slowly globally.
In this case (Theorem~\ref{t:localMixingModel}, in Section~\ref{s:modelResults}, below) we also prove polynomial time complexity bounds, but obtain more explicit constants, and control their dimensional dependence.

We numerically illustrate some aspects of the  performance of ASMC algorithm in Section~\ref{s:numerics}. In particular we highlight how the algorithm adjusts the mass in each valley, which can change as temperature changes and we investigate how the accuracy of the algorithm depends on the number of levels $M$, under fixed computational budget. 
A reference implementation is provided in~\cite{HanIyerEA25code}.

\subsection{Literature review}\label{s:litReview}
We first recall the widely used sampling techniques for nice (e.g. log-concave) measures and then discuss the literature on sampling multimodal distributions. 

\subsubsection{General sampling algorithms}
Perhaps the simplest practical  technique for drawing samples of the target distribution $\pi \propto \exp(-U)$ is based on rejection sampling.
One first draws samples of distribution $\mu$ from which exact samples can be obtained easily (say a Gaussian or a Lebesgue measure on a square), and is such that $\pi$ is absolutely continuous with respect to $\mu$.
One then accepts these samples with probability proportional to $\frac{d\pi}{d\mu}$.
If the measure $\pi$ is much more concentrated than $\mu$ the acceptance probability becomes very small.
For a Gibbs distribution in~$d$-dimensions at temperature~$\eta$, the acceptance rate is typically proportional to $1/\eta^d$, making the cost of this method prohibitively expensive.

Thus in high-dimensions, a different approach is needed.
Most of the widely used methods are based on a stochastic process whose invariant measure is $\pi$.
The largest class of these are Markov Chain Monte Carlo (MCMC) methods which include the seminal Metropolis-Hastings algorithm, Langevin Monte Carlo, Metropolis adjusted Langevin Algorithm (MALA), Hamiltonian Monte Carlo (HMC) and others~\cite{sanzalonso2024coursemontecarlomethods}.
We now briefly recall some of the main algorithms, as any of these can be used in the step~\ref{i:mcmc}  of ASMC (the Markov transition step), provided it rapidly mixes within the modes of the distribution.

The Langevin Monte Carlo (LMC) algorithm  relies on updating individual particles following the overdamped Langevin equation~\eqref{e:Langevin}.
The law of the solution, denoted by $\mu^\epsilon_{t}$, satisfies the Fokker--Planck equation and converges to the stationary distribution~$\pi_\epsilon$ exponentially as~$t \to \infty$.
If the energy function~$U$ is uniformly convex, and satisfies $\alpha I \leq \Hess U$, then it is known that the~$2$-Wasserstein distance converges exponentially with rate $\alpha$ (i.e.\ $W_2(\mu^\epsilon_t, \pi_\epsilon) \leq \exp(-\alpha t) W_2(\mu_0, \pi_\epsilon)$). 
To use this algorithm in practice, one needs to discretize the SDE, which is often done using the explicit  Euler--Maruyama scheme.
Convergence of the time discretized SDE were proved by Vempala and Wibisono~\cite{VempalaWibisono19} with respect to KL divergence, and by Chewi et al.~\cite{Che+24LMC} in~$W_2$.

The general Langevin dynamics allows for inertial effects and is modeled by a system of an SDE for momentum and ODE for position.
This property is the foundation of popular Hamiltonian Monte Carlo (HMC) algorithm, which extends the configuration space to include the momentum variable $p$, and considers Hamiltonian dynamics whose invariant measure is~$\pi^H \propto \exp \big(-U(x) - \frac12 |p|^2 \big)$.
Observe that the first marginal of $\pi^H$ is exactly the target Gibbs measure~$\pi_1$ (with temperature~$\epsilon = 1$).
To numerically obtain samples from~$\pi^H$, the HMC algorithm alternates between the flow of the Hamiltonian dynamics in the phase space, and drawing a new random momentum, whose marginal distribution is a standard Gaussian. The optimal convergence rate for the idealized (i.e.\ one with an exact Hamiltonian dynamics solver) HMC was proved by Chen and Vempala in~\cite{ChenVempala22}.

While the dynamics above have $\pi_1$ as the invariant measure at the continuum level, this is not preserved at the level of numerical schemes resulting in bias that the estimates above control.
The original Metropolis--Hastings algorithm~\cite{MetropolisRosenbluthEA53,Hastings70} offers an algorithm where the target measure $\pi_1$ is the invariant measure at the discrete level.
The algorithm proposes a new sample from a (simple) proposal distribution, and then accepts/rejects the proposal in a manner that ensures the desired target distribution is the invariant measure. This accept/reject step can be combined with several other methods.
In particular when added to LMC one obtains the popular MALA algorithm~\cite{Besag94, RobertsTweedie96, RobertsRosenthal02, Chewi23}.
Other algorithms in this direction include the proximal sampler~\cite{Pereyra16,chenetal2022proximalsampler,Chewi23}.

\subsubsection{Sampling from multimodal distributions}
Sampling from multimodal distributions, is challenge that none of the algorithms like LMC, MALA, or HMC can effectively overcome as their convergence rate becomes extremely slow as the separation between the modes increases. As a result, there is a broad spectrum of works studying algorithms that are suitable for sampling multimodal distributions.
\smallskip

\emph{Annealed importance sampling}, introduced by Neal \cite{Neal01}, involves drawing samples from a sequence of auxiliary distributions starting from starting from one that is easy to sample from, and ending with the target distribution.
Samples are moved from one distribution to the next by a reweighting procedure, and then improved by iterating a Markov chain. The algorithm outputs a set of weighted sample points representing the target distribution.
While this is extremely popular and versatile, one drawback is that the variance of the weights can become extremely large, and with most of the mass being distributed over only a few points~\cite[Chapter 9]{ChopinPapaspiliopoulos20}.
\smallskip

\emph{Sequential Monte Carlo (SMC)} was first developed to study of the average extension of molecular chains~\cite{HammersleyMorton54,RosenbluthRosenbluth55}.
Its use in sampling~\cite{DoucetFreitasEA01,ChopinPapaspiliopoulos20,SyedBouchardCoteEA24} can be seen as generalization of AIS, with the addition of a key resampling step that leads to balanced particle weights. SMC algorithms are enormously popular in a variety of applications and numerous modifications have been developed.

There are a number of works that consider convergence of SMC including obtaining central limit theorems, starting with the work of Chopin \cite{Chopin04}; see also \cite{ChopinPapaspiliopoulos20}.  
As remarked in Section 11.2.4 of \cite{ChopinPapaspiliopoulos20},  the variance of the error typically grows exponentially with the number of levels.
The variance can be controlled~\cite[Section 11.4]{ChopinPapaspiliopoulos20} under restrictive conditions that do not apply to multimodal distributions. 

The works of Schweitzer~\cite{schweizer2012SMC} and Paulin, Jasra, and Thiery~\cite{Paulin19SMC} are the first to rigorously consider the convergence of SMC for multimodal distributions and prove bounds on the variance of the error.
However their assumptions require a strong stability condition on the underlying Markov kernels, which can not be used to in the context of Theorem~\ref{t:mainIntro}.
Building on the coupling technique developed in \cite{Marion23SMC}, Matthews and Schmidler \cite{Mathews24SMC} prove finite sample error bounds for SMC in multimodal setting. Their assumption on the underlying Markov kernels is restrictive requiring knowledge of partition of the domain corresponding to the modes, and also can not be used in the context of Theorem~\ref{t:mainIntro}.

The recent work of of Lee and Santana-Gijzen \cite{lee2024SMC} takes a similar angle as our work in that it shows convergence results for SMC under assumptions of local mixing within the wells and boundedness of ratios of the densities of consecutive levels.
While these assumptions resemble the assumptions we make, there is a key difference: they require a sequence of interpolating measures where the mass in each component of the mixture is known and remains constant.
Devising such an interpolation sequence requires knowing the components of the mixture, which is not available in many practical problems and in particular precludes using interpolations based on adjusting the temperature in the Gibbs measure, such as the geometric annealing we study.

The main differentiating factor between our work, and the SMC papers mentioned above, is that we do not require structural assumptions on the underlying Markov kernels, and do not require any prior knowledge of the mixture components.
As such, our result, stated in Theorem \ref{thm: main}, is the first to provide polynomial time complexity bounds for ASMC using Langevin diffusions and a geometric annealing schedule.
\smallskip

\emph{Parallel, simulated, and related tempering methods.}
Parallel tempering was introduced in a form by Swendsen and Wang \cite{Swendsen1986tempering} and developed by Geyer in \cite{geyer1991markov}.
Simulated tempering introduced by Marinari and Parisi \cite{MarinariParisi92} and developed further by Geyer \cite{Geyer1995}.
These algorithms rely on Markov chains that run on a product space of the desired configuration space and various levels of the temperature.
Samples drawn at a particular value of the temperature may be modified into samples from either a higher, or a lower temperature. At the lowest temperature the marginal of the invariant measure on the product space is the target measure, while at the highest it is a measure where the Markov chain mixes rapidly. 

There are notable results on rigorously showing convergence of parallel and simulated tempering.
In particular, Woodard, Schmidler, and Huber~\cite{WoodardSchmidlerEA09} obtain conditions under which tempering methods are rapidly mixing.
When applied to sampling multimodal distributions the authors considered distributions which have separated modes, but require the variance near each mode to be of size one.
Thus their results do not address the low temperature regime that Theorem~\ref{t:mainIntro} applies to. In \cite{WoodardSchmidlerEA09a} the authors prove that the mixing of these tempering approaches slows exponentially with dimension  if components of multimodal measures have different variances. 
If all the modes have the same shape, Ge, Lee, and Risteski~\cite{GeLeeEA18,GeLeeEA20} show the convergence in TV norm of simulated tempering with error rates that are polynomial in inverse temperature and dimension, provided we have initial estimate on the ratio of the normalizing constants.
The precise degree of the polynomial, however, is not explicitly identified.

Further tempering methods in this family include tempered transitions introduced by Neal in \cite{Neal96a}, which rely on compositions of transitions steps that result in jumps at the lowest temperature and   tempered Hamiltonian Monte Carlo \cite{Nea11}.
Though, to the best of our knowledge, there are no results that apply in the setting of Theorem~\ref{t:mainIntro} and provide polynomial time complexity bounds.
\smallskip

\emph{Annealing without reweighting or resampling.}
There are a number of annealing approaches that evolve a measure from one that is easy to sample from, to the desired target distribution.
In particular, the annealed Langevin Monte Carlo considers Langevin dynamics with slowly changing stationary measure  \cite{guo2025provable,VacherChehabStromme25}. 
These papers show rigorous convergence results for target measures satisfying restrictive structure conditions.
In general the annealed LMC lacks a way to easily adjust the mass within a well at low temperatures.
As a result, the convergence rate is exponentially small in the inverse temperature, and this was rigorously shown for geometric tempering schedule in \cite{VacherChehabStromme25}.
\smallskip

\emph{Further approaches.}
Some recent papers explore new avenues to sampling multimodal distributions. These include approaches based on exploring ideas from diffusion models \cite{VacherChehabStromme25} where the authors show rigorous complexity bounds.
This method, however, suffers from the curse of dimensionality and the error bounds scale like $\delta^d$, where $\delta$ is allowable error and $d$ is dimension. 

 The work \cite{Pompe_20} proposed  a framework of MCMC algorithms for multimodal sampling, which combines an optimization step to find the modes with Markov transition steps. 
  They showed the weak law of large numbers for the Monte Carlo integral using samples generated by the  Auxiliary Variable Adaptive MCMC algorithm.

Another direction explored is to use ensemble methods that involve Markov Chains whose jump rates use estimation of the density of the measure represented by the particle configuration \cite{lu2019acceleratinglangevinsamplingbirthdeath,Lindsey22teleporting,Lu_2023}. 
These approaches can be seen as particle approximations of gradient flows of KL divergence in spherical Hellinger metric, which converge exponentially fast with rate that is independent of the height of the barrier.
However this method also suffers from the curse of dimensionality, as the kernel density estimation used to estimate density based on the configuration of particles introduces bias that becomes large in high dimensions.

A few methods modify~\eqref{e:Langevin} in a manner that allows particles to move between modes faster.
The authors of~\cite{engquist2024samplingadaptivevariancemultimodal} do this by modifying the diffusion, and the authors of~\cite{ReyBelletSpiliopoulos15,DamakFrankeEA20,ChristieFengEA23} do this by introducing an additional drift term.
In both cases the modified equation has terms that grow exponentially with the inverse temperature, and a numerical implementation is computationally expensive.

\subsection*{Plan of the paper}

In Section~\ref{s:results} we precisely state our algorithm, and state results guaranteeing convergence both for ASMC in an idealized scenario (Theorem~\ref{t:localMixingModel}), and for a double-well energy function (Theorem~\ref{thm: main}, which generalizes Theorem~\ref{t:mainIntro}).
For the idealized scenario we are able to obtain explicit constants, and track the dimensional dependence (Proposition~\ref{prop:CrEg}).
Numerical simulations illustrating relevant aspects of the performance of ASMC in model situations are shown in Section~\ref{s:numerics}.
We prove Theorems~\ref{t:localMixingModel} and~\ref{thm: main} in Sections~\ref{s:localModel} and~\ref{s:mainProof} respectively.
The proof of Theorem~\ref{t:localMixingModel} relies on a few lemmas which are proved in Section~\ref{s:localmodelproof}.

The proof of Theorem~\ref{thm: main} is a little more involved and relies on certain spectral properties and Langevin error estimates.
The error estimates are proved in Sections~\ref{s:langevinErrorProof} and Section~\ref{sec:ProofIteration}. In Section~\ref{s:CBV} we show that that regular enough energy functions satisfy the assumptions required for Theorems~\ref{t:localMixingModel} and~\ref{thm: main}, and prove the dimensional dependence stated in Proposition~\ref{prop:CrEg}.
In Section~\ref{sec:propertycheck} we verify the spectral properties required in the proof of Theorem~\ref{thm: main}.
Finally in Section~\ref{sec:multi} we sketch the modifications required to generalize Theorem~\ref{thm: main} to the multi-well setting.

\subsection*{Acknowledgments}
We thank Seungjae Son for helpful comments and for spotting a minor error in the first version of this manuscript.
We also thank the anonymous referee for several extremely helpful comments that substantially improved this paper.

\section{Main results}\label{s:results}
\subsection{Annealed Sequential Monte Carlo (ASMC)}\label{s:ASMC}
We now briefly introduce the ASMC algorithm, which is stated precisely as Algorithm~\ref{a:ASMC}, below.
In many situations of interest, the configuration space~$\mathcal X$ admits a decomposition into \emph{energy valleys}.
MCMC samplers (such as~\eqref{e:Langevin}) are typically confined to an energy valley for time~$e^{O(1/\epsilon)}$ before moving to a different valley (see for instance~\cite{Arrhenius89}).
Of course, waiting time~$e^{O(1/\epsilon)}$ to explore the state space is practically infeasible, and directly using an MCMC sampler is prohibitively slow at low temperatures.

Annealing and tempering (both terms having origin in metallurgy and describing heat treatment of metals) based algorithms, in particular ASMC we study, have been introduced to overcome the issue of slow global mixing of the MCMC algorithms. ASMC a special case of a \emph{sequential Monte Carlo} algorithm, as samples are drawn in sequence from an auxiliary family of distributions, starting from one that is easy to sample from and ending with the target distribution. The name ASMC stems from the fact that the auxiliary family of distributions used are obtained by starting from the Gibbs distribution at a high temperature, and then gradually lowering the temperature until the desired temperature is reached.

To use ASMC, we choose an \emph{annealing schedule}, which is a sequence of temperatures~$\eta_1 > \eta_2 \cdots > \eta_M$, chosen so that the MCMC sampler converges fast at temperature~$\eta_1$, the desired final temperature is~$\eta_M = \eta$.
Samples at temperature~$\eta_k$ are transformed to samples at temperature~$\eta_{k+1}$ by reweighting them with the ratio of densities~$\pi_{\eta_{k+1}} / \pi_{\eta_k}$.
To ensure the mass is spread across sample points, the weights are redistributed using a \emph{resampling} process.
The samples are then improved by iterating an MCMC sampler for a fixed amount of time, and then the above processes is repeated at the next temperature until the final temperature is reached. It is important to note that for the reweighting step, one \emph{does not} have access to the normalized densities~$\pi_{\eta_k}$ in practice, as the normalization constants are not known and are hard to compute.
However, using weights proportional to the ratio of the normalized densities~$\pi_{\eta_{k+1}} / \pi_{\eta_k}$ is \emph{equivalent} to using weights proportional to the ratio of the unnormalized densities~$\tilde \pi_{\eta_{k+1}} / \tilde \pi_{\eta_k}$.
The unnormalized densities are known, and are used in the reweighting step instead of the normalized densities.

We now describe the \emph{resampling} step: given points~$x^1_k$, \dots, $x^N_k$ which are (approximate) samples from~$\pi_{\eta_k}$, we obtain~$y^1_{k+1}$, \dots, $y^N_{k+1}$ by resampling from the  points~$\set{x^1_k, \dots, x^N_k}$ using the multinomial distribution with probabilities 
\begin{equation}\label{e:rkTilde}
  \P( y^i_{k+1} = x^j_k ) = \frac{\tilde r_k(x^j_k)}{\sum_{n=1}^N \tilde r_k(x^n_k)}
  ,
  \quad\text{where}\quad \tilde r_k
    \defeq \frac{\tilde \pi_{\eta_{k+1}}}{\tilde \pi_{\eta_k}}.
\end{equation}
Some points may be repeated or lost.
Nevertheless, an elementary heuristic (explained in Section~\ref{s:locMMProof}, after~\eqref{e:PyiEqXj}, below) suggests that the new points~$y^1_{k+1}$, \dots, $y^N_{k+1}$ should be good samples from~$\pi_{\eta_{k+1}}$.
\begin{algorithm}[htb]
  \caption{Annealed Sequential Monte Carlo (ASMC) to sample from~$\pi_\eta$.}
  \label{a:ASMC}
  \begin{algorithmic}[1]\reqnomode
    \Require Temperature~$\eta$, energy function~$U$, and Markov processes~$\set{Y_{\epsilon, \cdot}}_{\epsilon \geq \eta}$ so that the stationary distribution of~$Y_{\epsilon, \cdot}$ is~$\pi_\epsilon$.
    \item[\textbf{Tunable parameters:}]
    \Statex
      \begin{enumerate}[(1)]
        \item Number of levels~$M \in \N$, and annealing schedule~$\eta_1 > \cdots > \eta_M = \eta$.
	\item Number of realizations~$N \in \N$, and initial points~$y^1_1$, \dots, $y^N_1 \in \mathcal X$.
	\item Level running time~$T > 0$.
      \end{enumerate}

    \For{$k \in \set{1, \dots, M-1}$}

      \State\label{i:exploreValley}
	For each $i \in \set{1, \dots N}$, simulate~$Y_{\eta_k, \cdot}$ for time~$T$ starting at~$y^i_k$ to obtain~$x^i_k$.

      \State\label{i:resampling}
	Choose $(y^1_{k+1}, \dots, y^N_{k+1})$ by resampling from $\set{x^1_k, \dots, x^N_k}$ using the multinomial distribution with probabilities given by~\eqref{e:rkTilde}.

    \EndFor

    \State
      For each $i \in \set{1, \dots N}$, simulate~$Y_{\eta_M, \cdot}$ for time~$T$ starting at~$y^i_M$ to obtain~$x^i$.
    \State\label{i:RTlast}%
      \Return $(x^1, \dots, x^N)$.
  \end{algorithmic}
\end{algorithm}

\begin{remark}
Instead of resampling at every step, modern, practical algorithms typically  control the variance of the weights using more sophisticated resampling procedures. A popular approach is to introduce a measure of the quality of the weight distribution and only resample when the quality becomes lower than a desired threshold, which is called \emph{adaptive resampling}. For this and other resampling approaches see, for instance, the books \cite[Sections 10.2]{ChopinPapaspiliopoulos20} or ~\cite[Chapter 3.4]{Liu08}.
\end{remark}

We now provide a brief heuristic explanation as to why one may be able to obtain good quality samples in polynomial time using Algorithm~\ref{a:ASMC}.
First, since~$\eta_1$ is large and the process~$Y_{\eta_1, \cdot}$ mixes quickly, and so the distribution of~$x^1_1$, $x^2_1$, \dots, $x^N_1$ will be close to the Gibbs measure~$\pi_{\eta_1}$.
Now the resampling step may produce degenerate samples with several repeated points.
However the fraction of points in each energy valley will be comparable to the~$\pi_{\eta_{2}}$-mass of the same valley.
In the situation we consider, the main bottleneck to fast mixing is moving mass between valleys.
Since the samples at temperature~$\eta_2$ have approximately the right fraction of mass in each energy valley, the distribution after running~$Y_{\eta_2, \cdot}$ for time~$T$ will be close to the Gibbs distribution~$\pi_{\eta_2}$.
Repeating this argument should iteratively yield good samples at the desired final temperature~$\eta_M$.

A rigorous proof of the above quantifying the convergence rate, however, requires some care.
The number of levels~$M$ is large (grows linearly in the inverse temperature), and the error going from level~$k$ to~$k+1$ accumulates \emph{multiplicatively}.
Nevertheless, we will show that if~$\eta_1, \dots, \eta_M$ according to the geometric annealing schedule, then the total error accumulates slowly enough that Algorithm~\ref{a:ASMC} produces good samples in time that is polynomial in~$1/\eta$.
Carrying out the details of this heuristic for a double-well energy function using Langevin diffusions as the MCMC sampler (as described in Theorem~\ref{t:mainIntro}) is technical, and requires several model specific bounds that distract from the main idea.
Thus, we first consider an illustrative model problem where we can study Algorithm~\ref{a:ASMC}, and then revisit it in the context of Theorem~\ref{t:mainIntro}.

\subsection{ASMC for a Local Mixing Model}\label{s:modelResults}
We now present an idealized scenario where we can analyze Algorithm~\ref{a:ASMC} quantitatively, and obtain explicit constants in our error estimates.
Suppose the number of components of the multimodal measure, ~$J \geq 2$, and the domain~$\mathcal X$ can be partitioned into~$J$ domains~$\Omega_1$, \dots, $\Omega_J$.
We are interested situations where we have access to a process~$Y_{\epsilon, \cdot}$ that mixes quickly in each domain~$\Omega_j$, however, transitions very slowly between domains and hence mixes slowly overall.

To model this behavior, for every~$\epsilon > 0$ let~$\chi_\epsilon \in (0, 1)$ denote probability of staying in the same domain after time~$1$.
Let~$Y_{\epsilon, \cdot}$ be the discrete time Markov process defined as follows.
At time~$n \in \N$, let~$j$ be the unique element of~$\set{1, \dots, J}$ such that~$Y_{\epsilon, n} \in \Omega_j$.
Flip an independent coin that lands heads with probability~$\chi_\epsilon$ and tails with probability~$1 - \chi_\epsilon$.
If the coin lands heads, we choose~$Y_{\epsilon, n+1} \in \mathcal X$ independently from the distribution~$\pi_\epsilon$.
If the coin landed tails, we choose~$Y_{\epsilon, n+1} \in \Omega_j$ independently from the distribution with density
\begin{equation}
  \frac{\pi_\epsilon \one_{\Omega_j} }{\pi_\epsilon(\Omega_j)}
  \,.
\end{equation}
In other words, $Y_{\epsilon, \cdot}$ is the Markov process whose one step transition density is
\begin{equation}\label{e:YTransitionKernel}
  p^\epsilon_1(x, y) =
      (1 - \chi_\epsilon) \pi_\epsilon(y)
      +  \chi_\epsilon \sum_{j = 1}^J \one_{\set{x,y \in \Omega_j}} 
      \frac{\pi_\epsilon(y)}{\pi_\epsilon(\Omega_j)}
	\,.
\end{equation}

Notice that the expected transition time between domains is the bottleneck to mixing, and is of order~$1/(1 - \chi_\epsilon)$.
One situation of interest, is when $\chi_\epsilon$ is extremely close to~$1$ (for instance $\chi_\epsilon \approx \exp(e^{-O(1/\epsilon)})$).
This models the behavior that arises several applications of interest, including Langevin dynamics driven by the gradient of an energy function with multiple wells, and this is studied in detail in Section~\ref{s:langevin}, below.
In such situations waiting for time~$1/(1-\chi_\epsilon)$ is prohibitively expensive when~$\epsilon$ is small, and can not be done in practice.

Suppose now we are interested in computing Monte Carlo integrals with respect to the Gibbs distribution~$\pi_\eta$ for some small temperature~$\eta > 0$.
A direct Monte Carlo approach simulating~$Y_{\eta, \cdot}$ is unfeasible as it requires simulating~$Y_{\eta, \cdot}$ for time~$O(1 / (1 - \chi_\eta))$, which very long when~$\eta$ is small.
We now show that Algorithm~\ref{a:ASMC}, with a judicious choice of parameters, makes this time an order of magnitude smaller.

\begin{theorem}\label{t:localMixingModel}
  Suppose for some~$0 \leq \eta_{\min} < \eta_{\max} \leq \infty$ we have
  \begin{equation}\label{e:CBV}
    C_\LBV \defeq \sum_{j = 1}^J \int_{\eta_{\min}}^{\eta_{\max}}
    \abs{\partial_\epsilon \ln \pi_\epsilon(\Omega_j)} \, d\epsilon < \infty
    .
  \end{equation}
  For any finite~$\eta_1 \in (\eta_{\min}, \eta_{\max}]$, $\delta, \eta, \nu > 0$ with~$\eta \in [\eta_{\min}, \eta_1)$, and constants~$C_T$, $C_N > 0$ choose~$M, N, T \in \N$ so that%
  \footnote{
    Here~$\tmix(Y_{\epsilon, \cdot}, \delta)$ denotes the~$\delta$-mixing time of the process~$Y_{\epsilon,\cdot}$ (see for instance~\cite{LevinPeres17}), and measures the~$\TV$-rate of convergence of~$Y_{\epsilon, \cdot}$ to the stationary distribution~$\pi_\epsilon$.
    Explicitly, if~$p^\epsilon_n$ denotes the~$n$-step transition density of~$Y_{\epsilon, \cdot}$, then the~$\delta$-mixing time is given by
    \begin{equation}
      \tmix( Y_{\epsilon, \cdot}, \delta )
	\defeq \inf\set[\Big]{ n \in \N \st \sup_{x \in \mathcal X} \norm{p_n^\epsilon(x, \cdot) - \pi_\epsilon(\cdot) }_{L^1} \leq 2\delta }
	.
    \end{equation}
  }
  \begin{equation}\label{eq: parametersTNloc}
    M \geq \ceil*{\frac{1}{\nu \eta}},
    \quad
    N \geq \frac{C_N M^2}{\delta^2},
    \quad\text{and}\quad
    T \geq \tmix\paren[\Big]{ Y_{\eta_1, \cdot}, \frac{\delta}{C_T}}
    ,
  \end{equation}
  and choose~$\eta_2, \dots, \eta_M$ so that~$\eta_M = \eta$ and~$1/\eta_1$, \dots, $1/\eta_M$ are linearly spaced.

  For every~$\delta,\nu > 0$, there exists (explicit) constants~$C_N = C_N(U/\eta_1, J, \nu)$, and $C_T = C_T(U/\eta_1, J, \nu)$ such that if the process~$Y_{\epsilon, \cdot}$ in Algorithm~\ref{a:ASMC} have transition density~\eqref{e:YTransitionKernel}, and if the parameters to Algorithm~\ref{a:ASMC} are chosen as in~\eqref{eq: parametersTNloc}, then for every bounded test function~$h$, and arbitrary initial data~$\set{x^i_0}$, the points~$(x^1, \dots, x^N)$ returned by Algorithm~\ref{a:ASMC}  satisfy
  \begin{equation}\label{e:MCErrorModel}
    \norm[\Big]{
	\frac{1}{N}\sum_{i=1}^{N}h(x^i)-\int_{\mathcal X}h(x)\pi_{\eta}(x)\,d x
      }_{L^2(\P)}
	  < \|h\|_{\osc} \delta
	.
  \end{equation}
\end{theorem}

We prove Theorem~\ref{t:localMixingModel} in Section~\ref{s:localModel}, below.

\begin{remark}
  In~\eqref{e:MCErrorModel} above, we clarify that the Monte Carlo sum~$\frac{1}{N} \sum_1^N h(x^i)$ is a random variable, as the points~$x^i$ are random, and the notation~$\norm{\cdot}_{L^2(\P)}$ denotes the~$L^2(\P)$ norm with respect to the underlying probability measure~$\P$.
  Explicitly, if~$X$ is a random variable, then~$\norm{X}_{L^2(\P)} = (\E X^2)^{1/2}$.
\end{remark}

\begin{remark}\label{r:etaMinIndep}
  Inequality~\eqref{e:CTCN} bounds the constants~$C_N$ and~$C_T$ explicitly in terms of~$J$, $\nu$, $C_\LBV$ and the constant~$C_r$ (defined in~\eqref{eq:defCr}, below).
  The reason we state Theorem~\ref{thm: main} allowing~$\eta_{\min} > 0$ is because we will show that for a class of double well potentials, the constant~$C_\LBV$ can be bounded	independent of~$\eta_{\min}$.
  Explicitly, if the energy function~$U$ has wells whose depth differ by~$O(\eta_{\min})$, then the mass in each well can be bounded away from zero, independent of~$\eta_{\min}$, and we will show (Corollary~\ref{cor:eqcbvholds}, below) $C_N$ and~$C_T$ are polynomials in the inverse temperature~$1/\eta$, with a dimension independent degree and coefficients that are independent of~$\eta_{\min}$.

  If, on the other hand, the difference between the well depths is much larger than the minimum temperature~$\eta_{\min}$, then the as~$\epsilon \to \eta_{\min}$ some of the domains~$\Omega_j$ will contain an exponentially small fraction fraction of the total mass.
  Hence the multimodal nature of the target distribution degenerates, and the sampling from this distribution requires the simulation of rare events.
  This goes beyond the scope of the present work and Theorem~\ref{t:localMixingModel} does not apply.
\end{remark}

\begin{remark}
  The constant~$C_{\LBV}$ controls the change in the mass~$\pi_k(\Omega_j)$ between levels.
  The change in the mass of each well with temperature has played a role in related convergence results, and was used by Woodard et\ al.\ in~\cite{WoodardSchmidlerEA09}.
  Indeed, Corollary~3.1 in~\cite{WoodardSchmidlerEA09} controls the spectral gap of a tempering chain by the quantity~$\gamma(\mathcal A)$, defined by
\begin{equation}\label{eq:gammaA}
  \gamma(\mathcal A) \defeq \min_{j=1,\dots,J}\prod_{k=1}^{M-1}\min\set[\Big]{1,\frac{\pi_{k}(\Omega_{j})}{\pi_{k+1}(\Omega_j)}}
  .
\end{equation}
  To relate this to~$C_\LBV$, we note that the  AM-GM inequality implies a lower bound for~$\gamma(\mathcal A) \geq \exp(-C_\LBV)$ for any selected annealing sequence~$\eta_k$.
\end{remark}

\begin{remark}
Theorem~\ref{t:localMixingModel} shows that the \emph{averaged} empirical measure is~$\TV$ close to the Gibbs distribution.
  Explicitly, the averaged empirical measure~$\mu$ is defined by
  \begin{equation}
    \mu(A) \defeq
      \frac{1}{N} \E \sum_{i=1}^N \delta_{x^i}(A)
      = \frac{1}{N} \sum_{i = 1}^N \P( x^i \in A )
    ,
  \end{equation}
  where~$x^1$, \dots, $x^N$ are the points returned by Algorithm~\ref{a:ASMC}.
  Now Theorem~\ref{t:localMixingModel} and Jensen's inequality immediately imply~$\abs{\mu(A) - \pi_\eta(A)} \leq \delta$ for every Borel set~$A$, and hence
  \begin{equation}
    \norm{\mu - \pi_\eta}_\TV \leq \delta
    .
  \end{equation}
\end{remark}

\subsubsection*{Computational Complexity}
We now  estimate the computational complexity of Monte Carlo integration using Theorem~\ref{t:localMixingModel} and compare it to the direct approach using the process~$Y_{\eta, \cdot}$.
We only address the asymptotic behavior of the complexity in the small temperature regime, and we do not quantify the dependence of the complexity on other factors such as the shape of~$U$.

We assume the computational complexity of simulating the process~$Y^\epsilon$ for time~$T$ is~$O(T)$.
Using an alias method~\cite{Vose91} the complexity of the resampling step is~$O(N)$, which makes the computational complexity of Algorithm~\ref{a:ASMC} to be of order~$M N T$.
To estimate~$T$, we need to estimate the~$\delta$-mixing time of the process~$Y_{\epsilon, \cdot}$ at the initial temperature~$\epsilon = \eta_1$.
For this, we use~\eqref{e:YTransitionKernel} to deduce that the~$n$-step transition density of~$Y_{\epsilon, \cdot}$ is
\begin{equation}
  p^\epsilon_n(x, y) =
      (1 - \chi_\epsilon^n) \pi_\epsilon(y)
      +  \chi_\epsilon^n \sum_{j = 1}^J \one_{\set{x,y \in \Omega_j}} 
      \frac{\pi_\epsilon(y)}{\pi_\epsilon(\Omega_j)}
	\,,
\end{equation}
which implies
\begin{equation}
  \tmix( Y_{\epsilon, \cdot}, \delta )
    \leq \frac{\ln (\delta/2)}{\ln \chi_\epsilon}.
\end{equation}
Thus Theorem~\ref{t:localMixingModel} implies the computational complexity of running Algorithm~\ref{a:ASMC} to achieve the Monte Carlo error~\eqref{e:MCErrorModel} is
\begin{equation}\label{e:ASMCModelCost}
  \text{complexity(Algorithm~\ref{t:localMixingModel})}
  =
  O(MNT)
    \leq  \frac{C(U)\abs{\ln \delta}}{\eta^3 \delta^2 \abs{\ln \chi_{\eta_1}} }
  ,
\end{equation}
for some~$U$-dependent constant~$C(U)$.

On the other hand, achieving the same Monte Carlo error by simulating independent realizations of~$Y_{\eta, \cdot}$ has a computational cost of
\begin{equation}\label{e:MCModelCost}
  \text{complexity(Direct Monte Carlo)}
  = O\paren[\Big]{\frac{\ln \delta}{\delta^2 \abs{\ln \chi_\eta}}}
  .
\end{equation}
We note that the complexity of ASMC in~\eqref{e:ASMCModelCost} involves the mixing time at the initial temperature~$\eta_1$, which is small, and a polynomial in the final temperature~$\eta$.
In contrast, the complexity of direct Monte Carlo~\eqref{e:MCModelCost} involves the mixing time at the \emph{final} temperature~$\eta$, which is typically exponential in~$1/\eta$.

\subsection{ASMC for a double-well energy function}\label{s:langevin}

We now study Algorithm~\ref{a:ASMC} when the configuration space~$\mathcal X$ is the~$d$-dimensional torus~$\T^d$.
Here the Gibbs measure~$\pi_\epsilon$ arises naturally as the stationary distribution of the \emph{overdamped Langevin equation}~\eqref{e:Langevin}.

When~$U$ is convex, the process~$Y$ mixes quickly even in high dimensions~\cite{bakry2014analysis}, and this provides a very efficient way to sample from the Gibbs distribution~$\pi_\epsilon$.
When~$U$ is not convex, however, the process~$Y_{\epsilon, \cdot}$ mixes extremely slowly.
In fact, the well known Arrhenius law \cite{Arrhenius89} states
that in general it takes time~$t \approx e^{C / \epsilon}$ before the distribution of~$Y_{\epsilon, \cdot}$ becomes close to the Gibbs distribution~$\pi_\epsilon$.
At low temperatures, this is too long to be practical.

The reason Langevin dynamics mixes so slowly is because the drift in~\eqref{e:Langevin} pulls trajectories towards local minima of~$U$.
In order to escape an energy valley, the noise term in~\eqref{e:Langevin} has to go against the drift for an~$O(1)$ amount of time, which happens with \emph{exponentially small} probability.
In each energy valley, however, the energy function ~$U$ is essentially convex which makes the process~$Y_{\epsilon, \cdot}$ mix quickly in valleys.
We also note that the situation considered in Section~\ref{s:modelResults} is an idealized model for the dynamics of~\eqref{e:Langevin}.

We study Algorithm~\ref{a:ASMC} for target distributions corresponding to  double-well energy functions and show that appropriate choice of parameters allows one to compute integrals with respect to the Gibbs distribution, with time complexity that is polynomial in the inverse temperature.

We again remark that the assumption that~$U$ is a double-well energy function is mainly to simplify the presentation, and the generalization to energy functions with more wells is straightforward.

\begin{theorem}\label{thm: main}
  Suppose for some~$0 \leq \eta_{\min} < \eta_{\max} \leq \infty$, the function~$U$ is a double-well function that satisfies Assumptions~\ref{a:criticalpts}, \ref{assumption: nondegeneracy} and~\ref{a:massRatioBound} in Section~\ref{s:mainProof} below.
  Let~$\hat \gamma_r \geq 1$ be the ratio of the saddle height to the energy barrier, defined precisely in~\eqref{e:gammaHatRDef}, below.
  Given~$\eta_1 \in (\eta_{\min}, \eta_{\max}]$ finite, $\alpha, \delta, \eta, \nu > 0$ with~$\eta \in [\eta_{\min}, \eta_1)$, and constants~$C_T$, $C_N > 0$ choose~$M, N \in \N$, and~$T \in \R$ so that%
  \begin{equation}\label{e:MTN}
    M \geq \ceil*{ \frac{1}{\nu \eta} },
    \quad
    T\geq C_T\paren[\Big]{ M^{(1 + \alpha) \hat \gamma_r}+\log\paren[\Big]{\frac{1}{\delta}}+\frac{1}{\eta} }
    \quad\text{and}\quad
    N \geq \frac{C_{N} M^2}{\delta^2},
  \end{equation}
  and choose~$\eta_2, \dots, \eta_M$ so that~$\eta_M = \eta$ and~$1/\eta_1$, \dots, $1/\eta_M$ are linearly spaced.

  For every~$\alpha, \delta, \nu > 0$, there exist constants~$C_T = C_T(\alpha, \nu, U/\eta_1)$ and~$C_{N}(\nu, U/\eta_1)$ such that if the process~$Y_{\epsilon , \cdot}$ in Algorithm~\ref{a:ASMC} is given by~\eqref{e:Langevin}, and the parameters to Algorithm~\ref{a:ASMC} are chosen as in~\eqref{e:MTN}, then for every bounded test function~$h$, and arbitrary initial data~$\set{x^i_0}$, the points~$(x^1, \dots, x^N)$ returned by Algorithm~\ref{a:ASMC}  satisfy
  \begin{equation}\label{e:MCerror}
      \norm[\Big]{
	\frac{1}{N}\sum_{i=1}^{N}h(x^i)-\int_{\mathbb T^d}h(x)\pi_{\eta}(x)\,d x}_{L^2(\P)}
	  < \|h\|_{\osc} \delta
	.
  \end{equation}
\end{theorem}

We remark that Assumptions~\ref{a:criticalpts}--\ref{a:massRatioBound} are nondegeneracy assumptions, and do not require symmetry, or similarity of the shape of the wells.
The proof of Theorem~\ref{thm: main} follows the same general strategy as that of Theorem~\ref{t:localMixingModel}, however the details more technically involved.
In Theorem~\ref{t:localMixingModel} the main idea is to show that if at level~$k$, the initial mass distributions in the domains~$\Omega_1$, \dots, $\Omega_J$ is distributed according to~$\pi_{\eta_k}$, then the process~$Y_{\eta_k, \cdot}$ will correct the shape and quickly give a distribution that is close to the Gibbs distribution~$\pi_{\eta_k}$.
To show this in the context of~\eqref{e:Langevin}, we consider a spectral decomposition based on eigenvalues of the generator of~\eqref{e:Langevin}.
We will show that if the projection of the initial distribution onto the \emph{second} eigenspace is small, then the Langevin dynamics will quickly correct the shape and yield a distribution close to the Gibbs measure.
The proof of this involves several technical lemmas controlling the shape of the eigenfunctions and introduces dimensional pre-factors that are not explicit.
This takes up the bulk of the paper and begins in Section~\ref{s:mainProof}, below.

\subsubsection*{Time and computational complexity.}
We now briefly discuss the computational cost of integration using Theorem~\ref{thm: main}.
Suppose~$U$ is a double-well function with wells of equal depth, so that~$\hat \gamma_r$ is exactly~$1$.
As mentioned earlier, the resampling step costs~$O(N)$ and so the time complexity of running Algorithm~\ref{a:ASMC} (for $\eta<1$ with  $\nu=1$) to achieve the error tolerance~\eqref{e:MCerror} is 
\begin{equation}\label{e:CostASMCLangevin}
  O(MTN) \leq \frac{\tilde C_d}{\eta^3 \delta^2 }\, T 
    \leq \frac{C_d}{\eta^3 \delta^2} \paren[\Big]{
      \frac{1}{\eta^{1 + \alpha}}+\log\paren[\Big]{\frac{1}{\delta}}
    }
    \,,
\end{equation}
for some dimensional constants~$\tilde C_d,C_d$. The second inequality gives us the precise, polynomial, time complexity of the algorithm.
The significance of the first inequality is that the computational complexity of the algorithm is, up to dimensional constant, $\frac{1}{ \eta^3 \delta^2}$ times the computational complexity of the numerical algorithm which mixes the distribution sufficiently well within the wells. 

In order to use Algorithm~\ref{a:ASMC} in practice, one has to time discretize~\eqref{e:Langevin} and consider the bias induced by this discretization.
Obtaining the precise errors for numerical discretizations of LMC and other algorithms is an active area of research, and we refer the reader to the notes by Chewi~\cite{Chewi23} for comprehensive overview. 
Obtaining rigorous computational complexity of ASMC is a challenging open problem, as the wells are not exactly log-concave and one would need to control various terms in our proof up to discretization error. We remark, however, that in our formulation the drift in~\eqref{e:Langevin} is independent of the temperature~$\epsilon$, and so for our purposes, the number of iterations required to simulate~\eqref{e:Langevin} for a given length of time~$T$ is proportional to time~$T$ and independent of the temperature~$\epsilon$.
Let us also remark that, assuming the scaling for smooth log-concave wells can be reached the estimates of Chewi \cite[Theorem 4.1.2]{Chewi23}
suggest that the total computational complexity of the algorithm, in terms of the number of evaluations of $\nabla U$,  would be $c_U d \frac{1}{\delta^2}$ times the time complexity, when applied to integrating bounded, Lipschitz continuous functions. 
The restriction to a smaller class of test functions is needed since the Wasserstein error controlled in~\cite[Theorem 4.1.2]{Chewi23} needs to control the integration error. We remark that better error bounds can be obtained by using different discretizations of LMC and by using MALA (see \cite{chewietal2021mala}).

For comparison we note the cost of using rejection sampling to achieve a comparable error is~$C_d / \eta^d$, which is huge when $\eta$ is small and the dimension is large.
Also, the cost of using LMC requires simulating~$1/\delta^2$ realizations of~\eqref{e:Langevin} for a time that is comparable to the~$\delta$-mixing time.
By the Arrhenius law~\cite{Arrhenius89,BovierGayrardEA05} this is~$e^{O(1/\eta)}$, which is much larger than~\eqref{e:CostASMCLangevin} when~$\eta$ is small.

\subsection{Numerical experiments}\label{s:numerics}
In practice, one typically needs to sample  from the target distribution with density proportional to~$e^{-V}$, for some given energy function~$V$.
In situations of interest the energy~$V$ has deep valleys and the associated Gibbs measure has several components (modes). 
To apply ASMC, we choose a temperature~$\eta>0$, which is small enough so that the Gibbs measure with energy function
\begin{equation}
  U \defeq \eta V
  ,
\end{equation}
is easy to sample from.
Then we run Algorithm~\ref{a:ASMC} with the energy function~$U$,
with initial temperature~$\eta_1 = 1$, and final temperature~$\eta$ to deliver samples from the Gibbs measure with density proportional to~$e^{-V}$.
A reference implementation is provided in~\cite{HanIyerEA25code}.

\begin{figure}[htb]
    \centering
    \includegraphics[width=0.5\linewidth]{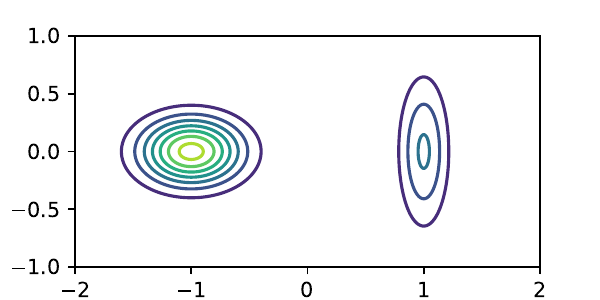}
    \caption{\small%
      Contour plot of the anisotropic Gaussian mixture in $\R^2$, defined in \eqref{e:eqGaussian}, and  used in experiments for Figure \ref{fig:stevN}.
      } 
    \label{fig:Ovaltest}
\end{figure}
\begin{figure}[htb]
    \centering
    \includegraphics[width=0.45\linewidth]{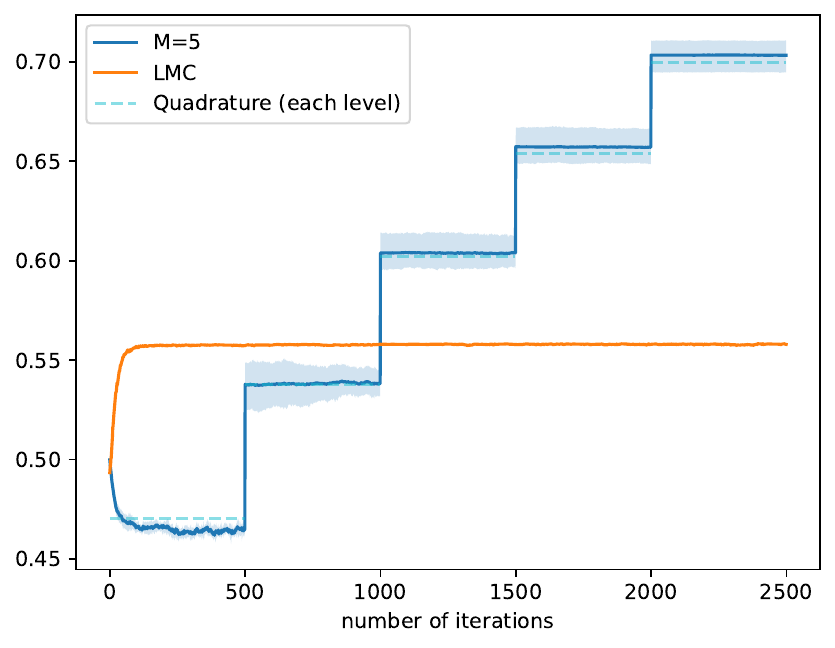}
    \quad
    \includegraphics[width=0.48\linewidth]{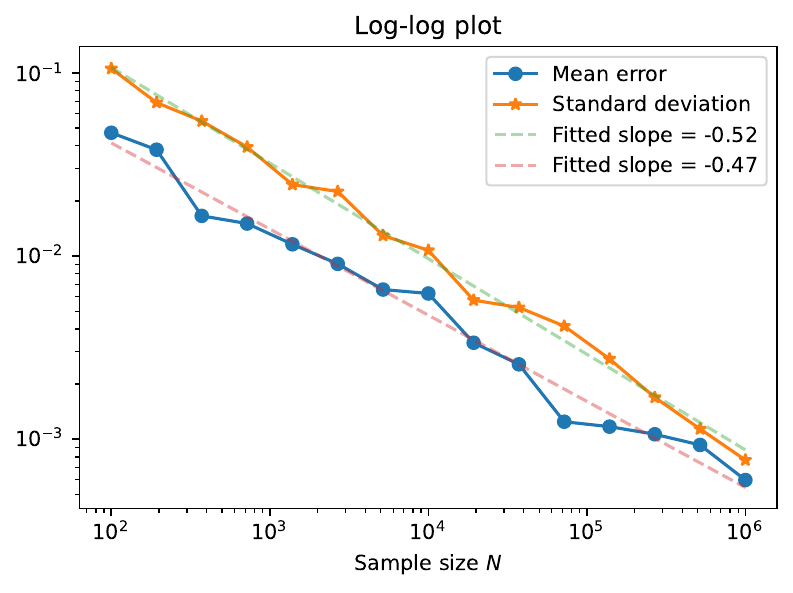}
    \caption{\small
      Left: A Monte Carlo integral computed using ASMC, LMC, and quadrature in 2D. Shaded regions indicate the 25\%-75\% quintile range over 100 independent Monte Carlo runs.
      Right: A Log-log plot of the mean error and standard deviation using ASMC as the number of particles varies.
      }
    \label{fig:stevN}
\end{figure}

For the first illustration, we consider a two-dimensional distribution.
In this case integrals with respect to the Gibbs distribution can also be effectively computed using quadrature, and can be used as a reference for our numerical simulations.
We choose the Gibbs measure to be a mixture of two dimensional, anisotropic Gaussians given by
\begin{equation}\label{e:eqGaussian}
  \pi = \sum_{i = 1}^2 a_i G_{\mu_i, \Sigma_i}
  .
\end{equation}
Here~$G_{\mu, \Sigma}$ is the PDF of the two dimensional Gaussian with mean~$\mu$ and covariance matrix~$\Sigma$.
We choose parameters $a_1=0.7$, $a_2=0.3$, $\mu_1=-e_1$, $\mu_2=e_1$ and
\begin{equation}
  \Sigma_1=\begin{pmatrix}
0.09 & 0\\
0 & 0.04
\end{pmatrix},\quad  \Sigma_2=\begin{pmatrix}
0.02 & 0\\
0 & 0.18
\end{pmatrix}
\end{equation}
A contour plot of~$\pi$ is shown on the left of Figure~\ref{fig:Ovaltest}.
The left panel  of Figure~\ref{fig:stevN} shows the results of numerical simulations computing the Monte Carlo integral of the indicator function of a separating hyperplane using samples from Algorithm~\ref{a:ASMC}.
For comparison, we also show the results of computing the same integral using direct LMC, and using quadrature.
To generate this plot we used $N=10^4$, time step $0.001$, $M=5$, $T=500$. For confirmation, we verify the mean error and standard deviation decrease like~$1/\sqrt{N}$, and show our results in Figure~\ref{fig:stevN}.

\begin{figure}[hbt]
    \centering
    \includegraphics[width=0.47\linewidth]{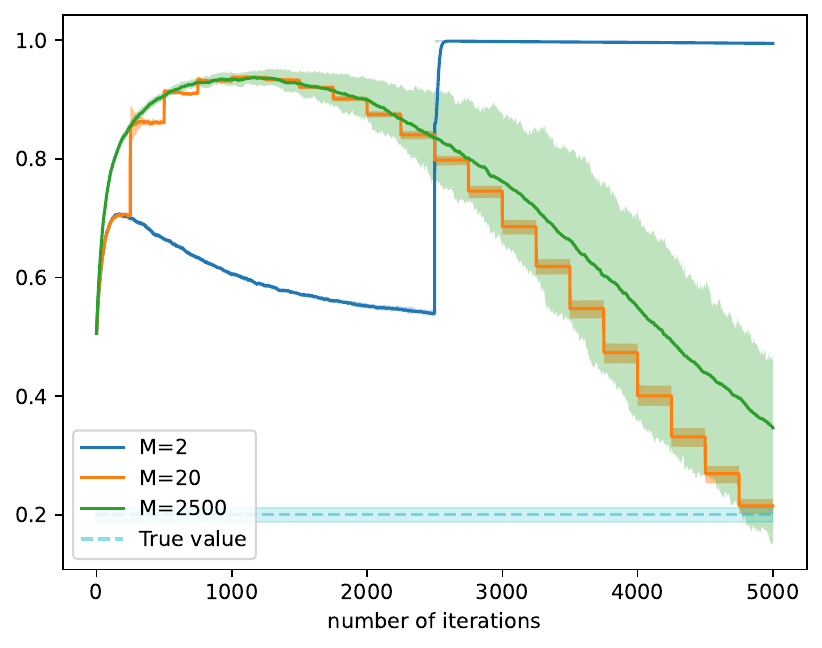}
    \quad
    \raisebox{1ex}{\includegraphics[width=0.48\linewidth]{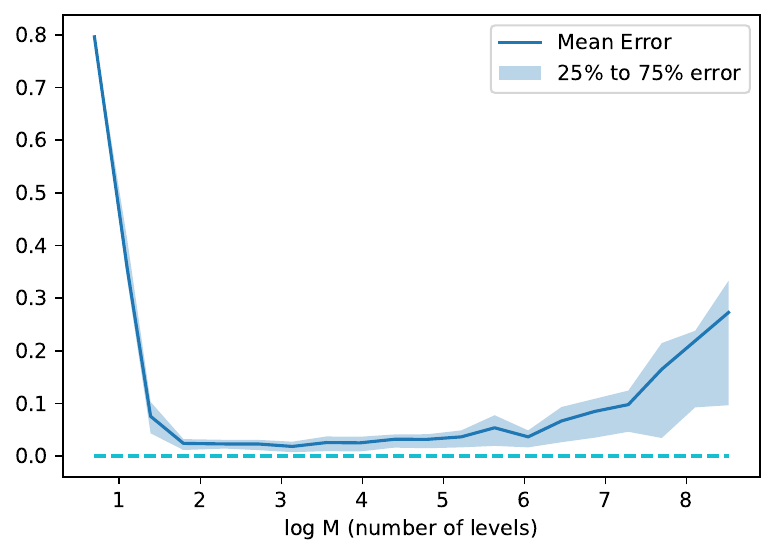}}
    \caption{\small%
    Mean error of an integral in dimension~$d=20$ computed using ASMC as~$M, T$ vary, while holding~$MT$ constant. Shaded regions indicate the 25\%-75\% quintile range.
    Left: A plot of the Monte Carlo integral vs the number of iterations for a few values of~$M$.
    Right: A plot of the mean error vs~$\log M$.
    }
    \label{fig:asMvaries}
\end{figure}

Our next experiment, illustrated in Figure~\ref{fig:asMvaries}, focuses on the trade-off between increasing the number of levels and the number of time steps under a fixed computational budget. 
We use samples obtained by Algorithm~\ref{a:ASMC} to compute a Monte Carlo integral in dimension~$20$. The target measure is a mixture of Gaussians given by \eqref{e:eqGaussian}
with parameters $a_1=0.2$, $a_2=0.8$, $\mu_1=-e_1$, $\mu_2=e_1$ and
\begin{equation}
  \Sigma_1=\frac{1}{16}I_{d},\quad  \Sigma_2=\frac{1}{25}I_{d},\quad d=20.
\end{equation}
We vary the number of levels~$M$ and the level running time~$T$, while keeping the total number of iterations~$M T$ constant.
To generate the plots we used a total of $5000$ iterations per run, sample size $N=10^4$ and time step $0.001$, and~$100$ independent Monte Carlo runs per choice of $M$
 and $T$. We observe that ASMC produces good results for intermediate values of~$M$, but performs poorly when~$M$ is too large or too small when compared to~$d$.
We note that this is not surprising. When $M=1$  ASMC becomes the rejection sampler, and with $M$ being small it is closely approximating a rejection sampler with a few intermediate levels. Since the jumps in temperature are large the resulting bias is large. When $M$ is very large and $T$ is quite small the Markov transitions do not have a chance to mix even within the wells. Thus the procedure basically only involves importance reweighting and resampling, thus  leading to most of the mass concentrated at few nodes, and large error. 

\section{Error Estimates for the Local Mixing Model (Theorem \ref{t:localMixingModel})}\label{s:localModel}

\subsection{Notation and convention}
Before delving into the proof of Theorem~\ref{t:localMixingModel} we briefly list notational conventions that will be used throughout this paper.
\begin{enumerate}[(i)]
  \item
    We will always assume~$C >0$ is a finite constant that can increase from line to line, provided it does not depend on the temperature~$\eta$.

  \item
    We use the convention that the expectation operator~$\E$ has lower precedence than multiplication.
    That is~$\E XY$ denotes the expectation of the product~$\E [XY]$, and~$\E X^2$, denotes the expectation of the square~$\E [X^2]$.

  \item 
    When taking expectations and probabilities, a subscript will denote the conditional expectation / conditional probability.
    That is~$\E_X Y = \E(Y \given X)$ denotes the conditional expectation of~$Y$ given the~$\sigma$-algebra generated by~$X$.
  \item 
    When averaging functions of Markov processes, a superscript will denote the initial distribution.
    That is $\E^\mu f(Y_t)$ denotes~$\E f(Y_t)$ given~$Y_0 \sim \mu$.
    When~$\mu = \delta_y$ is the Dirac~$\delta$-measure supported at~$y$, we will use~$\E^y$ to denote~$\E^{\delta_y}$.

  \item
    We interchangeably use~$\pi_\epsilon$ to denote the measure and the density.
    That is for~$x \in \T^d$, $\pi_\epsilon(x)$ is given by~\eqref{e:piNu}, however for Borel sets~$A$, $\pi_\epsilon(A)$ denotes~$\int_A \pi_\epsilon(x) \, dx$.
\end{enumerate}

\subsection{Description of the Constants in Theorem \ref{t:localMixingModel}}\label{s:constantsLoc}

As remarked earlier, the constants~$C_N$ and~$C_T$ in Theorem~\ref{t:localMixingModel} are explicit.
Since these determine the efficiency of Algorithm~\ref{a:ASMC}, we state them precisely before embarking on the proof of Theorem~\ref{t:localMixingModel}.
First we need auxiliary constant~$C_r = C_r( U/\eta_1, \nu)$ that will be used to bound the ratio of the densities at each level.
For~$k \in \set{1, \dots, M}$, by a slight abuse of notation we define
\begin{equation}\label{e:piKDef}
  \pi_k \defeq \pi_{\eta_k} 
  ,
  \quad
  \tilde \pi_k \defeq \tilde \pi_{\eta_k}
  \quad\text{and}\quad
  Z_k \defeq Z_{\eta_k}
\end{equation}
where~$\pi_{\eta_k}$, $\tilde \pi_{\eta_k}$ and~$Z_{\eta_k}$ are defined by~\eqref{e:piNu} with~$\epsilon = \eta_k$.
Next we define
\begin{equation}\label{e:defrk}
  r_k \defeq \frac{\pi_{k+1}}{\pi_k}
\end{equation}
to be the ratio of \emph{normalized} densities at levels~$k+1$ and~$k$.
In practice, we do not have access to~$r_k$ as we do not have access to the normalization constants~$Z_k$.
This is why Algorithm~\ref{a:ASMC} is formulated using the ratio of unnormalized densities~$\tilde r_k$ defined in~\eqref{e:rkTilde}.
The auxiliary constant~$C_r$ mentioned above is defined by
\begin{equation}\label{eq:defCr}
    C_r \defeq \max_{1\leq k\leq M-1}\|r_k\|_{L^\infty(\mathcal{X})}
    .
\end{equation}

Clearly~$C_r \to 1$ as~$\nu \to 0$.
However, choosing~$\nu$ very small increases the number of levels~$M$ and hence the computational cost of Algorithm~\ref{a:ASMC}.
A bound for~$C_r$, which may be easier to check in practice, is
\begin{equation}\label{eq:CrExplicit}
  C_r
    \leq \inf_{c>0}(1+s_c)\exp\paren{ c\nu }
  \end{equation}
  where
  \begin{equation}\label{def:scU0}
  s_c \defeq
    \frac{\int_{\set{U_0 > c}}e^{-U_0}\, d x}{\int_{\set{U_0\leq c}}e^{-U_0}\, d x}<\infty
    \quad\text{and}\quad  U_0\defeq \frac{U-\inf U}{\eta_1}.
\end{equation}
We prove~\eqref{eq:CrExplicit} in Lemma~\ref{lem: choose eta to make r bound}, below.

Now, the proof of Theorem~\ref{t:localMixingModel} will show that constants~$C_T$ and~$C_N$ are given by
\begin{equation}\label{e:CTCN}
  C_{T}\defeq 4JC_r\paren{2C_{\beta}+1},\quad C_{N}\defeq J^2\paren{2C_{\beta}+1}^2(1+C_r)^2.
\end{equation}
where
\begin{equation}\label{e:CBetaLoc}
  C_{\beta} \defeq \exp(2C_r C_{\LBV}).
\end{equation}

\subsubsection*{Dimensional dependence}

Suppose now $\mathcal{X}=\mathbb{R}^d$.
For a certain class of energy functions, it is possible to make the constants $C_{T}$, $C_{N}$ independent of $d$ by choosing a geometric annealing schedule with~$M$ linear in~$d$.
One such class of energies are those which can be separated into a sum of two functions -- one which depends on the first~$\tilde d$ coordinates and may have multiple local minima; and the other only depends on the last~$d - \tilde d$ coordinates and is convex.

Explicitly, suppose there exists an integer $\tilde{d} \leq d$ such that the function~$U$ is of the form
\begin{equation}\label{e:sepPot}
  U_0(x)=\tilde U_0(x_1,\dots,x_{\tilde{d}})+V_0(x_{\tilde{d}+1},\dots, x_d).
\end{equation}
Here $\tilde{U}_0$ is an any function for which $e^{-\tilde{U}_0}$ is integrable, and may have several local minima.
The function $V_0$ is assumed to be a convex function for which there exist constants $\alpha_0>0, k_0>1$, $\alpha_{u},\alpha_b\in\mathbb{R}$ and a point $x_0\in\mathbb{R}^{d-\tilde{d}}$ such that for all $x\in \mathbb{R}^{d-\tilde{d}}$, we have
\begin{equation}\label{e:AlmostPoly}
\alpha_0|x-x_0|^{k_0}+\alpha_{b}\leq V_0(x)\leq \alpha_0|x-x_0|^{k_0}+\alpha_{u}.
\end{equation}
One class of functions that have this structure are Gaussian mixtures of points that are located on a~$\tilde d$ dimensional hyperplane, and whose covariance matrices in the perpendicular direction are all equal.
For such energies we have the following dimension independent bounds.
\begin{proposition}\label{prop:CrEg}
Assume that~$U_0$ satisfies and \eqref{e:sepPot} and \eqref{e:AlmostPoly}.   Choose 
\begin{equation}
  M \geq \ceil*{ \frac{d}{ \eta} }
\end{equation}
and $\eta_k$ such that $1/\eta_1$, \dots, $1/\eta_M$ are linearly spaced.
  Then $C_{T}$ and $C_{N}$ in \eqref{e:CTCN} can be bounded above in terms of~$\alpha_0$, $\alpha_b$, $\alpha_u$, $k_0$, $U_0$, but \emph{independent} of $d$.  
\end{proposition}

Proposition~\ref{prop:CrEg} will be proved using asymptotics for the incomplete gamma function and the proof is Section~\ref{s:CBV}, below.

\subsection{Proof of Theorem \ref{t:localMixingModel}}\label{s:locMMProof}
In order to prove Theorem~\ref{t:localMixingModel}, we note that Algorithm~\ref{a:ASMC} consists of repeating two steps: (local) exploration using using the process~$Y_{\epsilon, \cdot}$ (Algorithm~\ref{a:ASMC}, step~\ref{i:exploreValley}), and then resampling (Algorithm~\ref{a:ASMC}, step~\ref{i:resampling}).
We now state lemmas for the errors accumulated in each of these steps.

To quantify the Monte Carlo error made by running the process~$Y_{\epsilon, \cdot}$ in the (local) exploration step, we introduce the following notation.
Given~$\epsilon, t > 0$ and a bounded test function~$h$, define the Monte Carlo error~$\Err_{\epsilon, t}(h)$ by
\begin{equation}\label{eq:defErr}
   \Err_{\epsilon,t}(h)\defeq
    \norm[\Big]{\frac{1}{N}\sum_{i=1}^{N}h(Y^{i}_{\epsilon, t})-\int_{\mathbb T^d}h\pi_{\epsilon}\,d x}_{L^2(\P)}
   ,
\end{equation}
where~$Y^{i}_{\epsilon, \cdot}$ are~$N$ independent realizations of a Markov process with transition density~\eqref{e:YTransitionKernel}.
\begin{lemma}\label{lem:mixNeedMass}
  Given~$N$ (random) points~$y^1$, \dots, $y^N$, let~$Y^i_{\epsilon, \cdot}$ be~$N$ independent realizations of the Markov process with transition density~\eqref{e:YTransitionKernel} and initial distribution~$Y^i_{\epsilon, 0} = y^i$.
  Then for any bounded test function~$h$, and any~$T \in \N$ we have
\begin{equation}\label{eq:mc}
     \Err_{\epsilon,T}(h)
      \leq  \chi^{T}_{\epsilon}\norm[\bigg]{
	\sum_{j=1}^{J}
	  \paren[\Big]{
	    1 - \frac{\mu_0(\Omega_j)}{\pi_\epsilon(\Omega_j)}
	  }
	  \int_{\Omega_j} h\pi_{\epsilon}\, dx
	}_{L^2(\P)} +\frac{\norm{h}_{\osc}}{2\sqrt{N}}
    ,
  \end{equation}
  where~$\mu_0$ is the empirical measure
  \begin{equation}\label{e:mu0}
    \mu_0 = \frac{1}{N} \sum_{i = 1}^N \delta_{y^i}
    .
  \end{equation}
  Consequently,
\begin{equation}\label{e:ErrTErr0}
    \Err_{\epsilon,T}(h)
      \leq  \frac{\norm{h}_\osc}{2} \paren[\Big]{
	\chi^{T}_{\epsilon}
	  \sum_{j=1}^{J}
	    \Err_{\epsilon, 0}(\one_{\Omega_j})
	   +\frac{1}{\sqrt{N}}
	}
    .
  \end{equation}
\end{lemma}

We clarify that~$\mu_0(\Omega_j)$ is random as the initial points~$y^{i}$ are themselves random.
The second term on the right of~\eqref{e:ErrTErr0} is the standard Monte Carlo error which can be made small by making~$N$ large.
To make the first term small, we have two options:
The first option is to wait for the mixing time of~$Y_{\epsilon, \cdot}$, and obtain smallness from the~$\chi_\epsilon^T$ factor.
The second is to ensure~$\sum_{j} \Err_{\epsilon , 0}(\one_{\Omega_j})$ is small.
In our situation the first option is undesirable as it requires~$T \gg 1 / |\ln \chi_\epsilon|$, which is too large to be practical.
Instead we use the second option, and make~$\sum_{j} \Err_{\epsilon , 0}(\one_{\Omega_j})$ small by ensuring the fraction of initial points in each domain~$\Omega_j$ is close to~$\pi_\epsilon(\Omega_j)$.
\smallskip

We now turn to the resampling step.
Suppose we have~$N$ i.i.d.\ samples~$x^1$, \dots, $x^N$ from a distribution with an unnormalized probability density function~$\tilde p$.
Let~$\tilde q$ be another unnormalized probability density function, such that $\set{\tilde q >0} \subseteq \set{p > 0}$.
  Choose~$(y^1, \dots, y^N)$ to be a \emph{resampling} of the points $\paren{x^1, \dots, x^N}$ using the multinomial distribution with probability
\begin{equation}\label{e:PyiEqXj}
  \P( y^i = x^j ) = \frac{\tilde r(x^j)}{\sum_{i=1}^N \tilde r(x^i)},
  \quad\text{where}\quad
  \tilde r \defeq \frac{\tilde q}{\tilde p}
  .
\end{equation}

Of course, some of the points~$x^{i}$ may be chosen multiple times and the points~$y^{1}$, \dots, $y^N$ may not be distinct.
Nevertheless, a simple heuristic argument suggests that when~$N$ is large the distribution of each of the points~$y^{i}$ will have a density proportional to~$\tilde q$.
Indeed, suppose~$\mathcal X$ is finite, $N \gg \abs{\mathcal X}$ and~$p, q$ are the normalized probability distributions corresponding to~$p, q$ respectively.
Then each~$x \in \mathcal X$ occurs amongst the points~$\set{x^1, \dots, x^N}$ roughly $N p(x)$ times, and so
\begin{equation}
  \P( y^i = x )
    \approx \frac{\tilde r(x) N p(x)}{ \sum_{x' \in \mathcal X} \tilde r(x') N p(x') }
    = \frac{\tilde q(x)}{ \sum_{x' \in \mathcal X} \tilde q(x') }
    \approx q(x)
    \,.
\end{equation}

To make the above quantitative, and usable in our context, some care has to be taken.
The points~$y^{i}$ are only conditionally independent given~$x^1$, \dots, $x^N$; they are not unconditionally independent, and it is hard to estimate the unconditional joint distribution.
We will instead obtain a Monte Carlo estimate which both quantifies the error and is sufficient for our purposes.

\begin{lemma}\label{l:resampling}
  Suppose~$x^1$, \dots, $x^N$ are~$N$ (not necessarily i.i.d.) random points in~$\mathcal X$.
  Let~$\tilde p$, $\tilde q \colon \mathcal X \to [0, \infty)$ be two unnormalized probability density functions, and choose $y^1,\dots,y^N$ independently from~$\set{x^1, \dots, x^N}$ according to~\eqref{e:PyiEqXj}.
  Then for any test function $h \in L^\infty(\mathcal X)$,  we have
  \begin{align}
    \MoveEqLeft[10]
      \norm[\Big]{ \frac{1}{N} \sum_{i=1}^N h(y^{i}) - \int_{\mathcal X} h q \, dx }_{L^2(\P)}
      \leq
	\frac{1}{\sqrt{N}} \norm[\Big]{h-\int_{\mathcal X} h q \, dx}_{L^\infty}
    \\
	&+ \norm[\Big]{h-\int_{\mathcal X} h q\, dx}_{L^\infty} \norm[\Big]{1- \frac{1}{N}\sum_{i=1}^N r(x^i)}_{L^2(\P)}
      \\
      \label{e:resampling}
	&+ \norm[\Big]{ \frac{1}{N}\sum_{i=1}^N r(x^i)\Big(h(x^i)  -\int_{\mathcal X} h q \, dx \Big)}_{L^2(\P)}
      .
  \end{align}
  Here~$r$ is the ratio
  \begin{equation}\label{e:rdef}
    r \defeq \frac{q}{p}
    ,
    \quad\text{where}\quad
    p = \frac{\tilde p}{\int_{\mathcal X} \tilde p \, dx}
    \quad\text{and}\quad
    q = \frac{\tilde q}{\int_{\mathcal X} \tilde q \, dx}
    .
  \end{equation}
\end{lemma}

Note Lemma~\ref{l:resampling} does not assume~$x^1$, \dots, $x^N$ are independent, or even that they have distribution~$p$.
If, however, the points~$x^1$, \dots, $x^N$ give good Monte Carlo estimates for integrals with respect to~$p$, then the right hand side of~\eqref{e:resampling} will be small.
Explicitly, in the typical situation where $x^i\sim p$ are i.i.d., we will have
\begin{equation}
  \norm[\Big]{\frac{1}{N} \sum_{i=1}^N g(x^i) - \int_{\mathcal X} g p \, dx}_{L^2(\P)}^2
    \leq \frac{C \var(g)}{N}
    .
\end{equation}
for any bounded test function~$g$.
Combined with the fact that
\begin{equation}
  \int_{\mathcal X} r p \, dx = \int_{\mathcal X} q \, dx = 1
  \quad\text{and}\quad
  \int_{\mathcal X} h r p \, dx = \int_{\mathcal X} h q \, dx
  ,
\end{equation}
this shows the right hand side of~\eqref{e:resampling} is~$O(1/\sqrt{N})$.
\medskip

We now use Lemma \ref{lem:mixNeedMass} and Lemma \ref{l:resampling} to derive a recurrence relation for the Monte Carlo error between levels~$k$ and~$k+1$ in Algorithm~\ref{a:ASMC}.

\begin{lemma}\label{lem:iterativelocal}
For each  $k=1,\dots,M-1$,
    \begin{align}
      \max_{1 \leq \ell \leq J} \Err_{k+1,0}(\one_{\Omega_{\ell}})&\leq \frac{1+\norm{r_k}_{\osc}}{\sqrt{N}}
   \\\label{eq:iterativelocal}
      &\quad + \paren[\Big]{1+2\sum_{j=1}^{J}\abs[\Big]{\frac{\pi_{k+1}(\Omega_j)}{\pi_{k}(\Omega_j)}-1}} \cdot\max_{1 \leq \ell \leq J} \Err_{k,0}(\one_{\Omega_{\ell}}).
    \end{align}
\end{lemma}
Here~$r_k$ is the ratio of the normalized densities defined in~\eqref{e:defrk}, and by a slight abuse of notation we use~$\Err_{k, \cdot}$ to denote~$\Err_{\eta_k, \cdot}$.

The proof of Theorem~\ref{t:localMixingModel} now reduces to solving the recurrence relation~\eqref{eq:iterativelocal} and using Lemma~\ref{lem:mixNeedMass}.
Notice that~\eqref{eq:iterativelocal} involves a bound on~$\norm{r_k}_\osc$, and the maximum of this as~$k$ varies is precisely the constant~$C_r$ defined in~\eqref{eq:defCr}.
A bound on~$C_r$ that may be easier to obtain in practice is~\eqref{eq:CrExplicit}, which we prove in Lemma~\ref{lem: choose eta to make r bound}, below.
Momentarily postponing the proofs of the above lemmas, we prove Theorem~\ref{t:localMixingModel}.
\begin{proof}[Proof of Theorem \ref{t:localMixingModel}]
  Applying \eqref{e:ErrTErr0} with $\epsilon=\eta_{M}$ gives
\begin{equation}\label{eq:mc2}
    \Err_{M,T}(h)\leq  \frac{\norm{h}_{\osc}}{2}\paren[\Big]{J\max_{j=1,\dots,J}\Err_{M,0}(\one_{\Omega_j})+\frac{1}{\sqrt{N}}}.
  \end{equation}

 We will show that the right hand side of~\eqref{eq:mc2} is bounded above by~$\delta \norm{h}_\osc$. For the first term, a direct calculation using~\eqref{eq:iterativelocal} immediately shows that
\begin{align} 
  \max_{j=1,\dots,J}\Err_{M,0}(\one_{\Omega_j}) &\leq \paren[\Big]{\prod_{\ell=2}^{M-1}\Theta(\ell,\ell+1)}\max_{j=1,\dots,J}\Err_{2,0}(\one_{\Omega_j})
       \\
       &\qquad+\sum_{k=2}^{M-1}\frac{1+\norm{r_k}_{\osc}}{\sqrt{N}} \prod_{\ell=k+1}^{M-1}\Theta(\ell,\ell+1)
       \label{eq:maxerr1}
       \end{align}
       where
       \begin{equation}
           \label{eq:defThetakloc}
           \Theta(\ell,\ell+1)\defeq 1+2\sum_{j=1}^{J}\abs[\Big]{\frac{\pi_{\ell+1}(\Omega_j)}{\pi_{\ell}(\Omega_j)}-1}.
       \end{equation}
       To finish the proof, we now need to estimate the terms $\prod_{\ell=k}^{M-1}\Theta(\ell,\ell+1)$ and $\max_{1 \leq j \leq J}\Err_{2, 0} (\one_{\Omega_j})$.
 
  \restartsteps
  \step[Estimating $\prod_{\ell=k}^{M-1}\Theta(\ell,\ell+1)$]
  Notice that for every $j=1,\dots, J$, and every~$k=1,\dots, M$, we have
  \begin{equation}\label{eq:massrationboundedbyCr}
      0<\frac{\pi_{k+1}(\Omega_j)}{\pi_{k}(\Omega_j)}\leq \norm{r_{k}}_{L^\infty}\overset{\eqref{eq:defCr}}{\leq} C_r.
  \end{equation}
  Using the fact that
  \begin{equation}\label{e:ylogy}
    \abs{y - 1} \leq (1 \varmax y) \abs{\ln y}
    ,
  \end{equation}
  for any $k= 1,\dots, M-1$, we obtain
\begin{align}\label{eq:leqCbetaloc}
  \MoveEqLeft
  \prod_{\ell=k}^{M-1}\Theta(\ell,\ell+1)
    \overset{\text{AM-GM}}{\leq} \paren[\bigg]{1+  \frac{2}{M-k}\sum_{j=1}^{J}\sum_{\ell=k}^{M-1}\abs[\Big]{\frac{\pi_{\ell+1}(\Omega_j)}{\pi_{\ell}(\Omega_j)}-1}}^{M-k}
    \\
    &\overset{\mathclap{\eqref{eq:massrationboundedbyCr}, \eqref{e:ylogy}}}{\leq} \quad\; \paren[\bigg]{1+ \frac{2}{M-k}\sum_{j=1}^{J}\sum_{\ell=k}^{M-1}C_r\abs[\Big]{\log\paren[\Big]{\frac{\pi_{\ell+1}(\Omega_j)}{\pi_{\ell}(\Omega_j)}}}}^{M-k}\\
    &=\paren[\bigg]{1+ \frac{2}{M-k}\sum_{j=1}^{J}\sum_{\ell=k}^{M-1}C_r\abs[\Big]{\int_{\eta_{\ell+1}}^{\eta_{\ell}} \partial_\epsilon \ln \pi_\epsilon( \Omega_j ) \, d\epsilon }}^{M-k}\\
  &\overset{\mathclap{\eqref{e:CBV}}}{\leq} \:\exp(2C_r C_{\LBV})\overset{\eqref{e:CBetaLoc}}{=} C_{\beta}.
\end{align}

  \step[Estimating~$\Err_{2,0}(\one_{\Omega_j})$]
  Applying Lemma \ref{l:resampling} with
$p=\pi_1$,
  $q=\pi_2$,
  $h=\one_{\Omega_j}$,
  and~$x^{i}=Y_{1,T}^{i}$,
to obtain
  \begin{equation}\label{e:err20}
    \Err_{2,0}(\one_{\Omega_j})\leq \frac{1}{\sqrt{N}}+\Err_{1,T}(r_1)+\Err_{1,T}(r_1(\one_{\Omega_j}-\pi_2(\Omega_j))).
\end{equation}

  Now the processes $Y^{1}_{1,\cdot}$, \dots, $Y^N_{1, \cdot}$ are all independent.\footnote{
    For~$k \geq 2$ the processes~$Y^{1}_{k, \cdot}$, \dots, $Y^N_{k, \cdot}$ are no longer independent as the initial distributions are not independent.
  }
  Thus for any bounded test function $h$,
\begin{align}
  \MoveEqLeft
  (\Err_{1,T}(h))^2 =
    \E \paren[\Big]{
    \frac{1}{N} \sum_1^N h(Y_{1,T}^i)
      - \int_{\mathcal X} h \, \pi_1 \, dx
      }^2\\
    &=
    \E \paren[\bigg]{
     \frac{1}{N}\sum_1^N \paren[\big]{
       h(Y_{1,T}^i) - \E h(Y_{1,T}^i)
      }
      + \frac{1}{N} \sum_1^N 
       \E h(Y_{1,T}^i)
       -\int_{\mathcal X} h \, \pi_1 \, dx
      }^2
    \\
  &\leq \frac{1}{N}\norm{h}^2_{L^{\infty}}+\paren[\Big]{\frac{1}{N}\sum_{i=1}^{N}\E 
 h(Y_{1,T}^i)
    - \int_{\mathcal X} h \, \pi_1 \, dx
    }^2\\
    &\leq \frac{1}{N}\norm{h}^2_{L^{\infty}}+\norm{h}^2_{L^{\infty}} \paren[\Big]{\frac{1}{N}\sum_{i=1}^{N}\norm{p^1_{T}(y_1^i, \cdot) -\pi_1}_{L^1}}^2
    ,
\end{align}
and hence
\begin{equation}\label{eq:err1Th}
   \Err_{1,T}(h)\leq \frac{1}{\sqrt{N}}\norm{h}_{L^{\infty}}+\norm{h}_{L^{\infty}} \frac{1}{N}\sum_{i=1}^{N}\norm{p^1_{T}(y_1^i, \cdot) -\pi_1}_{L^1}
    .
\end{equation}

 Notice that the choice of $C_T$ and $C_N$ in~\eqref{e:CTCN}, implies
\begin{equation}\label{e:TNloc}
T\geq \tmix\paren[\Big]{ Y_{\eta_1, \cdot}, \frac{\tilde{\delta}}{4C_r}},
  \quad\text{and}\quad
  \frac{1}{\sqrt{N}}\leq \frac{1+C_r}{\sqrt{N}}
      \overset{\eqref{eq: parametersTNloc}}{\leq}\frac{\tilde{\delta}}{M}
\end{equation}
where 
\begin{equation}\label{e:TildeDeltaloc}
      \tilde{\delta}=
    \frac{\delta}{ J\paren{2C_{\beta}+1}}.
  \end{equation}
Using \eqref{e:TNloc} in~\eqref{eq:err1Th} with  $h = r_1$ gives
\begin{equation}\label{eq:err1r1}
    \Err_{1,T}(r_1)\leq \frac{1}{\sqrt{N}}\norm{r_1}_{L^{\infty}}+\frac{\tilde{\delta}}{2C_r}\norm{r_1}_{L^{\infty}}\overset{\eqref{eq:defCr},\eqref{e:TNloc}}{\leq} \frac{\tilde{\delta}}{M}+\frac{\tilde{\delta}}{2}.
\end{equation}
Similarly,
\begin{equation}\label{eq:err1r1omega}
    \Err_{1,T}(r_1(\one_{\Omega_j}-\pi_2(\Omega_j)))\leq\frac{\tilde{\delta}}{M}+\frac{\tilde{\delta}}{2}.
\end{equation}
  Plugging \eqref{eq:err1r1},\eqref{eq:err1r1omega} and \eqref{e:TNloc} into~\eqref{e:err20} yields
  \begin{equation}\label{e:Err2loc}
      \max_{j=1,\dots,J}\Err_{2,0}(\one_{\Omega_j})\leq \frac{3\tilde{\delta}}{M}+\tilde{\delta}.
  \end{equation}

 Now using~\eqref{eq:maxerr1}, we obtain
\begin{align} 
\max_{j=1,\dots,J}\Err_{M,0}(\one_{\Omega_j})
  \quad~~ &\overset{\mathclap{\eqref{eq:leqCbetaloc},~\eqref{e:Err2loc}}}{\leq} \quad~~ C_{\beta}\paren[\Big]{\frac{3\tilde{\delta}}{M}+\tilde{\delta}}+\sum_{k=2}^{M-2}C_{\beta}\frac{\tilde{\delta}}{M} +\frac{\tilde{\delta}}{M}
    \\
    &= \paren[\Big]{
      2 C_\beta + \frac{1}{M}
    } \tilde \delta
    \overset{\eqref{e:TildeDeltaloc}}{\leq}\frac{\delta}{J}
       . \label{eq: estmaxm}
\end{align}

Using~\eqref{e:TNloc} and~\eqref{eq: estmaxm} in~\eqref{eq:mc2} implies
\begin{equation}
  \Err_{M, T}(h) \leq \frac{\|h\|_{\osc}}{2}\paren[\Big]{\frac{J \delta}{J}+\frac{\tilde{\delta}}{M}}\overset{\eqref{e:TildeDeltaloc}}{<}\delta\|h\|_{\osc}
  .
\end{equation}
This proves~\eqref{e:MCErrorModel}, concluding the proof.
\end{proof}

\section{Error Estimates for a double-well energy (Theorem~\ref{thm: main})}\label{s:mainProof}

The aim of this section is to prove Theorem~\ref{thm: main} and obtain error estimates when ASMC is used to sample from a double-well energy function on a~$d$-dimensional torus.

\subsection{Assumptions and Notation}\label{s:assumptions}
We begin by precisely stating the assumptions that were used in Theorem~\ref{thm: main}.
The first assumption requires~$U$ to be a regular, double-well function with nondegenerate critical points.

\begin{assumption}\label{a:criticalpts}
  The function~$U\in C^6(\mathbb T^d,\mathbb R)$, has a nondegenerate Hessian at all critical points, and has exactly two local minima located at~$x_{\min, 1}$ and~$x_{\min, 2}$.
  We normalize~$U$ so that
\begin{equation}\label{e:Upositive}
  0 = U(x_{\min, 1} )  \leq U( x_{\min, 2} )
  .
\end{equation}
\end{assumption}

Our next assumption concerns the saddle between the local minima~$x_{\min, 1}$ and~$x_{\min, 2}$.
Define the saddle height between~$x_{\min,1}$ and $x_{\min,2}$ to be the minimum amount of energy needed to go from the global minimum~$x_{\min, 1}$ to~$x_{\min, 2}$, and is given by
\begin{equation}\label{e:UHatDef}
    \hat{U} = \hat{U}(x_{\min,1},x_{\min,2}) \defeq \inf_\omega \sup\limits_{t\in [0,1]}U(\omega(t))
    .
\end{equation}
Here the infimum above is taken over all continuous paths~$\omega \in C( [0, 1]; \mathbb T^d)$ such that~$\omega(0) =x_{\min,1}$, $\omega(1) = x_{\min,2}$.
To prove Theorem~\ref{thm: main} we need to assume a nondegeneracy condition on the saddle.

\begin{assumption}\label{assumption: nondegeneracy}
   The saddle height between $x_{\min,1}$ and $x_{\min,2}$ is attained at a unique critical point $s_{1,2}$ of index one.
   That is, 
   the first eigenvalue of $\Hess U(s_{1,2})$ is negative and the others are positive.
\end{assumption}

We can now define the ratio~$\hat \gamma_r$ that appeared in~\eqref{e:MTN}, above.
The \emph{energy barrier}, denoted by~$\hat{\gamma}$, is defined to be the minimum amount of energy needed to go from the (possibly local) minimum~$x_{\min, 2}$ to the global minimum~$x_{\min, 1}$.
In terms of~$s_{1, 2}$, the energy barrier~$\hat \gamma$ and the saddle height are given by
\begin{equation}\label{e:gammaHatDef}
    \hat{\gamma}\defeq U(s_{1,2})-U(x_{\min,2}),
    \quad\text{and}\quad
    \hat U = U(s_{1,2})
    .
\end{equation}
The ratio $\hat{\gamma}_r$ is the ratio of the saddle height~$\hat U$ to the energy barrier~$\hat \gamma$, given by
\begin{equation}\label{e:gammaHatRDef}
    \hat \gamma_r\defeq \frac{\hat{U}}{\hat{\gamma}}
    .
\end{equation}

Finally, we require the distribution~$\pi_\eta$ to be truly multimodal in the temperature range of interest.
That is, we require the mass in the \emph{basins of attraction} around each of the local minima~$x_{\min, 1}$ and~$x_{\min, 2}$ to be bounded away from~$0$.
We recall the basin of attraction around~$x_{\min, i}$, denoted by~$\Omega_i$, is the set of all initial points for which the gradient flow of~$U$ eventually reaches~$x_{\min, i}$.
Precisely, $\Omega_{i}$ is defined by
\begin{equation}
  \Omega_i
    \defeq
    \set[\Big]{
      y\in\mathbb T^d \st
      \lim_{t\to\infty}y_t=x_{\min,i}, \text{ where }\dot{y}_t=-\nabla U(y_t)
	\text{ with } y_0=y
      }
    ,
\end{equation}
and our multimodality condition is as follows.

\begin{assumption}\label{a:massRatioBound}
  There exists~$0 \leq \eta_{\min} < \eta_{\max} \leq \infty$, a constant~$C_m$ such that
   \begin{equation}\label{e:massRatioBound}
     \inf_{\substack{
       \epsilon\in [\eta_{\min},\eta_{\max}]\\
       0 < \epsilon < \infty
     }}
      \pi_{\epsilon}(\Omega_i)
	\geq \frac{1}{C_m^2}
	.
   \end{equation}
\end{assumption}

  We will show (Lemma~\ref{lem: upper bound on mass ratio of two wells}, below) that~\eqref{e:massRatioBound} is satisfied if the wells have \emph{nearly equal depth}.
  That is, if~$U(x_{\min, 2}) - U(x_{\min, 1}) \leq O(\eta_{\min})$, then one can show~\eqref{e:massRatioBound} holds for some constant~$C_m$ that is independent of~$\eta_{\min}$.
  We state this precisely as the following lemma.

\begin{lemma}\label{lem: upper bound on mass ratio of two wells}
  Suppose $U$ satisfies Assumptions~\ref{a:criticalpts}, \ref{assumption: nondegeneracy}, and there exists a temperature~$\eta_{\min} \geq 0$ and constant $C_{\ell}>0$ such that
\begin{equation}\label{eq:depthdifference}
  U(x_{\min, 2}) - U(x_{\min, 1} )\leq C_{\ell}\eta_{\min}
  .
\end{equation}
  Then for any finite~$\eta_{\max} > \eta_{\min}$ there exists a constant~$C_m = C_m( U, \eta_{\max}, C_\ell)$, independent of~$\eta_{\min}$ such that~\eqref{e:massRatioBound} holds.
\end{lemma}

\begin{remark}
  We note that the condition~\eqref{eq:depthdifference} implies the finiteness condition~\eqref{e:CBV} that was used in Theorem~\ref{t:localMixingModel}.
  This is shown in Corollary~\ref{cor:eqcbvholds}, below, and was previously referred to in Remark~\ref{r:etaMinIndep}.
\end{remark}

\subsection{Proof of Theorem \ref{thm: main}}

In this subsection, we explain the main idea behind the proof of Theorem \ref{thm: main}.
For simplicity and without loss of generality we assume $\eta_1=1$.
We begin by rewriting our algorithm in a manner that that is convenient for the proof.
Fix~$T > 0$ and~$N \in \N$ that will be chosen later.

\restartsteps
\step
  We start with~$N$ arbitrary points~$y^1_1$, \dots, $y^N_1$.

\step[Langevin step]
  For each~$k \in \set{1, \dots, M}$, and~$i \in \set{1, \dots N}$, let~$X^i_{k, \cdot}$ be the solution to the overdamped Langevin equation~\eqref{e:Langevin} with initial data~$X^i_{k, 0} = y^i_k$, driven by independent Brownian motions.

\step[Resampling step]
  Given the processes~$\set{X^i_{k, \cdot} \st i \leq N, k \leq M-1}$ we choose the points~$\set{y^1_{k+1}, \dots, y^N_{k+1}}$ independently from~$\set{X^1_{k, T}, \dots, X^N_{k, T}}$ so that
  \begin{equation}
    \P( y^i_{k+1} = X^j_{k, T} ) = \frac{\tilde r_k(X^j_{k, T})}{\sum_{i=1}^N \tilde r_k(X^i_{k,T})}
    .
  \end{equation}
  Here~$\tilde r_k$ is the ratio defined by~\eqref{e:rkTilde}.

  \subsubsection{Spectral Properties}
  We now briefly recall a few standard facts about the overdamped Langevin dynamics~\eqref{e:Langevin} that will be used in the proof. Let~$L_{\epsilon}$ be the generator of~\eqref{e:Langevin}, whose action on smooth test functions is defined by
\begin{equation}\label{e:Ldef}
  L_{\epsilon} f \defeq-\epsilon\Delta f +\grad U\cdot\grad f
  .
\end{equation}
  Let~$L_\epsilon^*$ be the dual operator defined by
  \begin{equation}\label{e:LStarDef}
    L_\epsilon^* f = -\dv( \grad U f ) - \epsilon \lap f
    \,.
  \end{equation}
  It is well known~\cite[Chapter 8]{Oksendal03} that if~$Y_{\epsilon, \cdot}$ solves~\eqref{e:Langevin} then its density~$f_t \defeq \PDF(Y_{\epsilon, t})$ satisfies the \emph{Fokker-Planck equation}, a.k.a. the \emph{Kolmogorov forward equation}
  \begin{equation}\label{e:KForward}
    \partial_t f + L_\epsilon^* f = 0
    .
  \end{equation}
  One can readily check that the Gibbs distribution~$\pi_\epsilon$ is a stationary solution of~\eqref{e:KForward}, and hence must be the stationary distribution of~\eqref{e:Langevin}.
A direct calculation shows that
\begin{equation}\label{e:fOverPiNu}
  \partial_t \paren[\Big]{\frac{f}{\pi_\epsilon}}
    + L_\epsilon \paren[\Big]{\frac{f}{\pi_\epsilon}}
    = 0
    .
\end{equation}

The mixing properties of Langevin dynamics can be deduced directly from the spectral properties of the operator~$L_\epsilon$, as we now explain.
It is well known (see for instance~\cite[Chapter 8]{Kolokoltsov00}) that on the weighted space~$L^2(\pi_{\epsilon})$ the operator~$L_\epsilon$ is self-adjoint and has a discrete spectrum with eigenvalues
\begin{equation}
  0 = \lambda_{1, \epsilon}
    < \lambda_{2, \epsilon}
    \leq \lambda_{3, \epsilon}
    \cdots
\end{equation}
with corresponding~$L^2(\pi_\epsilon)$ normalized eigenfunctions~$\psi_{1, \epsilon}$, $\psi_{2, \epsilon}$, etc.
The first eigenvalue~$\lambda_{1, \epsilon} = 0$ corresponds to the constant eigenfunction~$\psi_{1, \epsilon} \equiv 1$.
The proof of Theorem~\ref{thm: main} crucially relies on certain spectral properties of the operator~$L_\epsilon$, which we now list.

\begin{property}[Eigenvalue bounds]\label{p:spectral}
  For every~$H > \hat U$, there exists a constant~$C_H$ such that for every~$\epsilon \in (0, 1]$ we have
\begin{equation}\label{e:egvalLEpsilon}
  \lambda_{2,\epsilon}
    \geq C_H \exp\paren[\Big]{-\frac{H}{\epsilon}}
  .
\end{equation}
  Also, there exists~$\Lambda$ such that for all~$\epsilon \in (0, 1]$ such that
  \begin{equation}\label{e:lambdaiLower}
    \lambda_{i,\epsilon}\geq \Lambda,\quad 
    \text{for all } i\geq 3
    .
  \end{equation}
\end{property}

We recall the function~$g \colon (0, 1] \to [0, \infty)$ has bounded variation if~$V_{(0, 1]}(g) < \infty$.
Here~$V_{(0, 1])}(g)$ is the first variation of~$g$, defined by
\begin{equation}
  V_{(0, 1]}(g) \defeq
    \sup\set[\Big]{ \sum_{i=1}^k |g(x_i) - g(x_{i-1})| \st 0<x_0 \leq \cdots \leq x_k \leq 1}
  .
\end{equation}

\begin{property}[Eigenfunction variation]\label{p:var}
  There exists a function~$g \colon (0, 1] \to [0, \infty)$ of bounded variation, such that the following holds.
  For every~$\gamma < \hat \gamma$, there exist constants~$C_1=C_1(\gamma)$ and~$C_2=C_2(\gamma)$ such that for every~$0 < \epsilon' < \epsilon \leq 1$ we have
  \begin{gather}
    \label{e:PsiEpPi}
    \abs[\bigg]{\int_{\mathbb T^d}\psi_{2,\epsilon}\pi_{\epsilon'}\, d x }
      \leq 
    C_1\paren[\Big]{\exp\paren[\Big]{-\frac{\gamma}{\epsilon}}+\abs[\big]{g(\epsilon')-g(\epsilon)}}
  \\
    \label{e:PsiEpEpP}
    \abs[\bigg]{
      \int_{\mathbb T^d}\psi_{2,\epsilon}\psi_{2,\epsilon'}\pi_{\epsilon'}\, d x
    }
    \leq 1+C_2\paren[\Big]{\exp\paren[\Big]{-\frac{\gamma}{\epsilon}}+\abs[\big]{g(\epsilon')-g(\epsilon)}}
    .
  \end{gather}
\end{property}

\begin{property}[Eigenfunction bounds]\label{property:eigenfunctions}
  There exists constant~$C_{\psi}$, independent of~$\epsilon$ such that 
  \begin{equation}\label{eq:defCpsi}
    \sup\limits_{0< \epsilon\leq 1}\|\psi_{2,\epsilon}\|_{L^{\infty}(\mathbb{T}^d)}\leq C_{\psi}.
  \end{equation}
\end{property}

\begin{remark}[Dimensional dependence of the constants in Theorem~\ref{thm: main}]
  The dimensional dependence of the constants~$C_N$ and~$C_T$ in Theorem~\ref{thm: main} arises mainly through the constants~$\Lambda$, $C_H$, $C_1$, $C_2$ and~$C_\psi$ in Properties~\ref{p:spectral}--\ref{property:eigenfunctions}.
  Given these constants, $C_N$ and~$C_T$ can be chosen explicitly (see~\eqref{eq:defCN}--\eqref{eq:defCT}, below), and we elaborate on this shortly before proving Theorem~\ref{thm: main} in Section~\ref{s:thmMainProof}, below.
\end{remark}

We now briefly discuss the significance of Properties~\ref{p:spectral}--\ref{property:eigenfunctions}.
It is well known that the convergence of a (reversible) Markov process can be studied using the spectral decomposition (see for instance Chapter 12 in~\cite{LevinPeres17}).
In particular, the rate of convergence is controlled by the \emph{spectral gap}, which in our case is simply~$\lambda_{2, \epsilon}$.
The spectral gap is not easy to estimate in general, and is an active area of study~\cite{
  Mic15,
  bakry2014analysis
}.

In our situation, Property~\ref{p:spectral} allows for the spectral gap to be exponentially small.
(The lower bound~\eqref{e:egvalLEpsilon} is in fact sharp (see~\cite[Chapter 8, Proposition 2.2]{Kolokoltsov00}, or~\eqref{eq:egvalub}, in Section~\ref{sec:Cpsi}, below), and hence the spectral gap is in fact exponentially small.)
Hence the standard argument using the spectral decomposition apriori only gives extremely slow convergence.
Indeed, equation~\eqref{e:fOverPiNu} implies
\begin{align}
  \norm[\Big]{\frac{f_t}{\pi_\epsilon} - 1}_{L^2(\pi_\epsilon)}^2
    &= \norm[\Big]{e^{-L_\epsilon t}\paren[\Big]{\frac{f_0}{\pi_\epsilon}} - 1}_{L^2(\pi_\epsilon)}^2
  \\
    \label{e:fOverPi}
    &\leq
      \exp\paren[\big]{-2\lambda_{2, \epsilon} t } \ip{f_0, \psi_{2, \epsilon}}_{L^2(\T^d)}^2
      + e^{-\Lambda t} \norm[\Big]{ \frac{f_0}{\pi_\epsilon} - 1 }_{L^2(\pi_\epsilon)}^2
    .
\end{align}
The extremely small spectral gap~\eqref{e:egvalLEpsilon} implies the first term decays extremely slowly, and is the bottleneck to convergence.

However, Property~\ref{p:spectral} also asserts that the third eigenvalue is large (i.e.\ bounded independent of~$\epsilon$), and this will give fast convergence \emph{provided} we control the projection onto the second eigenspace.
Variants of this idea have been used by several authors in many contexts to accelerate convergence~\cite{
  ConstantinKiselevEA08,
  KwokLauEA13,
  FengIyer19
}.
A warm start to Langevin dynamics will also control the projection on the second eigenspace, and this was recently used by Koehler, Lee, and Vuong~\cite{koehler2024efficientlylearningsamplingmultimodal} in a related multimodal sampling problem.

In our situation, the projection onto the second eigenspace is controlled by the term~$\ip{f_0, \psi_{2, \epsilon}}_{L^2(\T^d)}$.
We will see shortly (in Section~\ref{sec:propertycheck}, below) that the coefficient~$\ip{f_0, \psi_{2, \epsilon}}_{L^2(\T^d)}$ is essentially the mass imbalance.
More precisely~$\ip{f_0, \psi_{2, \epsilon}}_{L^2(\T^d)}$ is comparable to the difference between~$f_0$-mass and the~$\pi_\epsilon$-mass of each well.
This provides a valuable insight into the convergence of~\eqref{e:KForward}:
\emph{if the mass imbalance is small, then~$f_t$ converges to~$\pi_\epsilon$ quickly.}

For double well potentials, both spectral bounds in Property~\ref{p:spectral} are well known.
The lower bound~\eqref{e:lambdaiLower} can be found in Propositions~2.1, 2.2 in Chapter~8 of~\cite{Kolokoltsov00}.
The lower bound in~\eqref{e:egvalLEpsilon} is more recent, and was proved in~\cite{MenzSchlichting14}.
We elaborate on this in Section~\ref{sec:prop69}, below, and show that Assumption~\ref{a:criticalpts} and~\ref{assumption: nondegeneracy} imply Property~\ref{p:spectral} (see Proposition~\ref{l: lower_bound_next_eigenvalue}).
The generalization of Property~\ref{p:spectral} when~$U$ has~$J$ wells will be described in Section~\ref{sec:multi}.
\smallskip

Property~\ref{p:spectral} alone is not sufficient to prove Theorem~\ref{thm: main}, and we now discuss why Property~\ref{p:var} is required.
In order to obtain time complexity bounds for ASMC, we recall that the ASMC algorithm moves samples through a sequence of intermediate distributions.
The above insight tells us that in order to have any hope of obtaining a polynomial time-complexity bound, we need the intermediate levels (at least when $k$ is large) to be started with a small mass imbalance.
To ensure this, we need the distributions~$\pi_{\epsilon}$, and the corresponding eigenfunctions~$\psi_{2, \epsilon}$ to vary in a controlled way as~$\epsilon$ decreases.
The bounds we need in our proof are precisely~\eqref{e:PsiEpPi} and~\eqref{e:PsiEpEpP} listed in Property~\ref{p:var}.
We will later prove (Proposition~\ref{prop:pvar}, below) that Assumption~\ref{a:criticalpts} and~\ref{a:massRatioBound} imply Property~\ref{p:var} holds with~$g(\epsilon) \defeq \pi_\epsilon(\Omega_1)$.

  Finally Property~\ref{property:eigenfunctions} is a uniform bound on~$\norm{\psi_{2, \epsilon}}_{L^\infty}$, which is needed for technical reasons.
  We prove that Assumptions~\ref{a:criticalpts} and~\ref{a:massRatioBound} imply Property~\ref{property:eigenfunctions} in Section~\ref{sec:Cpsi}, below.

\subsubsection{Error estimates}
Returning to the proof of Theorem~\ref{thm: main},
we need an estimate on the Monte Carlo error when using~$N$ independent realizations to compute the integral of a test function.
We recall the error estimate for approximation the integral with resect to the Gibbs measure~$\pi_\epsilon$ by empirical average over the output of the 
standard Langevin Monte Carlo algorithm (LMC)
\begin{equation}
  \int_{\T^d} h \pi_\epsilon \, dx
    \approx
    \frac{1}{N}\sum_{i=1}^{N}h(Y^{i}_{\epsilon, t})
    ,
\end{equation}
where~$Y^{1}_{\epsilon , \cdot}$, \dots, $Y^{N}_{\epsilon, \cdot}$ are~$N$ independent solutions to~\eqref{e:Langevin}.
The right hand side approaches the left hand side as~$N, t \to \infty$. 
Our first lemma controlling the error is as follows.

\begin{lemma}\label{l:langevinError}
  Assume that for each $i \in \set{1, \dots, N}$, $\PDF(Y_{\epsilon,0}^{i}) = q^{i}_{\epsilon,0}$.  
Then  for any bounded test function $h$,
  \begin{equation}\label{e:langevinError}
    \Err_{\epsilon,T}(h)
    \leq e^{-\lambda_{2, \epsilon} T} \abs[\Big]{\int_{\mathbb T^d}h\psi_{2,\epsilon}\pi_{\epsilon} \, dx} \Err_{\epsilon,0}(\psi_{2,\epsilon})
    +\frac{1}{2\sqrt{N}}\|h\|_{\osc}
    + \mathcal E_{\epsilon,T}(h)
\end{equation}
where
\begin{equation}\label{eq: def of tilde epsilon}
  \mathcal E_{\epsilon,T}(h)\defeq \|h\|_{\osc}e^{-\Lambda T}\max_{i=1,\dots,N}\norm[\Big]{\frac{q^{i}_{\epsilon,0}}{\pi_{\epsilon}}}_{L^{\infty}(\pi_{\epsilon})}^{\frac12} .
\end{equation}
\end{lemma}

Notice that in order to use Lemma~\ref{l:langevinError} and Lemma~\ref{lem: first level}, we need an estimate for~$\norm{q_{k,0}/\pi_{k}}_{L^{\infty}(\pi_{k})}$.
This is addressed in the following lemma.
\begin{lemma}\label{lem: probability density L infty norm between levels}
  For every $2\leq k\leq M$, $1 \leq i \leq N$, let $q^{i}_{k,0}$ be the probability density function of~$X^{i}_{k,0}$.
  For any~$T_0 > 0$, there exists a constant $C_q = C_q(U, T_0)$ such that if $T\geq T_0$, then
  \begin{equation}\label{eq: probability density L infty norm between levels}
      \max_{i=1,\dots,N} \norm[\Big]{\frac{q^{i}_{k,0}}{\pi_{k}}}_{L^{\infty}(\pi_{k})}\leq C_{q}\exp\Big(\|U\|_{\osc}\Big(\frac{1}{\eta_k}-1\Big)\Big)
       .
  \end{equation}
\end{lemma}
Clearly the right hand side of~\eqref{eq: probability density L infty norm between levels} is exponentially large, which, at first sight, is concerning.
We will, however, only use~\eqref{eq: probability density L infty norm between levels} in~\eqref{eq: def of tilde epsilon} which has an exponentially small~$e^{-\Lambda T}$ factor.
Choosing~$T \geq O(1/\eta_k)$ will allow us to control it.
\medskip
We will now use Lemma \ref{l:resampling} and Lemma \ref{l:langevinError} to derive Monte Carlo error estimates between levels~$k$ and~$k+1$ in Algorithm~\ref{a:ASMC}.
Recall in Algorithm~\ref{a:ASMC}, $M$ is chosen according to~\eqref{e:MTN}, $\eta_1 = 1$, $\eta_M = \eta$, and the reciprocals~$1/\eta_1$, \dots, $1/\eta_M$ are linearly spaced.
That is~$\eta_k$ is chosen according to
\begin{equation}\label{e:chooseEtaK}
  \eta_k
    \defeq \frac{(M - 1) \eta}{(M-1) \eta + (k-1)(1 - \eta) }
    .
\end{equation}

For simplicity of notation, we use a subscript of~$k$ on the error, eigenvalue and eigenfunction to denote the corresponding quantities at~$\epsilon = \eta_k$.
Explicitly, we write
\begin{equation}
    \lambda_{2,k}\defeq \lambda_{2,\eta_k}
    , \quad
    \psi_{2,k}\defeq \psi_{2,\eta_k}
    , \quad\text{and}\quad
    \Err_{k, 0}(\psi_{2, k})\defeq \Err_{\eta_k, 0}(\psi_{2, k})
    .
\end{equation}
In order to prove Theorem~\ref{thm: main}, we need additionally ensure that the normalized eigenfunctions~$\psi_k$ and~$\psi_{k+1}$ are close.

The main idea behind the proof of Theorem~\ref{thm: main} is to first estimate ~$\Err_{k+1, 0}(\psi_{2, k+1})$ in terms of~$\Err_{k, 0}(\psi_{2, k})$, and then use Lemma~\ref{l:langevinError} to obtain~\eqref{e:MCerror}.
Obtaining this
recurrence relation, however, requires a fair amount of technical work. We state this
in the next lemma.
\begin{lemma}\label{l:iteration}
 Fix $\delta>0$, let~$C_q = C_q(U, 1)$ be the constant from Lemma~\ref{lem: probability density L infty norm between levels} with~$T_0 = 1$.
 Suppose that Properties~\ref{p:spectral},~\ref{p:var} and~\ref{property:eigenfunctions} hold. Define $\tilde{C}_{N}$ by
\begin{equation}\label{eq:deftildeCN}
    \tilde{C}_{N}\defeq 4\paren[\Big]{C_{\psi}\paren[\Big]{1+\frac{3}{2}C_r}}^2
    .
\end{equation}

For any~$\alpha> 0$, let
\begin{equation}\label{eq:Hgammachoice}
    H\defeq(1+\alpha)^{\frac12}\hat{U}>\hat{U},\quad \text{and}\quad\gamma\defeq \frac{\hat{\gamma}}{(1+\alpha)^{\frac12}}<\hat{\gamma}
  ,
\end{equation}
and choose
\begin{align} C_{\alpha}&\defeq \frac{1}{C_{H}} \log(C_{\psi}C_r+C_r^{\frac12}), \label{eq:deftildeCalpha}\\
  C_{\beta}&\defeq \exp\paren[\Big]{(V_{(0,1]}(g) +1)(C_2+C_{\psi} C_1)} \label{eq: defCbeta}
\end{align}
  where~$C_r$ is the constant from~\eqref{eq:defCr} and~$C_1, C_2, C_\psi$ and~$g$ are from Properties~\ref{p:var} and~\ref{property:eigenfunctions}.
If
   \begin{align}\label{eq:TNpropN}
            N&\geq\tilde{C}_{N} \frac{M^2}{\delta^2}
	    ,\\
        T&\geq\max\Bigl\{\frac{1}{\Lambda}\Big(\log(\frac{1}{\delta})+ \frac{\|U\|_{\osc}}{2\eta}+\log(M)+\log(6C_rC_{\psi}C_q^{\frac12})\Big),1,\\
	\label{eq:TNpropT}
        &\qquad~C_{\alpha}M^{(1+\alpha)\hat{\gamma}_r}\Bigr\},
    \end{align}
  then there exists (explicit) constants~$\beta_k, c_k$ such that for every~$k \in \set{2, \dots, M-1}$ we have
  \begin{equation}\label{e:iteration}
  \Err_{k+1,0}(\psi_{2,k+1})
     \leq
      \beta_k \Err_{k,0}(\psi_{2,k})
      + c_k
    ,
\end{equation}
and
  \begin{equation}\label{e:CBeta}
    \prod_{j=k}^{M-1} \beta_j \leq C_{\beta}
    \quad\text{and}\quad
    c_k \leq \frac{\delta}{M}
    .
  \end{equation}
\end{lemma}

The proof of Theorem~\ref{thm: main} now reduces to solving the recurrence relation~\eqref{e:iteration} and using Lemma~\ref{l:langevinError}. Finally, we will need an estimate for $\Err_{2,0}(\psi_{2,2})$ which we obtain as the mixing time when~$k = 1$ is of order~$1$.
\begin{lemma}\label{lem: first level}
  Define $\tilde C_{1} = \tilde C_1(U)$ by
  \begin{align}
    \tilde{C}_1&\defeq 1+\max\Bigl\{ \frac{1}{\Lambda}
    \max\bigl\{ 1,\log(12C_{\psi}C_{r}C_q^{\frac12})\bigr\},
    \\
    \label{eq:tildeC1}
    &\qquad\frac{\exp(H)}{C_{H}}\max\bigl\{ 1,\log(4C^2_{\psi}(C_r+1))\bigr\}\Bigr\}
    ,
  \end{align}
  where the constants~$C_r$, $C_\psi$, and~$H$ are defined
  in~\eqref{e:piKDef},
  \eqref{eq:defCpsi},
  and~\eqref{eq:Hgammachoice} respectively, and the constant~$C_q$ is from the statement of Lemma~\ref{l:iteration}.
  For any $\delta>0$, we have
  \begin{equation}\label{e:ErrPsi2}
    \Err_{2,0}(\psi_{2,2})
       \leq \delta
    \quad\text{provided}\quad
    T\geq \tilde C_{1}\Big(\log\paren[\Big]{\frac{1}{\delta}}+1\Big),\quad N\geq \frac{\tilde{C}_{N}}{\delta^2}
    .
  \end{equation}
  Here $\tilde{C}_{N}$ is the same constant as in Lemma \ref{l:iteration}.
\end{lemma}

\subsubsection{Proof of Theorem~\ref{thm: main}}\label{s:thmMainProof}
We are now well-equipped to prove Theorem~\ref{thm: main}.
Our proof will show that the constants~$C_N$ and~$C_T$ in~\eqref{e:MTN} are given by
\begin{equation}\label{eq:defCN}
    C_N\defeq 16C_{\beta}^2\tilde{C}_N\overset{\eqref{eq:deftildeCN}}{=} 64C_{\beta}^2\paren[\Big]{C_{\psi}\paren[\Big]{1+\frac{3}{2}C_r}}^2,
\end{equation}
and
\begin{equation}\label{eq:defCT}
    C_T\defeq \max\Bigl\{ C_{\alpha}, \tilde{C}_1\paren{1+\log(4C_{\beta})}, \frac{\norm{U}_{\osc}}{2\Lambda},
   \frac{1}{\Lambda}\paren[\big]{\log(6C_rC_{\psi}C_q^{\frac12}))+\log(4C_{\beta})}\Bigr\}.
\end{equation}
From this we see that the dimensional dependence of the constants~$C_N$ and~$C_T$ is only through the dimensional dependence of the constants arising in Properties~\ref{p:spectral}--\ref{property:eigenfunctions}, and the constant~$C_q$.
Recall in Lemma~\ref{lem: probability density L infty norm between levels}, the constant~$C_q(U, T_0)$ is the~$L^1 \to L^\infty$ bound for solutions to~\eqref{e:fOverPiNu} with~$\epsilon = 1$ after time~$T_0$.
In the statement of Lemma~\ref{l:iteration} we chose~$T_0 = 1$, and set~$C_q = C_q(U, 1)$.
We note, however, that~$C_q(U, T_0)$ is decreasing with~$T_0$, and by making~$T_0$ larger we can bound~$C_q$ independent of dimension.
As a result, the main unknown dimensional dependence in the constants~$C_T$ and~$C_N$ is through the constants in Properties~\ref{p:spectral}--\ref{property:eigenfunctions}.

\begin{proof}[Proof of Theorem~\ref{thm: main}]
 Fix $\alpha, \delta>0$, and define
  \begin{equation}\label{e:TildeDelta}
      \tilde{\delta}=
    \frac{\delta}{4C_{\beta}},
  \end{equation}
where $C_{\beta}$ is the constant in Lemma \ref{l:iteration}. Choose $M$ as in \eqref{e:MTN} and it can be directly checked that, if~$T,N$ satisfy~\eqref{e:MTN} with~$C_N$, $C_T$ defined as in~\eqref{eq:defCN} and~\eqref{eq:defCT}, then we have
  \begin{align}
    \label{eq: parametersT}
    T &\geq \max\Bigl\{
      C_{\alpha}M^{(1+\alpha)\hat{\gamma}_r}
      ,~
      \tilde C_1\Big(\log\paren[\Big]{\frac{1}{\tilde{\delta}}}+1\Big)
      ,~1,
    \\
       &\qquad
       \frac{1}{\Lambda}\Big(\log\paren[\Big]{\frac{1}{ \tilde{\delta}}}+\frac{\norm{U}_\osc}{2\eta} +\log(M)+\log(6C_rC_{\psi}C_q^{\frac12})\Big)
       \Bigr\}
       ,\\
    \label{eq: parametersN}
	   N&\geq\frac{\tilde{C}_{N} M^2}{\tilde{\delta}^2}
	   .
    \end{align}

  We claim~\eqref{e:MCerror} holds for any bounded test function~$h \in L^\infty(\T^d)$.
  To see this, we use Lemma~\ref{l:langevinError}, to obtain
    \begin{align}
   \Err_{M,T}(h) 
      &\leq \abs[\Big]{\int_{\mathbb T^d}h\psi_{2,k}\pi_{k} \, dx} e^{-\lambda_{2, k} T}\Err_{M,0}(\psi_{2,M})
    \\
      \label{e:ErrMT1}
      &\qquad+\frac{1}{2\sqrt{N}}\|h\|_{\osc}+\mathcal E_{\eta,T}(h)
    .
\end{align}
 We will now show that the right hand side of~\eqref{e:ErrMT1} is bounded above by~$\delta \norm{h}_\osc$.
  For the first term, a direct calculation using~\eqref{e:iteration} immediately shows that for~$T, N$ as in~\eqref{eq: parametersT}, \eqref{eq: parametersN} we have
\begin{align}\label{eq: est E psiM}
  \Err_{M,0}(\psi_{2,M}) &\leq \paren[\Big]{\prod_{j=2}^{M-1}\beta_j}\Err_{2,0}(\psi_{2,2})
       +\sum_{k=2}^{M-2}c_k \Big(\prod_{j=k+1}^{M-1}\beta_j\Big)+c_{M-1}\\
       &\overset{\mathclap{\eqref{e:CBeta},~\eqref{e:ErrPsi2}}}{\leq} \quad~~ C_{\beta}\tilde{\delta}+\sum_{k=2}^{M-2}C_{\beta}\frac{\tilde{\delta}}{M} +\frac{\tilde{\delta}}{M}
       \leq 2C_{\beta}\tilde{\delta} \overset{\eqref{e:TildeDelta}}{=} \frac{\delta}{2}
       .
\end{align}

Next, we see
  \begin{equation}\label{e:hLinf}
    \frac{1}{\sqrt{N}}
      \overset{\eqref{eq: parametersN}}{\leq}
      \frac{\tilde{\delta}}{M}
      \overset{\eqref{e:TildeDelta}}{\leq}
      \frac{\delta}{4}
      .
\end{equation}
Finally,
  \begin{align}
    \mathcal E_{\eta,T}(h)\overset{\eqref{eq: def of tilde epsilon}}{=}&\|h\|_{\osc} e^{-\Lambda T}\max_{i=1,\dots,N}\norm[\Big]{\frac{q^{i}_{M,0}}{\pi_{M}}}^{\frac12}_{L^{\infty}(\pi_{M})}\\
    \overset{\eqref{eq: probability density L infty norm between levels}}{\leq}&
    C_q^{\frac12}\|h\|_{\osc} \exp\paren[\Big]{\frac{\|U\|_{\osc}}{2\eta}} e^{-\Lambda T}\\
    \label{e:calEh}
   \overset{\eqref{eq: parametersT}}{\leq}&\|h\|_{\osc}\tilde{\delta}< \|h\|_{\osc}\frac{\delta}{4}
    .
\end{align}

Using~\eqref{eq: est E psiM}, \eqref{e:hLinf} and~\eqref{e:calEh} in~\eqref{e:ErrMT1} implies
\begin{equation}
  \Err_{M, T}(h) \leq \delta\|h\|_{\osc}
  .
\end{equation}
This proves~\eqref{e:MCerror}, concluding the proof.
\end{proof}

It remains to prove Lemmas~\ref{l:langevinError}, \ref{l:iteration},~\ref{lem: probability density L infty norm between levels} and \ref{lem: first level}, which will be done in subsequent sections.

\section{Proof of Lemmas for the Local Mixing Model} \label{s:localmodelproof}

In this section, we prove Lemmas~\ref{lem:mixNeedMass}, \ref{l:resampling} and~\ref{lem:iterativelocal} that were used Section~\ref{s:localModel} to prove Theorem~\ref{t:localMixingModel}.
We also prove the bound for~$\norm{r_k}_{L^\infty}$ stated in~\eqref{eq:CrExplicit}, that may be easier to use in practice.
Since the ideas used in Lemma~\ref{lem:mixNeedMass} and~\ref{lem:iterativelocal} are related, we prove Lemma~\ref{l:resampling} first.

\subsection{The Resampling Error (Lemma~\ref{l:resampling})}
Notice that the points~$y^{1}$, \dots, $y_N$ chosen according to~\eqref{e:PyiEqXj} are identically distributed, but need not be independent.
However, given the points~$x^{1}$, \dots, $x_N$, the points~$y^{1}$, \dots, $y_N$ are (conditionally) independent.
The main idea behind the proof of Lemma~\ref{l:resampling} is to split the error into the sum of a conditional mean, and a conditional standard deviation, and use conditional independence of~$y^{1}$, \dots, $y_N$.
\begin{proof}[Proof of Lemma \ref{l:resampling}]
    For simplicity of notation, let
  \begin{equation}\label{eq: defTildeh}
      \mathbf{x}\defeq \{x^1,\dots,x^{N}\},\quad\text{and}\quad \tilde{h}\defeq h-\int_{\mathcal{X}} hq\, dx
      ,
\end{equation}
  and let~$\E_{\mathbf x}$ denote the conditional expectation given the~$\sigma$-algebra generated by~$\mathbf x$.
    By the tower property,
    \begin{align}
    \E \paren[\Big]{ \frac{1}{N} \sum_{i=1}^N h(y^{i}) - \int_{\mathcal{X}} hq\, dx }^2
      &= \E \paren[\Big]{ \frac{1}{N} \sum_{i=1}^N \tilde{h}(y^{i})}^2
    \\
      \label{e:Ehy1}
       &=  \E\Big[\E_{\mathbf{x}}\Big(\frac{1}{N}\sum_{i=1}^N
       \tilde{h}(y^{i})\Big)^2\Big]
       .
    \end{align}
    We write
   \begin{equation}\label{e:J1J2}
 \E_{\mathbf{x}}\Big(\frac{1}{N}\sum_{i=1}^N
       \tilde{h}(y^{i})\Big)^2
      =J_1+J_2
      ,
\end{equation}
where
\begin{equation}
        J_1\defeq \Big(\frac{1}{N}\sum_{i=1}^{N}\E_{\mathbf{x}}\tilde{h}(y^{i})\Big)^2\\
	\quad\text{and}\quad
        J_2\defeq \E_{\mathbf{x}}\Big(\frac{1}{N}\sum_{i=1}^N
      \tilde{h}(y^{i})-\frac{1}{N}\sum_{i=1}^{N}\E_{\mathbf{x}}\tilde{h}(y^{i})\Big)^2
      .
\end{equation}

Notice that the points~$y^1$, \dots, $y^N$ are not independent; however, when conditioned on~$\mathbf x$, the points~$y^i$ are independent and identically distributed.
Thus
  \begin{equation}\label{e:J2}
 J_2\leq \frac{1}{N^2}\sum_{i=1}^{N}\var_{\mathbf{x}}(\tilde{h}(y^{i}))\leq \frac{1}{N}\|\tilde{h}\|^2_{L^{\infty}}
  .
\end{equation}

To bound~$J_1$, we note
  \begin{align}
    \frac{1}{N}\sum_{i=1}^{N}\E_{\mathbf{x}}\tilde{h}(y^{i}) &=\E_{\mathbf{x}}\tilde{h}(y^{1})
	  =\frac{ \sum_{i=1}^N
	 \tilde{h}(x^{i}) \tilde r(x^{i})}{\sum_{i=1}^N \tilde r(x^{i}) }
	 \overset{\eqref{e:rdef}}{=}
	 \frac{ \sum_{i=1}^N
	 \tilde{h}(x^{i}) {r}(x^{i})}{\sum_{i=1}^N {r}(x^{i}) }
    \\
      \label{e:J1sqrt}
      &= J_3 + J_4,
\end{align}
where
\begin{equation}
        J_3=\frac{ \sum_{i=1}^N
	 \tilde{h}(x^{i}) r(x^{i})}{\sum_{i=1}^N r(x^{i}) }\paren[\Big]{1-\frac{1}{N}\sum_{i=1}^N r(x^{i}) },
	 \quad\text{and}\quad
   J_4= \frac{1}{N}\sum_{i=1}^N
	  \tilde{h}(x^{i}) r(x^{i}).
\end{equation}
Clearly 
\begin{equation}\label{e:j3}
  \abs{J_3} \leq \norm{\tilde h}_{L^\infty} \paren[\Big]{1-\frac{1}{N}\sum_{i=1}^N r(x^{i}) }
  .
\end{equation}

Thus
\begin{align}
  \norm[\Big]{ \frac{1}{N} \sum_{i=1}^N h(y^{i}) - \int_{\mathcal{X}} hq\, dx }_{L^2(\P)}
    &\overset{\eqref{e:Ehy1}}\leq 
      (\E J_1)^{1/2} + (\E J_2)^{1/2}
  \\
    &\overset{\mathclap{\eqref{e:J1sqrt},~\eqref{e:J2}}}{\leq}
      \quad
      (\E J_3^2)^{1/2}
      + (\E J_4^2)^{1/2}
      + \frac{\norm{\tilde h}_{L^\infty}}{\sqrt{N}}
    \,.
\end{align}
Using the definition of~$\tilde h$ in~\eqref{eq: defTildeh} we obtain~\eqref{e:resampling} as desired.
\end{proof}

\subsection{The Monte Carlo error (Lemma~\ref{lem:mixNeedMass})}
In this section, we prove Lemma~\ref{lem:mixNeedMass} which provides an estimate for the error when using independent realizations of the process~$Y_{\epsilon, \cdot}$ to compute Monte Carlo integrals.
\begin{proof}[Proof of Lemma \ref{lem:mixNeedMass}]
  Since
  \begin{equation}
    \sum_{j=1}^{J}
      \paren[\Big]{
	1 - \frac{\mu_0(\Omega_j)}{\pi_\epsilon(\Omega_j)}
      }
      \pi_\epsilon(\Omega_j) = 0
      ,
  \end{equation}
  both sides of~\eqref{eq:mc} remain unchanged when a constant is added to the function~$h$.
  Thus without loss of generality we may replace~$h$ with~$h - \inf h + \frac{1}{2} \norm{h}_\osc$, and assume~$\norm{h}_{L^\infty} = \frac{1}{2} \norm{h}_\osc$.
Next we write
  \begin{equation}\label{e:mcHminusIntH}
     \frac{1}{N} \sum_1^N h(Y^i_{\epsilon, T})
        - \int_{\mathcal X} h \, \pi_\epsilon \, dx=I_1+I_2
    ,
\end{equation}
where
\begin{align}
    I_1\defeq
      \frac{1}{N} \sum_1^N \paren[\big]{ h(Y^i_{\epsilon, T}) - \E_0 h(Y^i_{\epsilon, T}) }
      ,\qquad
       I_2\defeq
       \frac{1}{N} \E_0 \sum_1^N h(Y^i_{\epsilon, T})
    -\int_{\mathcal X} h \, \pi_\epsilon \, dx
    ,
  \end{align}
  and~$\E_0$ denote the conditional expectation with respect to~$\sigma(y^1, \dots, y^N)$.

Now since $Y^{i}_{\epsilon,\cdot}$ are conditionally independent given~$y^1$, \dots, $y^N$,
  \begin{equation}\label{e:E0I12}
     \E_0 I_1^2
     = \frac{1}{N} \sum_1^N \E_0 \paren[\big]{ h(Y^i_{\epsilon, T}) - \E_0 h(Y^i_{\epsilon, T}) }^2
    \leq \frac{1}{N}\norm{h}_{L^{\infty}}^2
    =  \frac{1}{4N}\norm{h}^2_{\osc}
    .
\end{equation}

  Next using~\eqref{e:YTransitionKernel} and~\eqref{e:mu0} we compute
\begin{align}
  \frac{1}{N} \E_0 \sum_1^N h(Y^i_{\epsilon, T})
    &= \frac{1}{N} \sum_{i=1}^N \int p_{T}^{\epsilon}(y^i,z)h(z)\, dz
  \\
    &= 
	(1 - \chi^{T}_{\epsilon}) \int_{\mathcal X} h  \pi_\epsilon \,  d x +
	\chi^{T}_{\epsilon}\sum_{j = 1}^J \frac{\mu_0(\Omega_j)}{\pi_{\epsilon}(\Omega_j)}\int_{\Omega_j} h  \pi_\epsilon \,  d x
	,
\end{align}
and hence
\begin{equation}\label{e:I2simp}
  I_2 = \chi^{T}_{\epsilon}\sum_{j=1}^{J}
	 \paren[\Big]{\frac{\mu_0(\Omega_j)}{\pi_{\epsilon}(\Omega_j)} - 1}\int_{\Omega_j}h\pi_{\epsilon}\, d x
	 .
\end{equation}

Using~\eqref{e:E0I12} in~\eqref{e:mcHminusIntH} implies
  \begin{equation}
    \Err_{\epsilon, T}(h) = \norm{I_1 + I_2}_{L^2(\P)}
      \leq \norm{I_1}_{L^2(\P)} + \norm{I_2}_{L^2(\P)}
      \leq \frac{\norm{h}_\osc}{2 \sqrt{N}} + \norm{I_2}_{L^2(\P)}
      .
  \end{equation}
  Using~\eqref{e:I2simp} in the above implies~\eqref{eq:mc} as desired. Finally, \eqref{e:ErrTErr0} follows immediately from~\eqref{eq:mc} and the fact that for every~$j \in \set{1, \dots, J}$ we have
    \begin{equation}
      \abs[\Big]{\int_{\Omega_j}h\pi_{\epsilon}\, d x}
	\leq \norm{h}_\infty \int_{\Omega_j}\pi_{\epsilon}\, d x
	= \frac{\norm{h}_\osc}{2} \pi_\epsilon(\Omega_j)
	.
	\qedhere
    \end{equation}
\end{proof}

\subsection{A recurrence relation for the error (Lemma \ref{lem:iterativelocal})}
We now prove Lemma~\ref{lem:iterativelocal}, which obtains a recurrence relation for the Monte Carlo error between levels~$k$ and~$k+1$ in Algorithm~\ref{a:ASMC}.
\begin{proof}[Proof of Lemma \ref{lem:iterativelocal}]
  Fix $k \in \set{1,\dots,M-1}$, and $\ell\in \set{1,\dots, J}$.
Applying Lemma~\ref{l:resampling} with~$p=\pi_{k}$, $q=\pi_{k+1}$, $h=\one_{\Omega_{\ell}}$, and
$x^{i}=Y_{k,T}^{i}$ gives
\begin{equation}\label{eq:resamplingLoc}
    \Err_{k+1,0}(\one_{\Omega_{\ell}})
      \leq
	\frac{1}{\sqrt{N}} 
	+ \Err_{k,T}(r_k)
	+ \Err_{k,T}\paren[\big]{r_k(\one_{\Omega_{\ell}}-\pi_{k+1}(\Omega_{\ell}))}.
  \end{equation}

We now bound the last two terms on the right hand side of \eqref{eq:resamplingLoc}.
  Applying Lemma~\ref{lem:mixNeedMass} with~$h = r_k$, $\epsilon = \eta_k$ and using~\eqref{e:defrk} gives
  \begin{align}
    \Err_{k,T}(r_k)
      &\leq \frac{\norm{r_k}_{\osc}}{2\sqrt{N}}
	+ \chi_k^T
	  \norm[\bigg]{ \sum_{j = 1}^J
	    \paren[\Big]{1 - \frac{\mu_0(\Omega_j)}{\pi_k(\Omega_j)}} \pi_{k+1}(\Omega_j)
	  }_{L^2(\P)}
      \\
      &= \frac{\norm{r_k}_{\osc}}{2\sqrt{N}}
	+ \chi_k^T
	  \norm[\bigg]{\sum_{j = 1}^J
	    \paren[\Big]{1 - \frac{\mu_0(\Omega_j)}{\pi_k(\Omega_j)}}
	    (\pi_{k+1}(\Omega_j) - \pi_k(\Omega_j)
	  }_{L^2(\P)}
      \\
      \label{e:ErrRkT}
      &\leq \frac{\norm{r_k}_{\osc}}{2\sqrt{N}}
	+ \chi_k^T
	  \paren[\Big]{\sum_{j=1}^J
	    \abs[\Big]{\frac{\pi_{k+1}(\Omega_j)}{\pi_k(\Omega_j)} - 1}}
	  \max_{1 \leq j \leq J}  \Err_{k, 0}(\one_{\Omega_j})
      .
  \end{align}
  Here~$\mu_0$ is defined by~\eqref{e:mu0} with~$y^{i} = Y^i_{k, 0}$, and the second equality above is true because~$\sum_j \mu_0(\Omega_j) = 1$.

  For the last term on the right of~\eqref{eq:resamplingLoc} again apply Lemma~\ref{lem:mixNeedMass} with~$\epsilon = \eta_k$ and~$h = r_k ( \one_{\Omega_\ell} - \pi_{k+1}(\Omega_\ell) )$ to obtain
  \begin{equation}\label{eq:termrhloc}
      \Err_{k,T}\paren[\big]{r_k(\one_{\Omega_{\ell}}-\pi_{k+1}(\Omega_{\ell}))} \leq \frac{\norm{r_k}_{\osc}}{2\sqrt{N}}+\chi^{T}_{\epsilon}J_1
  \end{equation}
where
\begin{equation}
  J_1 \defeq \norm[\Big]{
      \sum_{j=1}^J
	\paren[\Big]{ 1 - \frac{\mu_0(\Omega_j)}{\pi_k(\Omega_j)} }
	\pi_{k+1}(\Omega_\ell)
	\paren{ \delta_{j, \ell} - \pi_{k+1}(\Omega_j) }
    }_{L^2(\P)}
    \,.
\end{equation}
Here $\delta_{j,\ell}$ is the Kronecker delta defined by
\begin{equation}
 \delta_{j,\ell}=\left\{ \begin{array}{ll}
      1&  \text{if }j=\ell,\\
     0 & \text{if }j\neq\ell.
 \end{array}\right. 
\end{equation}
Since~$\sum_j \mu_0(\Omega_j) = 1$ we note
\begin{align}
  J_1 &= \norm[\Big]{
      \sum_{j=1}^J
	\paren[\Big]{ 1 - \frac{\mu_0(\Omega_j)}{\pi_k(\Omega_j)} }
	\pi_{k+1}(\Omega_\ell)
	\paren{ \delta_{j, \ell} - \pi_{k+1}(\Omega_j) + \pi_k(\Omega_j)}
    }_{L^2(\P)}
    \\
      &\leq
      \pi_{k+1}(\Omega_\ell)
      \sum_{j=1}^J
	\abs[\Big]{ \frac{\delta_{j, \ell}}{\pi_k(\Omega_j)} + 1 - \frac{\pi_{k+1}(\Omega_j)}{\pi_k(\Omega_j)}}
      \max_{1 \leq j \leq J} \Err_{k, 0}(\one_{\Omega_j})
    \\
      \label{e:J1}
      &\leq
      \paren[\Big]{
	1 + \pi(\Omega_\ell)
	\sum_{j=1}^J
	  \abs[\Big]{ 1 - \frac{\pi_{k+1}(\Omega_j)}{\pi_k(\Omega_j)}}
      }
      \max_{1 \leq j \leq J} \Err_{k, 0}(\one_{\Omega_j})
    \,.
\end{align}
Using~\eqref{e:ErrRkT}, \eqref{eq:termrhloc} and~\eqref{e:J1} in~\eqref{eq:resamplingLoc} yields~\eqref{eq:iterativelocal} as desired.
\end{proof}

\section{Langevin error estimates (Lemmas~\ref{l:langevinError} and~\ref{lem: probability density L infty norm between levels})}\label{s:langevinErrorProof}
In this section we prove Lemmas~\ref{l:langevinError} and~\ref{lem: probability density L infty norm between levels}.
The proof of Lemma~\ref{l:langevinError} is based on a spectral decomposition, and is presented in Section~\ref{s:langevinError}, below.
The proof of Lemma~\ref{lem: probability density L infty norm between levels} is based on the maximum principle and is presented in Section~\ref{s:LinfBound}, below.

\subsection{The Monte Carlo error in the Langevin Step (Lemma~\ref{l:langevinError})}\label{s:langevinError}
The proof of Lemma \ref{l:langevinError} has three main steps.
First, we separate the error into the sum of the conditional mean and the conditional standard deviation.
The conditional standard deviation is~$O(1/\sqrt{N})$ and is easily bounded.
Then we decompose the conditional mean as the sum of  $\Err_{\epsilon,0}(\psi_{2,\epsilon})$ and a remainder using the spectral decomposition of $L_{\epsilon}$, and bound the remainder terms.
\begin{proof}[Proof of Lemma \ref{l:langevinError}] 
Since $\int_{\T^d}\psi_{2,\epsilon}\pi_{\epsilon}\, dx=0$, both sides of~\eqref{e:langevinError} remain unchanged when a constant is added to the function~$h$.
  Thus without loss of generality we may again replace~$h$ with~$h - \inf h + \frac{1}{2} \norm{h}_\osc$, and assume~$\norm{h}_{L^\infty} = \frac{1}{2} \norm{h}_\osc$.
As before, let $\E_0$ denote the conditional expectation given the~$\sigma$-algebra generated by~$\{Y_{\epsilon,0}^1,\dots,Y_{\epsilon,0}^{N}\}$.

\restartsteps
  \step
    By the tower property of conditional expectation, we have
\begin{equation}\label{eq: tower property level k}
       (\Err_{\epsilon,T}(h))^2
       = \E \Big[\E_0\Big(\frac{1}{N}\sum_{i=1}^{N}h(Y^{i}_{\epsilon, T})-\int_{\mathbb T^d}h\pi_{\epsilon}\,d x\Big)^2   \Big].
\end{equation}
Observe that
\begin{equation}
    \frac{1}{N}\sum_{i=1}^{N}h(Y^{i}_{\epsilon, T})-\int_{\mathbb T^d}h\pi_{\epsilon}\,d x=I_1+I_2,
\end{equation}
    where
    \begin{align}
            I_1&\defeq \frac{1}{N}\sum_{i=1}^{N}h(Y^{i}_{\epsilon,T})-\frac{1}{N}\sum_{i=1}^{N}\E_0h(Y^{i}_{\epsilon,T}),\\
            I_2&\defeq \frac{1}{N}\sum_{i=1}^{N}\E_0h(Y^{i}_{\epsilon,T})-\int_{\mathbb T^d}h\pi_{\epsilon}\,d x.
    \end{align}
  Notice $\E_0 I_1=0$, and $I_2$ is $\sigma(Y_{\epsilon,0}^1,\dots,Y_{\epsilon,0}^{N})$-measurable.
  Hence
\begin{equation}\label{eq: var+mean square level k}
    \E_0\Big(\frac{1}{N}\sum_{i=1}^{N}h(Y^{i}_{\epsilon, T})-\int_{\mathbb T^d}h\pi_{\epsilon}\,d x\Big)^2   
        =\E_0 I_1^2+I_2^2.
\end{equation}
Next, notice that after conditioning on $\{Y_{\epsilon,0}^1,\dots,Y_{\epsilon,0}^{N}\}$, the random variables $Y^{i}_{\epsilon,T}$ are independent.
Hence
\begin{equation}\label{eq: 1/N var level k}
    \E_0 I_1^2
     =\frac{1}{N^2}\sum_{i=1}^N \var_0(h(Y^{i}_{\epsilon,T}))\leq \frac{1}{N}\|h\|^2_{L^{\infty}}.
\end{equation}

Thus, we conclude
\begin{align}\label{e:errhstep1}
  (\Err_{\epsilon,T}(h))^2  
       &\overset{\eqref{eq: tower property level k}}{=} \E \Big[\E_0\Big(\frac{1}{N}\sum_{i=1}^{N}h(Y^{i}_{\epsilon, T})-\int_{\mathbb T^d}h\pi_{\epsilon}\,d x\Big)^2   \Big]\\
       &\overset{\eqref{eq: var+mean square level k}}{=}\E\E_0 I_1^2
        +\E I_2^2
        \overset{\eqref{eq: 1/N var level k}}{\leq}  \frac{1}{N}\|h\|^2_{L^{\infty}}+\E I_2^2.
\end{align}

\step
  In this step, we use a spectral decomposition to rewrite~$I_2$.
  Notice that $ h$ can be decomposed into components along the subspace spanned by $\{1,\psi_{2,\epsilon}\}$ and its orthogonal complement.
  This decomposition gives
  \begin{equation}
    h = \int_{\mathbb T^d}h\pi_{\epsilon}\,d x
      + f_0 + f_0^\perp
    ,
  \end{equation}
  where
\begin{align}
  \label{e:f0def}
  f_0(y) &\defeq 
    \paren[\Big]{\int_{\mathbb T^d}h\psi_{2,\epsilon}\pi_{\epsilon}\, dx}\psi_{2,\epsilon}(y)\\
  \label{e:f0perpdef}
  f_0^\perp(y)  &\defeq h(y)-\int_{\mathbb T^d}h\pi_{\epsilon}\,d x-\paren[\Big]{\int_{\mathbb T^d}h\psi_{2,\epsilon}\pi_{\epsilon}\,dx}\psi_{2,\epsilon}(y).
\end{align}
Therefore,
\begin{equation}\label{eq:I2sep}
  I_2 =\frac{1}{N}\sum_{i=1}^{N}\E_0h(Y^{i}_{\epsilon,T})-\int_{\mathbb T^d}h\pi_{\epsilon}\,d x
      =I_3+I_4,
\end{equation}
where
\begin{equation}
  I_3\defeq \frac{1}{N}\sum_{i=1}^{N}\E_0 f_0(Y^{i}_{\epsilon,T}),
  \quad\text{and}\quad
  I_4\defeq \frac{1}{N}\sum_{i=1}^{N}\E_0 f_0^\perp(Y^{i}_{\epsilon,T}).
\end{equation}

\substep[$I_3$ bound]
To bound bound~$I_3$, we first claim
\begin{equation}\label{e:Psi2T}
\E_0\psi_{2,\epsilon}(Y^{i}_{\epsilon,T})=e^{-\lambda_{2, \epsilon} T}\psi_{2,\epsilon}(Y^{i}_{\epsilon,0})
\,.
\end{equation}
To see this, recall that for any~$g_0 \in L^\infty(\T^d)$ the function~$g$ defined by
\begin{equation}
  g_t(y) = \E^y g_0( Y_{\epsilon, t} ),
\end{equation}
solves the Kolmogorov backward equation 
\begin{equation}\label{e:KBackward}
   \partial_t g + L_\epsilon g =0,
\end{equation}
with initial data~$g_0$.
Here~$Y_{\epsilon, \cdot}$ is a solution to the Langevin equation~\eqref{e:Langevin}.
Since~$\psi_{2, \epsilon}$ is the second eigenfunction of the operator~$L_\epsilon$ (defined in~\eqref{e:Ldef}), we see
\[
  \E^y \psi_{2, \epsilon}(Y_{\epsilon, t}) = e^{-\lambda_{2, \epsilon} t} \psi_{2, \epsilon}(y)
  ,
\]
which immediately implies~\eqref{e:Psi2T}.

Now~\eqref{e:f0def} and~\eqref{e:Psi2T} imply
\begin{equation}
\begin{split}
     I_3=&\frac{1}{N}\sum_{i=1}^{N}\E_0\Big[\paren[\Big]{\int_{\mathbb T^d}h\psi_{2,\epsilon}\pi_{\epsilon}\,dx}\psi_{2,\epsilon}(Y^{i}_{\epsilon,T})\Big]\\
     =&
    \paren[\Big]{\int_{\mathbb T^d}h\psi_{2,\epsilon}\pi_{\epsilon}\,dx}\frac{1}{N}\sum_{i=1}^{N}e^{-\lambda_{2, \epsilon} T}\psi_{2,\epsilon}(Y^{i}_{\epsilon,0}),
\end{split}
\end{equation}
and hence
\begin{align}
    (\E I_3^2)^{\frac12}&= \abs[\Big]{\int_{\mathbb T^d}h\psi_{2,\epsilon}\pi_{\epsilon}\,dx}e^{-\lambda_{2, \epsilon} T}\norm[\Big]{\frac{1}{N}\sum_{i=1}^{N}\psi_{2,\epsilon}(Y^{i}_{\epsilon,0})}_{L^2(\P)}\\
    &=\abs[\Big]{\int_{\mathbb T^d}h\psi_{2,\epsilon}\pi_{\epsilon}\,dx}e^{-\lambda_{2, \epsilon} T}\Err_{\epsilon,0}(\psi_{2,\epsilon}).
    \label{eq:I3}
\end{align}

\substep[$I_4$ bound]
To bound~$I_4$ we note
\begin{equation}\label{eq:I4}
  I_4=\frac{1}{N}\sum_{i=1}^{N}\E_0 f_0^\perp(Y^{i}_{\epsilon,T})=\frac{1}{N}\sum_{i=1}^{N}f_T^\perp(Y^{i}_{\epsilon,0})
    ,
\end{equation}
where~$f^\perp_t(y) = \E^y f^\perp_0(Y_{\epsilon, t})$ and, as before, $Y_{\epsilon, \cdot}$ is a solution of~\eqref{e:Langevin}.
Observe that, 
\begin{align}
 \E I_4^2
    &= \E\paren[\Big]{\frac{1}{N}\sum_{i=1}^{N} f_T^\perp(Y^{i}_{\epsilon,0})}^2
        \leq  \frac{1}{N}\sum_{i=1}^{N}\E f_T^\perp(Y^{i}_{\epsilon,0})^2
	= \frac{1}{N}\sum_{i=1}^{N}\int_{\mathbb T^d} (f_T^\perp)^2 q^{i}_{\epsilon,0}\, d x\\
  \label{eq:Ef2}
	&=\frac{1}{N}\sum_{i=1}^{N}\int_{\mathbb T^d} (f_T^\perp)^2\frac{q^{i}_{\epsilon,0}}{\pi_{\epsilon}}\pi_{\epsilon}\, d x
       \leq   \|f^\perp_T\|^2_{L^2(\pi_{\epsilon})}\max_{i=1,\dots,N}\norm[\Big]{\frac{q^{i}_{\epsilon,0}}{\pi_{\epsilon}}}_{L^{\infty}(\pi_{\epsilon})}.
\end{align}

To bound~$\norm{f_T^{\perp}}_{L^2(\pi_\epsilon)}^2$, we note that~$f^\perp$ solves the Kolmogorov backward equation~\eqref{e:KBackward}, and hence we have the spectral decomposition
\begin{equation}
    \norm{f_T^{\perp}}_{L^2(\pi_\epsilon)}^2=\sum_{i=1}^{\infty}e^{-2\lambda_{i,\epsilon}T}\abs[\Big]{\int f_0^{\perp}\psi_{i,\epsilon}\, \pi_\epsilon\, dx}^2\\
\end{equation}
Using~\eqref{e:f0def} and~\eqref{e:f0perpdef} the first two terms on the right vanish, and hence
the spectral decomposition gives
\begin{equation}\label{eq:fTL2}
    \norm{f_T^{\perp}}_{L^2(\pi_\epsilon)}^2
    =\sum_{i=3}^{\infty}e^{-2\lambda_{i,\epsilon}T}\abs[\Big]{\int f_0^{\perp}\psi_{i,\epsilon}\, \pi_\epsilon\, dx}^2\overset{\eqref{e:egvalLEpsilon}}{\leq} \|f_0^{\perp}\|^2_{L^2(\pi_{\epsilon})}e^{-2\Lambda T}.
\end{equation}

We will now bound $\|f_0^{\perp}\|_{L^2(\pi_{\epsilon})}$.
Notice
\begin{align}
  \|f_0^{\perp}\|_{L^2(\pi_{\epsilon})}&\leq
        \norm[\Big]{h-\int_{\mathbb T^d}h\pi_{\epsilon}}_{L^2(\pi_{\epsilon})}+\norm[\Big]{ \paren[\Big]{\int_{\mathbb T^d}h\psi_{2,\epsilon}\pi_{\epsilon}}   \psi_{2,\epsilon}}_{L^2(\pi_{\epsilon})}\\
	&\leq  \|h\|_{L^2(\pi_{\epsilon})}+\Big|\int_{\mathbb T^d}h\psi_{2,\epsilon}\pi_{\epsilon}\Big|\leq 2\|h\|_{L^2(\pi_{\epsilon})}\leq 2\|h\|_{L^{\infty}}.
\end{align}
Together with \eqref{eq:fTL2} this gives
\begin{equation}\label{eq:fL2}
    \|f_T^{\perp}\|_{L^2(\pi_{\epsilon})}\leq 2\|h\|_{L^{\infty}}e^{-\Lambda T}.
\end{equation}
Therefore, plugging \eqref{eq:fL2} into \eqref{eq:Ef2} yields
\begin{equation}\label{eq:E I4}
    (\E I_4^2)^{\frac12}\leq 2\|h\|_{L^{\infty}(\pi_{\epsilon})}e^{-\Lambda T}\max_{i=1,\dots,N}\norm[\Big]{\frac{q^{i}_{\epsilon,0}}{\pi_{\epsilon}}}^{\frac12}_{L^{\infty}(\pi_{\epsilon})}= \mathcal E_{\epsilon,T}(h).
\end{equation}

\step
Based on the previous steps,
\begin{multline}
  (\Err_{\epsilon,T}(h))^2  
        \overset{\eqref{e:errhstep1}}{\leq}  \frac{1}{N}\|h\|^2_{L^{\infty}}+\E I_2^2  \overset{\eqref{eq:I2sep}}{\leq}\frac{1}{N}\|h\|^2_{L^{\infty}}+\paren[\Big]{  (\E I_3^2)^{\frac12}+\E (I_4^2)^{\frac12}}^2 \\
\overset{\eqref{eq:I3},\eqref{eq:E I4}}{\leq} \frac{1}{N}\|h\|^2_{L^{\infty}}
       +\Big( \abs[\Big]{\int_{\mathbb T^d}h\psi_{2,\epsilon}\pi_{\epsilon}\, d x} e^{-\lambda_{2, \epsilon} T}\Err_{\epsilon,0}(\psi_{2,\epsilon})+\mathcal E_{\epsilon,T}(h)\Big)^{2}.
\end{multline}
Taking square root on both sides and using $\norm{h}_{L^\infty} = \frac{1}{2} \norm{h}_\osc$ finish the proof.
\end{proof}

\subsection{Growth of \texorpdfstring{$\norm{q^{i}_{k,0}/\pi_{k}}_{L^{\infty}}$}{q/pi} (Lemma \ref{lem: probability density L infty norm between levels})}\label{s:LinfBound}
In this section, we prove Lemma \ref{lem: probability density L infty norm between levels} which will be used in the proof of Lemma \ref{l:iteration}, and was also used in the proof of Theorem~\ref{thm: main} to obtain~\eqref{e:calEh}.
Let $q^{i}_{k,t}$ denote the probability density of $X^{i}_{k,t}$.
The proof Lemma~\ref{lem: probability density L infty norm between levels} involves controlling the growth of $\norm{q^{i}_{k,t}/\pi_{k}}_{L^{\infty}}$ in the Langevin step, and in the resampling step. In the Langevin step, $\norm{q^{i}_{k,t}/\pi_{k}}_{L^{\infty}}$ is nonincreasing due to maximum principle. In the resampling step, the growth of $\|{q_{k,0}}/{\pi_{k}}\|_{L^{\infty}(\pi_{k})}$ between levels is tracked using duality.

\begin{proof}[Proof of Lemma \ref{lem: probability density L infty norm between levels}]
  The proof contains three steps.
  Fix~$T_0 > 0$.
  First by standard~$L^1 \to L^\infty$ bounds for~\eqref{e:fOverPiNu} with~$\epsilon = 1$,
  see \cite{Aronson}, there exists a constant~$C_q = C_q( U, T_0 )$ such that
\begin{equation}\label{eq:leqCq}
    \max_{i=1,\dots,N}\norm[\Big]{\frac{q^{i}_{1,T}}{\pi_1}}_{L^{\infty}}\leq C_q
    .
\end{equation}
Next we will show that for every $k\geq 2$, 
\begin{equation}\label{eq: induction}
   \max_{i=1,\dots,N} \norm[\Big]{\frac{q^{i}_{k,0}}{\pi_{k}}}_{L^{\infty}}\leq\frac{1}{\prod_{\ell=1}^{k-1}\min r_{\ell}}\max_{i=1,\dots,N} \norm[\Big]{\frac{q^{i}_{1,T}}{\pi_1}}_{L^{\infty}},
\end{equation}
where $r_k$ is defined in \eqref{e:defrk}.
 
Finally we show by direct computation that
\begin{equation}\label{eq:prodminrk}
    \prod_{\ell=1}^{k-1}\min r_{\ell}
    \geq \exp\paren[\Big]{\paren[\Big]{1-\frac{1}{\eta}}\|U\|_{\osc}}.
\end{equation}
Combining \eqref{eq:leqCq}, \eqref{eq: induction} and \eqref{eq:prodminrk} completes the proof.
  We will now show each of the above inequalities.

\restartsteps
  \step[Proof of~\eqref{eq:leqCq}]
  Since~$X^i_{1, \cdot}$ solves~\eqref{e:Langevin} with~$\epsilon = \eta_1 = 1$, its density, denoted by~$q^i_{1, t}$, must solve~\eqref{e:KForward} (with~$\epsilon = \eta_1$).
  Since~$L_1^*$ (defined by~\eqref{e:LStarDef}, with~$\epsilon = 1$) is nondegenerate, parabolic regularity implies there exists a constant~$C = C(U, T_0)$ such that
  \begin{equation}
   \max_{i=1,\dots,N} \norm{q^i_{1, T_0}}_{L^\infty}
      \leq C \norm{q^i_{1, T_0/2}}_{L^1}
      = C\,.
  \end{equation}
  The last equality above followed because~$q^i_{1, t}$ is a probability density and so~$\norm{q^i_{1, t}}_{L^1} = 1$ for every~$t > 0$.
  Thus
  \begin{equation}\label{eq:T0infity}
     \norm[\Big]{\frac{q^i_{1, T_0}}{\pi_1}}_{L^\infty}
      \leq \frac{\norm{q^i_{1, T_0}}_{L^\infty}}{\min \pi_1}
      \leq \frac{C}{\min \pi_1}
  \end{equation}
  and we choose~$C_q = C_q(U, T_0)$ to be the right hand side of the above.
  (We note that our state space is the compact torus~$\T^d$ and so~$\min \pi_1 > 0$.)
  
  Since $q^i_{1,t} / \pi_1 $ solves the backward equation~\eqref{e:KBackward} with $\epsilon=\eta_1 = 1$, the maximum principle and~\eqref{eq:T0infity} imply~\eqref{eq:leqCq} for all~$T \geq T_0$.

  \step[Proof of~\eqref{eq: induction}]
  We claim that for all~$k \in \set{1, \dots, M-1}$, we have
\begin{equation}\label{eq:induc1}
    \max_{i=1,\dots,N} \norm[\Big]{\frac{q^{i}_{k+1,0}}{\pi_{k+1}}}_{L^{\infty}}\leq\frac{1}{\min r_{k}}\max_{i=1,\dots,N} \norm[\Big]{\frac{q^{i}_{k,T}}{\pi_k}}_{L^{\infty}}
    .
\end{equation}
  Next we note that $q^i_{k,t} / \pi_k $ satisfies~\eqref{e:KBackward} with~$\epsilon = \eta_k$.
  Thus, by the maximum principle we have
  \begin{equation}\label{eq:induc2}
    \norm[\Big]{\frac{q^{i}_{k,T}}{\pi_{k}}}_{L^{\infty}}
      \leq  \norm[\Big]{\frac{q^{i}_{k,0}}{\pi_{k}}}_{L^{\infty}},
  \end{equation}
  for every~$k \in \set{2,\dots,M}$ and every~$i \in \set{1, \dots, N}$.
  The bound \eqref{eq: induction} immediately follows from \eqref{eq:induc1} and \eqref{eq:induc2}.
  Thus it only remains to prove~\eqref{eq:induc1}.
  \smallskip

  We note that if~$X^1_{k, T}$, \dots, $X^N_{k, T}$ were i.i.d.\, then one has an explicit formula for~$q^i_{k+1, 0}$, from which~\eqref{eq:induc1} follows immediately.
  In our situation these processes are not independent, and so we prove~\eqref{eq:induc1} using duality, and without relying on an explicit formula.

For any test function $h\in L^{1}(\pi_{k+1})$ we have
   \begin{align}
     \MoveEqLeft
       \int_{\mathbb T^d}h(x) q^{i}_{k+1,0}(x)\, dx
	  =\E h(X^{i}_{k+1,0})
	\\
	  &=\E \E \paren[\big]{ h(X^{i}_{k+1,0}) \given X^{1}_{k,T},\dots,X^{N}_{k,T}}
	   =\E\paren[\bigg]{\frac{\sum_j h(X^{j}_{k,T})r_k(X^{j}_{k,T})}{\sum_j r_k(X^{j}_{k,T})}}
	\\
	  \label{eq: holder0}
	  &\leq\frac{1}{N\min r_k} \sum_{j=1}^N \E \abs[\big]{ h(X^{j}_{k,T})r_k(X^{j}_{k,T}) }
	  .
   \end{align}
Next, we note that for every $j=1,\dots,N$,
\begin{align}
  \E \abs[\big]{h(X^{j}_{k,T})r_k(X^{j}_{k,T}) }
	&=\int_{\mathbb T^d}|h|r_kq^{j}_{k,T}\, d x
        \overset{\eqref{e:defrk}}{=}\int_{\mathbb T^d}|h|\frac{r_kq^{j}_{k,T}}{r_k\pi_{k}}\pi_{k+1}\, d x
    \\
    \label{eq:Ehrk}
	&\leq \|h\|_{L^1(\pi_{k+1})}\Big\|\frac{q^{j}_{k,T}}{\pi_{k}}\Big\|_{L^{\infty}}.
\end{align}

Thus~\eqref{eq: holder0} and~\eqref{eq:Ehrk} imply
\begin{equation}
  \abs[\Big]{\int_{\mathbb T^d}h(x) q^{i}_{k+1,0}(x)\, dx }
    \leq 
      \frac{\norm{h}_{L^1(\pi_{k+1})}}{\min r_k}\max_{j=1,\dots,N}\Big\|\frac{q^{j}_{k,T}}{\pi_{k}}\Big\|_{L^{\infty}}
      ,
\end{equation}
from which~\eqref{eq:induc1} follows by duality.

\step[Proof of~\eqref{eq:prodminrk}]
Observe that 
\begin{equation}\label{e:rkExplicit}
  r_k(x) = \frac{Z_k}{Z_{k+1}} \exp\paren[\Big]{ -\paren[\Big]{\frac{1}{\eta_{k+1}} - \frac{1}{\eta_k} } U(x) }
\end{equation}
where~$Z_k \defeq Z_{\eta_k}$, and~$Z_{\eta_k}$ is the normalization constant in~\eqref{e:piNu}.
Hence, for all~$k \in \set{1, \dots, M-1}$ the minimum of~$r_k$ is attained at the same point, which we denote by~$x^*$.
Thus,
\begin{equation}\label{e:minr1}
    \prod_{\ell=1}^{k-1}\min r_{\ell}
      = \prod_{\ell=1}^{k-1} r_{\ell}(x^*)
	=\frac{\pi_{k}(x^*)}{\pi_1(x^*)}
    =\frac{Z_1}{Z_k}\exp\paren[\Big]{\paren[\Big]{1-\frac{1}{\eta_k}}U(x^*)}\\
  .
\end{equation}
Since~$U \geq 0$ by assumption, and $\eta_k < \eta_1 = 1$, we must have $Z_1 \geq Z_k$.
Using this in~\eqref{e:minr1} immediately implies~\eqref{eq:prodminrk} as desired.
This completes the proof of Lemma~\ref{lem: probability density L infty norm between levels}.
\end{proof}

\section{Iterating error estimates (Lemmas \ref{l:iteration} and \ref{lem: first level})}\label{sec:ProofIteration}

Lemma ~\ref{l:iteration} consists of two main parts: the derivation of recurrence relation~\eqref{e:iteration}, and obtaining the estimate~\eqref{e:CBeta} for~$\beta_k$, $c_k$.
We do each of these steps in Sections~\ref{sec: recurrence},~\ref{sec:estCj} and~\ref{sec: est coeff}.
We combine these and prove Lemma~\ref{l:iteration} at the end of Section~\ref{sec: est coeff}.

\subsection{Recurrence relation}\label{sec: recurrence}

We will now prove~\eqref{e:iteration} by combining the estimate for the Monte Carlo error (Lemma~\ref{l:langevinError}) and the resampling error (Lemma~\ref{l:resampling}).
For clarity, we state this as a new lemma and give explicit formulae for the constants~$\beta_k$ and~$c_k$ appearing in~\eqref{e:iteration}.

\begin{lemma}\label{lem: iterative scheme between levels}
  For each $2\leq k\leq M-1$, the inequality~\eqref{e:iteration} holds with~$\beta_k$ and~$c_k$ given by
\begin{align}
    \label{eq: defBetak}
    \beta_k &\defeq e^{-\lambda_{2, k} T}\Big(\Big|\int_{\mathbb T^d}r_k\psi_{2,k}\pi_{k}\, dx\Big|\cdot\|\psi_{2,k+1}\|_{L^{\infty}}
    \\
      &\qquad\qquad\qquad \mathbin{+}\Big|\int_{\mathbb T^d}\psi_{2,k+1}\psi_{2,k}\pi_{k+1}\, dx\Big|\Big)
    \\
  \label{eq: defck}
    c_k &\defeq 3\|\psi_{2,k+1}\|_{L^{\infty}}\norm{r_k}_{L^{\infty}}\Big(\frac{1}{2\sqrt{N}}+e^{-\Lambda T}\max_{i=1,\dots,N}\norm[\Big]{\frac{q^{i}_{k,0}}{\pi_{k}}}^{\frac12}_{L^{\infty}(\pi_{k})} \Big)\\
    &\qquad+\frac{1}{\sqrt{N}}\|\psi_{2, k+1}\|_{L^{\infty}}
    .
\end{align}
\end{lemma}

\begin{proof}[Proof of Lemma \ref{lem: iterative scheme between levels}]
Applying Lemma \ref{l:resampling} with
\[
p=\pi_{k},\quad q=\pi_{k+1},\quad h=\psi_{2,k+1},\quad
x^{i}=X_{k,T}^{i},\quad y^{i}=X^{i}_{k+1,0}
\]
 gives
\begin{align}
  \MoveEqLeft
    \Err_{k+1,0}(\psi_{2,k+1})=\norm[\Big]{
    \frac{1}{N} \sum_{i=1}^N \psi_{2,k+1}(X^{i}_{k+1,0})}_{L^2(\P)}
      \\
      &\overset{\eqref{e:resampling}}{\leq}
	\frac{1}{\sqrt{N}} \norm{\psi_{2,k+1}}_{L^\infty}
	+ \norm{\psi_{2,k+1}}_{L^\infty} \norm[\Big]{1- \frac{1}{N}\sum_{\ell=1}^N r_{k}(X_{k,T}^i)}_{L^2(\P)}
      \\
	&\qquad+ \norm[\Big]{\frac{1}{N}\sum_{i=1}^N \psi_{2,k+1}(X_{k,T}^i)r_{k}(X_{k,T}^i)}_{L^2(\P)}\\
    \label{eq: ErrK+1Decomp}
    &=
	\frac{1}{\sqrt{N}} \norm{\psi_{2,k+1}}_{L^\infty}
	+ \norm{\psi_{2,k+1}}_{L^\infty} \Err_{k,T}(r_k-1)
	+ \Err_{k,T}(\psi_{2,k+1}r_k)
    .
\end{align}

  We will now bound the last two terms on the right of~\eqref{eq: ErrK+1Decomp}.
For the term $\Err_{k,T}(r_k)$, we apply Lemma~\ref{l:langevinError} with
\[
\epsilon=\eta_k, \quad h=r_k,\quad q_{\epsilon,0}=q_{k,0}
,
\]
to obtain
\begin{align}
    \Err_{k,T}(r_k)
	&\overset{\eqref{e:langevinError}}{\leq}\frac{1}{2\sqrt{N}}\norm{r_k}_{L^{\infty}}+\abs[\Big]{\int_{\mathbb T^d}r_k\psi_{2,k}\pi_{k}\, dx}e^{-\lambda_{2, k} T}\Err_{k,0}(\psi_{2,k})\\
    \label{eq: ErrR}
	&\qquad+\norm{r_k}_{L^{\infty}}e^{-\Lambda T}\max_{i=1,\dots,N}\norm[\Big]{\frac{q^{i}_{k,0}}{\pi_{k}}}^{\frac12}_{L^{\infty}(\pi_{k})} 
\end{align}
where we use the fact that $\|r_k\|_{\osc}\leq \norm{r_k}_{L^{\infty}}$.

Similarly, for the  term $\Err_{k,T}(\psi_{2,k+1}r_k)$, we apply Lemma \ref{l:langevinError} with
\[
\epsilon=\eta_k, \quad h=\psi_{2,k+1}r_k,\quad q_{\epsilon,0}=q_{k,0}
\] 
to obtain
\begin{align}
\Err_{k,T}(\psi_{2,k+1}r_k)
    &\overset{\eqref{e:langevinError}}{\leq}\abs[\Big]{\int_{\mathbb T^d}\psi_{2,k+1}\psi_{2,k}\pi_{k+1}\, dx}e^{-\lambda_{2, k} T}\Err_{k,0}(\psi_{2,k})\\
	&\qquad+\frac{1}{\sqrt{N}}\norm{r_k}_{L^{\infty}}\norm{\psi_{2,k+1}}_{L^{\infty}}
    \\
      \label{eq: ErrPsiR}
      &\qquad +2\norm{r_k}_{L^{\infty}}\|\psi_{2,k+1}\|_{L^{\infty}}e^{-\Lambda T}\max_{i=1,\dots,N}\norm[\Big]{\frac{q^{i}_{k,0}}{\pi_{k}}}^{\frac12}_{L^{\infty}(\pi_{k})}
      ,
\end{align}
where we use the fact that
\begin{equation}
    \norm{r_k\psi_{2,k+1}}_{\osc}\leq 2 \norm{r_k\psi_{2,k+1}}_{L^{\infty}}\leq 2\norm{r_k}_{L^{\infty}}\norm{\psi_{2,k+1}}_{L^{\infty}}.
\end{equation}
  Plugging~\eqref{eq: ErrR} and~\eqref{eq: ErrPsiR} into~\eqref{eq: ErrK+1Decomp} and using~\eqref{eq: defBetak},~\eqref{eq: defck} yields~\eqref{e:iteration}, completing the proof.
\end{proof}

\subsection{Estimate of \texorpdfstring{$c_k$}{ck}}\label{sec:estCj}

We will now show how~$N$ and~$T$ can be chosen so that we obtain the bound for~$c_k$ in~\eqref{e:CBeta}. Precisely, we have the following lemma.
\begin{lemma}[Estimate of $c_k$]\label{lem: cj}
Suppose Properties~\ref{p:spectral} and~\ref{property:eigenfunctions} hold.  Fix $\delta>0$, let~$C_q = C_q(U, 1)$ be the constant from Lemma~\ref{lem: probability density L infty norm between levels} with~$T_0 = 1$, and $C_{\psi}$, $C_{r}$ be the constants defined in \eqref{eq:defCpsi} and \eqref{eq:defCr} respectively.
  If~$N$ and~$T$ are chosen such that
  \begin{align}
     & N\geq  \tilde{C}_{N}\frac{M^2}{\delta^2}, \label{eq:Nck}\\
     &T\geq \max\Bigl\{\frac{1}{\Lambda}\Big(\log\paren[\Big]{\frac{1}{\delta}}+ \frac{\|U\|_{\osc}}{2\eta}+\log M+\log(6C_{\psi}C_{r}C_q^{\frac12})\Big),1\Bigr\}, \label{eq:Tck}
  \end{align}
then 
\begin{equation*}
  c_k\leq  \frac{\delta}{M},\quad \text{for every } 2\leq k\leq M-1.
\end{equation*}
\end{lemma}

\begin{proof}[Proof of Lemma \ref{lem: cj}]
We first rewrite \eqref{eq: defck} as
\begin{equation}\label{eq:cksp}
    c_k=I_1+I_2
\end{equation}
where 
\begin{align*}
        I_1&=\frac{1}{\sqrt{N}}\|\psi_{2,k+1}\|_{L^{\infty}}\paren[\Big]{1+\frac{3}{2}\norm{r_k}_{L^{\infty}}},\\
        I_2&=3\|\psi_{2,k+1}\|_{L^{\infty}}\norm{r_k}_{L^{\infty}}e^{-\Lambda T}\max_{i=1,\dots,N}\norm[\Big]{\frac{q^{i}_{k,0}}{\pi_{k}}}^{\frac12}_{L^{\infty}(\pi_{k})}.
\end{align*}
Notice that the choice of $T$ and $N$ gives
\begin{equation}\label{eq:cki1}
I_1\overset{\eqref{eq:defCr},\eqref{eq:defCpsi}}{\leq} \frac{1}{\sqrt{N}} C_{\psi}\paren[\Big]{1+\frac{3}{2}C_r}
    \overset{\eqref{eq:deftildeCN},\eqref{eq:Nck}}{\leq}   \frac{\delta}{2M}
    ,
\end{equation}
and
\begin{equation}\label{eq:ckeLambdaT}
  e^{-\Lambda T} \overset{\eqref{eq:Tck}}{\leq} \delta \exp\paren[\Big]{-\frac{\|U\|_{\osc}}{2\eta}}\frac{1}{M}\frac{1}{6C_rC_{\psi}C_q^{\frac12}}.
\end{equation}

Therefore,
\begin{equation}\label{eq:cki2}
    I_2\overset{\eqref{eq:defCpsi},\eqref{eq:defCr},\eqref{eq: probability density L infty norm between levels}}{\leq} 3C_{\psi}C_re^{-\Lambda T}C^{\frac12}_q\exp\Big(\frac{\|U\|_{\osc}}{2}\Big(\frac{1}{\eta_k}-1\Big) \Big)\overset{\eqref{eq:ckeLambdaT}}{\leq} \frac{\delta}{2M}.
\end{equation}
Using~\eqref{eq:cki1} and~\eqref{eq:cki2} in~\eqref{eq:cksp} concludes the proof.
\end{proof}

\subsection{Estimate of \texorpdfstring{$\beta_k$}{betak}}\label{sec: est coeff}

Recall from~\eqref{e:iteration} the error grows by a factor of~$\beta_k$ at each level, and so to prove Theorem~\ref{thm: main} we need to ensure~$\prod \beta_k$ remains bounded.
The main result in this section (Lemma~\ref{lem: betaj}, below) obtains this bound and shows that the first inequality in~\eqref{e:CBeta} holds.  For simplicity of notation, let
\begin{equation}\label{eq:defTheta}
      \Theta(k,k+1)\defeq \Big|\int_{\mathbb T^d}r_k\psi_{2,k}\pi_{k}\, dx\Big|\cdot\|\psi_{2,k+1}\|_{L^{\infty}}+\Big|\int_{\mathbb T^d}\psi_{2,k+1}\psi_{2,k}\pi_{k+1}\, dx\Big|.
  \end{equation}
  and note
\begin{equation}
        \beta_k=e^{-\lambda_{2,k}T}\Theta(k,k+1).
    \end{equation}

  Here it is where Property~\ref{property:eigenfunctions} comes to play. Indeed, in the low temperature regime, the second assertion in Property~\ref{property:eigenfunctions} guarantees the first term in~\eqref{eq:defTheta} is small, and second term is close to~$1$.
This makes~$\Theta(k, k+1)$ close to~$1$, and hence the product bounded when $\lambda_{2,k}$ is close to zero.

We will bound $\prod_{j=k}^{M-1}\beta_j$ differently when the temperature~$\eta_k$ is low, and when it is high.
First, when the temperature is low, the exponential factor $e^{-\lambda_{2,k}T}$ is very close to $1$, and does not help much.
In this case will show that the product $\prod_{j=k}^{M-1}\Theta(j,j+1)$ stays bounded, by approximating $\Theta(k,k+1)$ in terms of the mass in each well and estimating the mass distribution using small temperature asymptotics.
When the temperature is high, the small temperature asymptotics are not valid anymore.
However, in this case $\lambda_{2,k}$ is not too small, and can be used to ensure $\beta_k\leq 1$ in a relatively short amount of time. We now bound~$\prod \beta_j$ to obtain the first inequality in~\eqref{e:CBeta}.
\begin{lemma}[Estimate of $\beta_j$]\label{lem: betaj}
Suppose Properties~\ref{p:spectral},~\ref{p:var} and~\ref{property:eigenfunctions} hold. If
 at each step $T\geq C_{\alpha}M^{(1+\alpha)\hat{\gamma}_r}$ with $C_{\alpha}$ defined as in~\eqref{eq:deftildeCalpha}, then for $2\leq k\leq M-1$, the first inequality~\eqref{e:CBeta} holds with~$C_{\beta}$ defined by~\eqref{eq: defCbeta}.
\end{lemma}

\begin{proof}
Given a fixed $\alpha>0$, the choice~\eqref{eq:Hgammachoice} gives
\begin{equation}\label{eq:HoverGamma}
    \frac{H}{\gamma}=(1+\alpha)\hat{\gamma}_r.
\end{equation}

Choose~$H$ according to~\eqref{eq:Hgammachoice}.
  By Property~\ref{p:spectral}, there  exists $C_H>0$ independent of $\epsilon$ such that
\begin{equation}\label{eq:leqgeqLambda}
  C_{H}\exp\paren[\Big]{-\frac{H}{\epsilon}} 
    \overset{\eqref{e:egvalLEpsilon}}{\leq}
  \lambda_{2,\epsilon},
    \quad \text{for all } \epsilon<1.
        \end{equation}
We choose a critical temperature $\eta_{\mathrm{cr}}>0$ so that
 \begin{equation}\label{eq:defEtacr}
  \exp
\Big(-\frac{{\gamma}}{\eta_{\mathrm{cr}}}\Big)=\frac{1}{M}.
 \end{equation}
We will prove the first inequality in~\eqref{e:CBeta} holds by splitting the analysis into two cases.

\restartcases

\case[$\eta\geq\eta_{\mathrm{cr}}$]
  In this case, for every~$k \geq 2$ we have $\eta_k \geq \eta \geq \eta_{\mathrm{cr}}$ and so
\begin{equation}\label{eq:expgeqCgammaM}
    \exp
\Big(-\frac{{\gamma}}{\eta_k}\Big)\geq \frac{1}{ M}.
\end{equation}
This implies
\begin{equation}
    \lambda_{2,k}\overset{\eqref{eq:HoverGamma},\eqref{eq:leqgeqLambda}}{\geq}  C_{H}\Big(\exp
\Big(-\frac{{\gamma}}{\eta_k}\Big)\Big)^{(1+\alpha)\hat{\gamma}_r}
\label{eq:EtaLargeLambdaLb}\overset{\eqref{eq:expgeqCgammaM}}{\geq} C_{H}\Big(\frac{1}{M}\Big)^{(1+\alpha)\hat{\gamma}_r}.
\end{equation}

Therefore, for 
\begin{equation}\label{eq:betakTlb}
T\geq C_{\alpha}M^{(1+\alpha)\hat{\gamma}_r},\quad    
\end{equation}
we have that for $k\geq 2$,
\begin{equation}\label{eq:lambdaT}
    \lambda_{2,k} T\overset{\eqref{eq:deftildeCalpha},\eqref{eq:EtaLargeLambdaLb}}{\geq}  \log\paren[\big]{C_r^{\frac12}(C_{\psi}+1)}.
\end{equation}

Notice that by using Cauchy-Schwarz directly, we have
\begin{equation}\label{eq:intPsiPisimple}
    \|\psi_{2,k}\|_{L^2(\pi_{k+1})}\leq \Big( \int_{\mathbb T^d}(\psi_{2,k})^2r_{k}\pi_{k}\, dx\Big)^{\frac12}\leq \|r_{k}\|^{\frac12}_{L^{\infty}(\pi_{k})}.
\end{equation}

Therefore, 
\begin{align}
    \MoveEqLeft e^{-\lambda_{2, k} T}\Theta(k,k+1)\\
&\overset{
\mathclap{\eqref{eq:defTheta}}}{=}e^{-\lambda_{2, k} T}\Big( \Big|\int_{\mathbb T^d}\psi_{2,k}\pi_{k+1}\, dx\Big|\|\psi_{2,k+1}\|_{L^{\infty}}+\Big|\int_{\mathbb T^d}\psi_{2,k+1}\psi_{2,k}\pi_{k+1}\, dx\Big|\Big)\\
&\overset{\mathclap{\text{Cauchy-Schwarz}}}{\leq}\qquad e^{-\lambda_{2, k} T}\|\psi_{2,k}\|_{L^2(\pi_{k+1})}\paren[\big]{\|\psi_{2,k+1}\|_{L^{\infty}}+1}\\
& \\
&\label{eq:elambdaTleq1}\overset{\mathclap{\eqref{eq:defCpsi},\eqref{eq:intPsiPisimple}}}{\leq} \qquad e^{-\lambda_{2, k} T}\norm{r_k}_{L^{\infty}}^{\frac12}\paren[\big]{C_{\psi}+1}\overset{\eqref{eq:defCr},\eqref{eq:lambdaT}}{\leq} 1. 
\end{align}

We conclude that, if $\eta\geq \eta_{\mathrm{cr}}$, and $T$ satisfies \eqref{eq:betakTlb},
then for every $k$,
\begin{equation}\label{eq:etalarge}
     \prod_{j=k}^{M-1} \beta_j= \prod_{j=k}^{M-1} \paren[\Big]{e^{-\lambda_{2, k} T}\Theta(k,k+1)}\overset{\eqref{eq:elambdaTleq1}}{\leq} 1.
\end{equation}

\case[$\eta<\eta_{\mathrm{cr}}$]
Define $k_0$ by
\begin{equation}
    k_0\defeq\min\{2\leq k\leq M-1 \: | \: \eta_k\leq \eta_{\mathrm{cr}}\}.
\end{equation}
We  first consider $k>k_0$, in which case we have~$\eta_k < \eta_{\mathrm{cr}}$. To bound $\Theta(k,k+1)$, we apply Properties~\ref{p:var} and~\ref{property:eigenfunctions} with~$\epsilon=\eta_k$ and~$\epsilon'=\eta_{k+1}$, to obtain
\begin{align}
       \Theta(k,k+1)&\overset{\mathclap{\eqref{eq:defTheta}}}{\leq}  1+\paren{C_2+C_{\psi} C_1}\paren[\Big]{\exp\paren[\Big]{-\frac{\gamma}{\eta_k}}+\abs[\big]{g(\eta_{k+1})-g(\eta_k)}}\\
       &\overset{\mathclap{\eqref{eq:defEtacr}}}{\leq}  1+\paren{C_2+C_{\psi} C_1}\paren[\Big]{\frac{1}{M}+\abs[\big]{g(\eta_{k+1})-g(\eta_k)}}. \label{eq:thetakk1}
\end{align}
By direct computation, for $k\geq k_0$,
\begin{align}
   \prod_{j=k}^{M-1}\beta_j&\leq
  \prod_{j=k}^{M-1}\Theta(j,j+1)\\
&\overset{\mathclap{\eqref{eq:thetakk1}}}{\leq} ~\prod_{j=k}^{M-1}\paren[\Big]{1+\paren{C_2+C_{\psi} C_1}\paren[\Big]{\frac{1}{M}+\abs[\big]{g(\eta_{k+1})-g(\eta_k)}}}\\
       &\overset{\mathclap{\text{AM-GM}}}{\leq}  ~~ \bigg(
	1+\frac{C_2+C_{\psi} C_1}{M}
    \\
      &\qquad\qquad
    \mathbin{+} \frac{1}{M-k}\sum\limits_{j=k}^{M-1}\paren[\big]{C_2+C_{\psi} C_1}\cdot \abs[\big]{g(\eta_{k+1})-g(\eta_k)}
       \bigg)^{M-k}\\
	&\leq \paren*{  1+\frac{V_{(0,1]}(g) (C_2+C_{\psi} C_1) }{M-k}+\frac{C_2+C_{\psi} C_1}{M}}^{M-k}\overset{\eqref{eq: defCbeta}}{\leq} C_{\beta},\label{eq:jgeqk0beta}
\end{align}
where the last inequality uses the fact that $M-k\leq M$. 

\medskip

Next, for the case or $k < k_0$, we observe that~$\eta_k \geq \eta_{\mathrm{cr}}$.
Now using the same argument as in the case $\eta\geq \eta_{\mathrm{cr}}$ we see,
\begin{equation}\label{eq: jleqk0}
  \prod_{j=k}^{k_0-1} \beta_j   \leq 1,
\end{equation}
provided $T$ satisfies \eqref{eq:betakTlb}.
This implies that
\begin{equation}\label{eq:etasmallksmall}
   \prod_{j=k}^{M-1} \beta_j
       = \paren[\Big]{\prod_{j=k}^{k_0-1} \beta_j}\paren[\Big]{\prod_{j=k_0}^{M-1} \beta_j}
        \overset{\eqref{eq: jleqk0}}{\leq}  \prod_{j=k_0}^{M-1}\beta_j\overset{\eqref{eq:jgeqk0beta}}{\leq}C_{\beta}.
\end{equation}

Combining \eqref{eq:etalarge}, \eqref{eq:jgeqk0beta} and \eqref{eq:etasmallksmall} completes the proof.
\end{proof}

Now we prove Lemma \ref{l:iteration}.
The proof follows immediately from Lemmas~\ref{lem: iterative scheme between levels}--\ref{lem: betaj}.

\begin{proof}[Proof of Lemma \ref{l:iteration}] Choose $\alpha>0$. For any given $\delta>0$, we choose~$T, N$ as in~\eqref{eq:TNpropN}--\eqref{eq:TNpropT}.
    Using Lemmas \ref{lem: iterative scheme between levels}, \ref{lem: cj} and~\ref{lem: betaj} we obtain~\eqref{e:iteration} and~\eqref{e:CBeta} as desired.
\end{proof}

\subsection{Mixing at the First Level (Lemma \ref{lem: first level})}
\label{sec: first level}

In this section we prove Lemma \ref{lem: first level} which bounds~$\Err_{2,0}(\psi_{2,2})$.
We obtain this bound without relying on the initial mass distribution, but instead using fast mixing at first level. 

\begin{proof}[Proof of Lemma \ref{lem: first level}]
    The proof is analogous to that of Lemma~\ref{lem: iterative scheme between levels} except when applying Lemma \ref{l:langevinError}, we choose $q_{\epsilon,0}=q_{1,T_0}$ with~$T_0 = 1$. Precisely, for any bounded test function $h$,
    \begin{equation}\label{e:1stlevel}
    \Err_{1,T}(h)
    \leq e^{-\lambda_{2, 1} (T-1)} \abs[\Big]{\int_{\mathbb T^d}h\psi_{2,1}\pi_{1} \, dx} \Err_{\epsilon,1}(\psi_{2,1})
    +\frac{1}{2\sqrt{N}}\|h\|_{\osc}
    + \mathcal E_{1,T}(h)
\end{equation}
where
\begin{equation}
  \mathcal E_{1,T}(h)= \|h\|_{\osc}e^{-\Lambda (T-1)}\max_{i=1,\dots,N}\norm[\Big]{\frac{q^{i}_{1,1}}{\pi_{\epsilon}}}_{L^{\infty}(\pi_{1})}^{\frac12} \overset{\eqref{eq:leqCq}}{\leq} \|h\|_{\osc}e^{-\Lambda (T-1)}C_q^{\frac12}.
\end{equation}
Next, we apply~\eqref{e:1stlevel} with~$h=r_1$ and~$h=r_1\psi_{2,1}$ and then plug into Lemma~\ref{l:resampling}. After a direct calculation, we obtain that for~$T \geq T_0=1$, 
    \begin{equation*}
    \Err_{2,0}(\psi_{2,2}) 
       \leq 
   \beta_1\Err_{1,T_0}(\psi_{2,1})+c_1,
\end{equation*}
where
\begin{align}\label{eq:beta1}
  \beta_1&= e^{-\lambda_{2, 1} (T-1)}\Big(\Big|\int_{\mathbb T^d}r_1\psi_{2,1}\pi_{1}\, dx\Big|\cdot\|\psi_{2,2}\|_{L^{\infty}}+\Big|\int_{\mathbb T^d}\psi_{2,2}\psi_{2,1}\pi_{2}\, dx\Big|\Big)\\
    &\overset{\mathclap{\eqref{eq:defCr},\eqref{eq:defCpsi}}}{\leq} ~~e^{-\lambda_{2,1}(T-1)}C_{\psi}(C_r+1).
\end{align}
and
\begin{align}\label{eq:c1}
  c_1&=3\|\psi_{2,2}\|_{L^{\infty}}\|r_1\|_{\osc}\Big(\frac{1}{2\sqrt{N}}+e^{-\Lambda (T-1)}C_q^{\frac12} \Big)
    +\frac{1}{\sqrt{N}}\|\psi_{2, 2}\|_{L^{\infty}}\\
    &\overset{\mathclap{\eqref{eq:defCr},\eqref{eq:defCpsi}}}{\leq}~~C_{\psi}\paren{1+\frac{3}{2}C_{r}}\frac{1}{\sqrt{N}}+3C_{\psi}C_{r}C_q^{\frac12}
   e^{-\Lambda (T-1)}.
\end{align}
  Here~$C_q = C_q(U, 1)$ is from Lemma~\ref{lem: probability density L infty norm between levels} and $C_{\psi}$, $C_{r}$ are the constants defined in \eqref{eq:defCpsi} and \eqref{eq:defCr} respectively. Therefore, for a given $\delta$, it can be checked the~$T,N$ in~\eqref{e:ErrPsi2} with~$\tilde{C}_1$ defined by~\eqref{eq:tildeC1} satisfy the following
\begin{align}
    \label{eq:Nfirst} N&\geq \tilde{C}_N\frac{1}{\delta^2}\overset{\eqref{eq:deftildeCN}}{=}4\paren[\Big]{C_{\psi}\paren[\Big]{1+\frac{3}{2}C_r}}^2\frac{1}{\delta^2},\\
      \label{eq:Tfirst}   T&~~\overset{\mathclap{\eqref{e:egvalLEpsilon},\eqref{eq:tildeC1}}}{\geq}~~ 1 + \max\Bigl\{\frac{1}{\Lambda}\Big(\log(\frac{1}{\delta})+ \log(12C_{\psi}C_{r}C_q^{\frac12})\Big),\\
        &\qquad\qquad \frac{1}{\lambda_{2,1}}\Big( \log(\frac{1}{\delta})+\log(4C^2_{\psi}(C_r+1))\Big) \Bigr\}.
\end{align}

It remains to check \eqref{e:ErrPsi2}. Using the fact that
 \begin{equation}\label{eq:Errpsi21}
     \Err_{1,T_0}(\psi_{2,1})\leq \|\psi_{2,1}\|_{L^{\infty}}\overset{\eqref{eq:defCpsi}}{\leq} C_{\psi},
 \end{equation}
  and
    \eqref{eq:beta1},
    \eqref{eq:c1},
    \eqref{eq:Nfirst}
  and~\eqref{eq:Tfirst}
 we obtain
\begin{align}
     \Err_{2,0}(\psi_{2,2}) 
       &\leq 
   \beta_1\Err_{1,T_0}(\psi_{2,1})+c_1\\
   & \leq
   e^{-\lambda_{2,1}(T-1)}C^2_{\psi}(C_r+1)+C_{\psi}(1+\frac{3}{2}C_{r})\frac{1}{\sqrt{N}}+3C_{\psi}C_{r}C_q^{\frac12}
   e^{-\Lambda (T-1)}\\
    & \overset{\mathclap{\eqref{eq:Nfirst},\eqref{eq:Tfirst}}}{<}\qquad
   \frac{\delta}{4}+\frac{\delta}{2}+\frac{\delta}{4}=\delta. \qedhere
\end{align}
\end{proof}

\section{Energy valley estimates}\label{s:CBV}

\subsection{The Mass Ratio (Lemma \ref{lem: upper bound on mass ratio of two wells})}
In this section we prove Lemma \ref{lem: upper bound on mass ratio of two wells}, whose main idea is that when $\epsilon\to 0$, the value of integral $\int_{\Omega_i}\exp(-U/\epsilon)\, dx$ is mainly determined by landscape near the local minima.

\begin{proof}[Proof of Lemma \ref{lem: upper bound on mass ratio of two wells}]

  We will prove that there exists $C>0$ independent of  $\eta_{\min}$ such that
\begin{equation}\label{e:massRatioBoundProof}
  \sup_{\epsilon\in [\eta_{\min},\eta_{\max}]}
      \frac{\pi_{\epsilon}(\Omega_1)}{\pi_{\epsilon}(\Omega_2)}
	\leq C,
      \quad
      \sup_{\epsilon\in [\eta_{\min},\eta_{\max}]}
      \frac{\pi_{\epsilon}(\Omega_2)}{\pi_{\epsilon}(\Omega_1)}
	\leq C
      .
   \end{equation}
   Combing this with the fact that
   \begin{equation}
     \sup_{\epsilon\in [\eta_{\min},\eta_{\max}]} \frac{1}{\pi_{\epsilon}(\Omega_i)} \leq \sup_{\epsilon\in [\eta_{\min},\eta_{\max}]} \frac{\pi_{\epsilon}(\Omega_1)+\pi_{\epsilon}(\Omega_2)}{\pi_{\epsilon}(\Omega_i)} \overset{\eqref{e:massRatioBound}}{\leq}1+C\defeq C_m^2,
    \end{equation}
    for~$i \in \set{1, 2}$, we obtain~\eqref{e:massRatioBound} as desired.

    To prove \eqref{e:massRatioBoundProof}, we note that Assumption~\ref{a:criticalpts} implies, via Laplace method as in \eqref{eq: est integral exp(-U/eta)}, that there exists~$\epsilon_2 > 0$ such that for all $\epsilon < \epsilon_2$ we have
\begin{equation*}
  \frac{\pi_{\epsilon}(\Omega_1)}{\pi_{\epsilon}(\Omega_2)}\overset{\eqref{eq: est integral exp(-U/eta)}}{=}\frac{(\mathrm{det}(\nabla^2U(x_{\min,1})))^{-\frac12}}{(\mathrm{det}(\nabla^2U(x_{\min,2})))^{-\frac12}\exp\paren[\big]{\frac{|U(x_{\min,1})-U(x_{\min,2})|}{\epsilon}}}+O(\epsilon),
\end{equation*}
  The assumption~\eqref{eq:depthdifference} further shows that for all $\epsilon \in [\eta_{\min}, \epsilon_2]$, we have
\begin{equation*}
    \begin{split}
        \frac{\pi_{\epsilon}(\Omega_1)}{\pi_{\epsilon}(\Omega_2)}&\leq\frac{(\mathrm{det}(\nabla^2U(x_{\min,1})))^{-\frac12}}{(\mathrm{det}(\nabla^2U(x_{\min,2})))^{-\frac12}}+O(\epsilon),\\
	\llap{\text{and}\qquad}
        \frac{\pi_{\epsilon}(\Omega_1)}{\pi_{\epsilon}(\Omega_2)}&\geq\frac{(\mathrm{det}(\nabla^2U(x_{\min,1})))^{-\frac12}}{(\mathrm{det}(\nabla^2U(x_{\min,2})))^{-\frac12}\exp(C_{l})}+O(\epsilon).
    \end{split}
\end{equation*}
  Thus, making~$\epsilon_2$ smaller if necessary, for every $\epsilon\in [\eta_{\min},\epsilon_2]$, we have
\[
\frac12\frac{(\mathrm{det}(\nabla^2U(x_{\min,1})))^{-\frac12}}{(\mathrm{det}(\nabla^2U(x_{\min,2})))^{-\frac12}\exp(C_{l})}\leq \frac{\pi_{\epsilon}(\Omega_1)}{\pi_{\epsilon}(\Omega_2)}\leq \frac{3}{2}\frac{(\mathrm{det}(\nabla^2U(x_{\min,1})))^{-\frac12}}{(\mathrm{det}(\nabla^2U(x_{\min,2})))^{-\frac12}}.
\]
  Since~$\epsilon \mapsto \pi_\epsilon(\Omega_1) / \pi_\epsilon(\Omega_2)$ is continuous and positive on the interval~$[\epsilon_2, \eta_{\max}]$, we obtain~\eqref{e:massRatioBoundProof} as desired.
\end{proof}

Lemmas~\ref{lem: upper bound on mass ratio of two wells} and~\ref{lem: mass in well as BV} immediately shows the following corollary, showing the finiteness condition~\eqref{e:CBV} holds provided Assumptions~\ref{a:criticalpts}, \ref{assumption: nondegeneracy} hold and the wells have nearly equal depth.
\begin{corollary}\label{cor:eqcbvholds}
  Assume the function $U$ satisfies Assumption \ref{a:criticalpts}, Assumption \ref{assumption: nondegeneracy} and there exist~$\eta_{\min} \geq 0$ and~$C_\ell < \infty$ such that~\eqref{eq:depthdifference} holds.
  Then for any finite $\eta_{\max} > \eta_{\min}$ the constant~$C_\LBV$ in~\eqref{e:CBV} can be bounded above in terms of~$U$, $C_\ell$ and $\eta_{\max}$, but independent of~$\eta_{\min}$.
\end{corollary}

\begin{proof}[Proof of Corollary \ref{cor:eqcbvholds}]
Notice that the inequality \eqref{eq:massBV} applied to $U/\eta_{\max}$ gives the existence of a constant~$\tilde{C}_{\BV}$, independent of $\eta$, such that for $i=1,2$, we have
    \begin{equation}\label{eq:intpartial}
        \int_{\eta}^{\eta_{\max}}\abs{\partial_\epsilon  \pi_\epsilon(\Omega_i)}\, d\epsilon\leq \tilde{C}_{\BV}.
    \end{equation}
    Therefore,
    \begin{equation}
        \int_{\eta}^{\eta_{\max}}
    \abs{\partial_\epsilon \ln \pi_\epsilon(\Omega_i)}\, d\epsilon=\int_{\eta}^{\eta_{\max}}
    \frac{\abs{\partial_\epsilon  \pi_\epsilon(\Omega_i)}}{\pi_\epsilon(\Omega_i)} \, d\epsilon\overset{\eqref{e:massRatioBound},\eqref{eq:intpartial}}{\leq} C_m^2\tilde{C}_{\BV}.
    \end{equation}
    Taking $\eta\to 0$ on the left hand side finishes the proof.
\end{proof}

\subsection{BV bounds on \texorpdfstring{$\pi_{\epsilon}(\Omega_i)$}{piepsilonomegai} (Lemma \ref{lem: mass in well as BV})}\label{sec:BV}
In this section, we will show that the masses in the two wells stay away from~$0$ do not oscillate too much. 
In particular, we will show that  the quantity~$\pi_{\epsilon}(\Omega_i)$ is of bounded variation as a function of~$\epsilon$ provided~$U$ satisfies Assumption~\ref{a:criticalpts} holds. Precisely, we have the following lemma.

\begin{lemma}\label{lem: mass in well as BV}
  If~$U$ satisfies Assumption~\ref{a:criticalpts} then there exists a constant $C_{\BV}$ such that such that for every $\eta\in (0,1)$, and every $i\in \set{1,2}$ we have
\begin{equation}\label{eq:massBV}
	 \int_{\eta}^{1}\abs{\partial_{\epsilon}\pi_{\epsilon}(\Omega_i)}\, d \epsilon\leq C_{\BV}.
       \end{equation}
\end{lemma}
Notice that~\eqref{eq:massBV} combined with our assumption Assumption~\ref{a:massRatioBound} implies~\eqref{e:CBV} holds.
This was the condition required to obtain error estimates for ASMC applied to the local mixing model (Theorem~\ref{t:localMixingModel}). We also remark that Lemma~\ref{lem: mass in well as BV} also holds when~$U$ has more than two local minima. The proof in the multi-well case is analogous to that of double-well case, but the calculation is  heavier.

We now prove Lemma \ref{lem: mass in well as BV} by calculating the derivative of $\pi_{\epsilon}(\Omega_i)$ with respect to $\epsilon$. It turns out that the derivative either has a sign or stay bounded as $\epsilon\to 0$, which implies that $\pi_{\epsilon}(\Omega_i)$ is a BV function.

\begin{proof}[Proof of Lemma \ref{lem: mass in well as BV}]
 We only prove \eqref{eq:massBV} for~$i = 1$.
  In fact, if \eqref{eq:massBV} holds for~$i = 1$, then from the identity $\pi_{\epsilon}(\Omega_2)=1-\pi_{\epsilon}(\Omega_1)$, we see that \eqref{eq:massBV} holds for $\pi_{\epsilon}(\Omega_2)$ for the same constant $C_{\BV}$.

Using the definition of $\pi_{\epsilon}(\Omega_1)$, and the fact that $U\in C^6$, we see that $\pi_{\epsilon}(\Omega_1)$ is differentiable in $\epsilon$.
  We compute
\begin{align}
   \partial_{\epsilon}\pi_{\epsilon}(\Omega_1)&=\partial_{\epsilon}\left( \frac{\int_{\Omega_1}\exp(-\frac{U}{\epsilon})\, d x}{\int_{\Omega_1}\exp(-\frac{U}{\epsilon})\, d x+\int_{\Omega_2}\exp(-\frac{U}{\epsilon})\, d x}\right)\\
        &=\frac{1}{\epsilon^2}\frac{\Big(\int_{\Omega_1}\exp(-\frac{U}{\epsilon})U\, d x\Big)\Big(\int_{\Omega_2}\exp(-\frac{U}{\epsilon})\, d x\Big)}{\Big(\int_{\Omega_1}\exp(-\frac{U}{\epsilon})\, d x+\int_{\Omega_2}\exp(-\frac{U}{\epsilon})\, d x \Big)^2}\\
        &\qquad-\frac{1}{\epsilon^2}\frac{\Big(\int_{\Omega_2}\exp(-\frac{U}{\epsilon})U\, d x\Big)\Big(\int_{\Omega_1}\exp(-\frac{U}{\epsilon})\, d x\Big)}{\Big(\int_{\Omega_1}\exp(-\frac{U}{\epsilon})\, d x+\int_{\Omega_2}\exp(-\frac{U}{\epsilon})\, d x \Big)^2}. 
\label{eq:partialPiEpsilon}
\end{align}
We now split the analysis into two cases: $U(x_{\min,1})=U(x_{\min,2})$ and $U(x_{\min,1})<U(x_{\min,2})$.

\restartcases
  \case[$U(x_{\min,1})=U(x_{\min,2})$]
We start by estimating the integrals involved in \eqref{eq:partialPiEpsilon}. 
According to Proposition B4 and Remark under the proof of Proposition B4, (page 289-290) in \cite{Kolokoltsov00},
\begin{equation}\label{eq:integral laplace method in one well with U}
    \int_{\Omega_i}\exp \paren[\Big]{-\frac{U}{\epsilon}}U\, d x=(2\pi \epsilon)^{\frac{d}{2}}\frac{\epsilon}{2}\Big(\mathfrak{c}_i+O(\epsilon)\Big)
    ,
\end{equation}
  where $\mathfrak{c}_i = \mathfrak c_i(U)$ depends on derivatives of $U$ up to order $3$ evaluated at $x_{\min,i}$, and is independent of $\epsilon$.
  The $O(\epsilon)$ involves constants that may depend on derivatives of $U$ up to order $6$ evaluated at $x_{\min,i}$.
  On the other hand, Proposition B2 in \cite{Kolokoltsov00} guarantees
\begin{equation}\label{eq: est integral exp(-U/eta)}
    \int_{\Omega_i}e^{-\frac{U}{\epsilon}} \, dx =(2\pi\epsilon)^{\frac{d}{2}}\frac{\exp(-\frac{U(x_{\min,i})}{\epsilon})}{\sqrt{\det \nabla^2 U(x_{\min,i})}} (1+O(\epsilon)),
    \quad\text{for } i \in 1,2
    .
\end{equation}

Therefore, when $\epsilon$ is sufficiently small, combining \eqref{eq:integral laplace method in one well with U} and \eqref{eq: est integral exp(-U/eta)},
\begin{align}
\partial_{\epsilon}\pi_{\epsilon}(\Omega_1)
    &=
    \begin{multlined}[t]
    \frac{1}{\epsilon^2}\frac{ \Big( \frac{\epsilon}{2}(\mathfrak{c}_1+O(\epsilon))(\mathrm{det}(\nabla^2U(x_{\min,2})))^{-\frac12} \Big)}{\Big( \sum_{i=1,2}(\mathrm{det}(\nabla^2U(x_{\min,i})))^{-\frac12}\Big)^2}\\
    -\frac{1}{\epsilon^2}\frac{ \Big( \frac{\epsilon}{2}(\mathfrak{c}_2+O(\epsilon))(\mathrm{det}(\nabla^2U(x_{\min,1})))^{-\frac12} \Big)}{\Big( \sum_{i=1,2}(\mathrm{det}(\nabla^2U(x_{\min,i})))^{-\frac12}\Big)^2}.
    \end{multlined}
\end{align}

We now discuss two different cases.

  \emph{Case I.1: $\mathfrak{c}_1\mathrm{det}(\nabla^2U(x_{\min,2})))^{-\frac12}-\mathfrak{c}_2\mathrm{det}(\nabla^2U(x_{\min,1})))^{-\frac12}=0$.}
  In this case, 
\[
\abs{\partial_{\epsilon}\pi_{\epsilon}(\Omega_1)}=O(1),\quad \epsilon \to 0.
\]
  Since the function  $\epsilon \mapsto \partial_{\epsilon}\pi_{\epsilon}(\Omega_1))$ is a continuous function on~$(0, 1]$, it must be bounded which implies~\eqref{eq:massBV}.

\emph{Case I.2: $\mathfrak{c}_1\mathrm{det}(\nabla^2U(x_{\min,2})))^{-\frac12}-\mathfrak{c}_2\mathrm{det}(\nabla^2U(x_{\min,1})))^{-\frac12}\neq 0$.}
Without loss of generality we assume~$\mathfrak{c}_1\mathrm{det}(\nabla^2U(x_{\min,2})))^{-\frac12}-\mathfrak{c}_2\mathrm{det}(\nabla^2U(x_{\min,1})))^{-\frac12} < 0$.
In this case there must exist some $\epsilon_{\mathrm{cr}}>0$ such that $\partial_{\epsilon}\pi_{\epsilon}(\Omega_1)) < 0$ for all $\epsilon\in (0, \epsilon_{\mathrm{cr}}]$.
Thus, for any $\eta\in (0,\epsilon_{\mathrm{cr}})$,
\begin{equation}\label{eq:epsilonsmall}
     \int_{\eta}^{\epsilon_{\mathrm{cr}}}\abs{\partial_{\epsilon}\pi_{\epsilon}(\Omega_1)}\, d \epsilon = - \int_{\eta}^{\epsilon_{\mathrm{cr}}}\partial_{\epsilon}\pi_{\epsilon}(\Omega_1)\, d \epsilon
     =\pi_{\eta}(\Omega_1)- \pi_{\epsilon_{\mathrm{cr}}}(\Omega_1)
     \leq 1
     .
\end{equation}
For $\epsilon\in [\epsilon_{\mathrm{cr}},1]$,  the function $\epsilon \mapsto \partial_{\epsilon}\pi_{\epsilon}(\Omega_1))$ is continuous and hence bounded.
This immediately implies~\eqref{eq:massBV}, concluding the proof of Case I.

\case[$U(x_{\min,1})<U(x_{\min,2})$]
Using~\eqref{eq:integral laplace method in one well with U} and \eqref{eq: est integral exp(-U/eta)} we see
\begin{align}
  \MoveEqLeft
\Big(\int_{\Omega_1}e^{-\frac{U}{\epsilon}}U\, d x\Big)\Big(\int_{\Omega_2}e^{-\frac{U}{\epsilon}}\, d x\Big)
-\Big(\int_{\Omega_2}e^{-\frac{U}{\epsilon}}U\, d x\Big)\Big(\int_{\Omega_1}e^{-\frac{U}{\epsilon}}\, d x\Big)\\
    &=(2\pi \epsilon)^d e^{-U(x_{\min,2})/\epsilon}
      \biggl(
      \begin{multlined}[t]
	\frac{O(\epsilon)}{\sqrt{\mathrm{det}(\nabla^2U(x_{\min,2}))}}
	\\
	  -\frac{(U(x_{\min,2})+O(\epsilon))}{\sqrt{\mathrm{det}(\nabla^2U(x_{\min,1}))}}(1+O(\epsilon))\biggr),
      \end{multlined}
\end{align}
which is negative when $\epsilon$ is small.
Using this in~\eqref{eq:partialPiEpsilon} implies $\partial_{\epsilon}\pi_{\epsilon}(\Omega_1) < 0$.
Using~\eqref{eq:epsilonsmall} and the same argument as in Case I.2 finishes the proof.
\end{proof}

\subsection{Uniform boundedness of \texorpdfstring{$r_k$}{rk}}

We recall that~$C_r$ defined by~\eqref{eq:defCr} is the maximum of the ratio of the \emph{normalized} densities.
Since this may be hard to estimate in practice, we now obtain a bound for~$C_r$ in a manner that may be easier to use in practice.
\begin{lemma}\label{lem: choose eta to make r bound}
  Suppose~$M$ is chosen by~\eqref{eq: parametersTNloc}, and choose~$\eta_2, \dots, \eta_M$ so that~$\eta_M = \eta$ and~$1/\eta_1$, \dots, $1/\eta_M$ are linearly spaced.
  If~$C_r = C_r( U/\eta_1, \nu )$ is defined by~\eqref{eq:defCr}, then~$C_r$ satisfies~\eqref{eq:CrExplicit}.
\end{lemma}

\begin{proof}[Proof of Lemma \ref{lem: choose eta to make r bound}]
Without loss of generality, we take $U=U_0$ , where $U_0$ is defined in \eqref{def:scU0}. Then we have $\eta_1=1$ and $U\geq 0$.

Observe that for every $k=1,\dots,M-1$, $x\in \mathcal{X}$, since $U\geq 0$, we have
\begin{equation}\label{e:rx}
     r_k(x) \overset{\mathclap{\eqref{e:rkExplicit}}}{=} \frac{Z_k}{Z_{k+1}}\exp\paren[\Big]{ -\paren[\Big]{\frac{1}{\eta_{k+1}} - \frac{1}{\eta_k} } U(x) }\leq \frac{Z_k}{Z_{k+1}}=\frac{\int_{\mathcal{X}} \exp\paren[\big]{\frac{-U}{\eta_{k}}} \, d y}{\int_{\mathcal{X}} \exp\paren[\big]{\frac{-U}{\eta_{k+1}}}\, d y}.
\end{equation}
   
Now we bound the ratio on the right hand side of \eqref{e:rx}. For $c\geq 0$, the constant $s_c$ defined in \eqref{def:scU0} now becomes
\begin{equation*}
    s_c \defeq
    \frac{\int_{\set{U > c}}e^{-U}\, d x}{\int_{\set{U\leq c}}e^{-U}\, d x}<\infty.
\end{equation*}

Then for $\epsilon<1$, 
 \begin{equation}
    \frac{\int_{\set{U > c}}e^{-\frac{U}{\epsilon}}\, d x}{\int_{\set{U \leq c}}e^{-\frac{U}{\epsilon}}\, d x}\leq\frac{\exp\paren[\big]{c\paren[\big]{1-\frac{1}{\epsilon}}}\int_{\set{U > c}}e^{-U}\, d x}{\exp\paren[\big]{c\paren[\big]{1-\frac{1}{\epsilon}}}\int_{\set{U \leq c}}e^{-U}\, d x}= \frac{\int_{\set{U > c}}e^{-U}\, d x}{\int_{\set{U \leq c}}e^{-U}\, d x}= s_c. \label{e:ratiosc}
 \end{equation}
Therefore, for $\eta_k<\eta_1$,
 \begin{align}
  \int_{\mathcal{X}} e^{-\frac{U}{\eta_{k}}}\, d x&= \int_{\set{U \leq c}} e^{-\frac{U}{\eta_{k}}}\, d x+\int_{\set{U > c}} e^{-\frac{U}{\eta_{k}}}\, d x\overset{\eqref{e:ratiosc}}{\leq} (1+s_c)\int_{\set{U \leq c}} e^{-\frac{U}{\eta_{k}}}\, d x\\
  &=(1+s_c)\int_{\set{U \leq c}} e^{-\frac{U}{\eta_{k+1}} + (\frac{U}{\eta_{k+1}}-\frac{U}{\eta_{k}})}\, d x\\
&\leq(1+s_c)\exp\paren[\Big]{ c\paren[\Big]{\frac{1}{\eta_{k+1}}-\frac{1}{\eta_{k}}}}\Big(\int_{\set{U \leq c}} e^{-\frac{U}{\eta_{k+1}}}\, d x \Big)\\
       &
       \leq (1+s_c)\exp\paren{ c\nu }\Big(\int_{\mathcal{X}} e^{-\frac{U}{\eta_{k+1}}}\, d x \Big). \label{e:intUetak}
\end{align}
   
  Here the last inequality is true because the choice of $M$ and $\eta_k$ (in~\eqref{e:MTN} and~\eqref{e:chooseEtaK} respectively) ensures
\begin{equation}\label{eq:etadiff}
   \frac{1}{\eta_{k+1}}-\frac{1}{\eta_{k}}\leq \nu
   .
\end{equation}
Since $r_k$ is always positive, using~\eqref{e:intUetak} in~\eqref{e:rx}, we obtain that for every $c>0$,
\begin{equation*}
    \sup_{1\leq k\leq M-1}\|r_k\|_{L^{\infty}}\leq (1+s_c)\exp\paren{ c\nu }.
\end{equation*}
Taking infimum on the right hand side 
  gives~\eqref{eq:defCr} with $C_r$ defined as in \eqref{eq:CrExplicit}.
\end{proof}

\subsection{Dimensional dependence of constants for separated energies}

\begin{proof}[Proof of Proposition \ref{prop:CrEg}]
  We only need to consider the case where $d\geq k_0$, where we recall~$k_0$ is the constant in~\eqref{e:AlmostPoly}.
  According to \eqref{e:CTCN} and \eqref{e:CBetaLoc}, it suffices to show that both $C_r$ and $C_{\LBV}$ are independent of $d$.
  For~$C_\LBV$, we note that~\eqref{e:sepPot} implies that we can take the domain $\Omega_j, j=1,\dots, J$ in the form
\begin{equation}
  \Omega_j=\tilde{\Omega}_j\times \mathbb{R}^{d-\tilde{d}},
\end{equation}
where $\tilde{\Omega}_j, j=1,\dots, J$ are subsets in $\mathbb{R}^{\tilde{d}}$, corresponding to the domain of measure proportional to $e^{-\tilde{U}_0}$.
  Then, Fubini's theorem shows that for any $j=1,\dots, J$ and any $\epsilon>0$, we have
  \begin{equation}
    \pi_{\epsilon}(\Omega_j)=
      \frac{\int_{\tilde{\Omega}_j} e^{-\tilde{U}_0/\epsilon}\, d x}
	{\int_{\mathbb{R}^{\tilde{d}}} e^{-\tilde{U}_0/\epsilon}\, d x}
      ,
  \end{equation}
  which is independent of $d$.
  Now using~\eqref{e:CBV} shows $C_{\LBV}$ is also independent of~$d$.
    
    Now it remains to find an upper bound of $C_r$ which is independent of $d$.
    Since $U_0$ is positive, we note
\begin{equation}
  \norm{r_k}_{L^{\infty}(\mathcal{X})} \overset{\eqref{e:defrk}}{\leq} T_1\cdot T_2,\quad\text{where}\quad T_1=\frac{\int_{\mathbb{R}^{\tilde{d}}}e^{-\frac{\tilde{U}_0}{\eta_{k}}}\, d x}{\int_{\mathbb{R}^{\tilde{d}}}e^{-\frac{\tilde{U}_0}{\eta_{k+1}}}\, d x},\quad T_2= \frac{\int_{\mathbb{R}^{d-\tilde{d}}}e^{-\frac{V_0}{\eta_{k}}}\, d x}{\int_{\mathbb{R}^{d-\tilde{d}}}e^{-\frac{V_0}{\eta_{k+1}}}\, d x}.
\end{equation}
Given our choice of $\eta_k$,  we bound $T_1$ by
\begin{equation}
  T_1 \overset{\eqref{e:intUetak}}{\leq}\inf_{c>0}(1+s_c(\tilde U_0))\exp\paren[\Big]{ \frac{c}{d} }
    \leq \inf_{c>0}(1+ s_c(\tilde U_0))e^c
    ,
\end{equation}
where
  \begin{equation}
    s_c(\tilde U_0) \defeq
  \frac{\int_{\set{\tilde{U}_0 > c}}e^{-\tilde{U}_0}\, d x}{\int_{\set{\tilde{U}_0\leq c}}e^{-\tilde{U}_0}\, d x}<\infty,
\end{equation}
  which implies that an upper bound of $T_1$ only depends on $\tilde{U}_0$ and is independent of~$d$.

Similarly for~$T_2$, we compute
\begin{equation*}
  T_2 \overset{\eqref{e:intUetak}}{\leq}\inf_{c>0}(1+s_c(V_0))\exp\paren[\Big]{ \frac{c}{d} },
  \quad\text{where}\quad s_c(V_0) \defeq
  \frac{\int_{\set{V_0 > c}}e^{-V_0}\, d x}{\int_{\set{V_0\leq c}}e^{-V_0}\, d x}
  .
\end{equation*}
Using \eqref{e:AlmostPoly} when $c>\alpha_{u}$, we compute
\begin{align}
  s_c(V_0)&\overset{\mathclap{\eqref{e:AlmostPoly}}}{\leq} e^{\alpha_{u}-\alpha_{b}}
  \frac{\int_{\set{V_0 > c}}e^{-\alpha_0|x-x_0|^{k_0}}\, d x}{\int_{\set{V_0\leq c}}e^{-\alpha_0|x-x_0|^{k_0}}\, d x}\\
  &\leq  e^{\alpha_{u}-\alpha_{b}}
  \frac{\int_{\set{\alpha_0|x-x_0|^{k_0}+\alpha_{u} > c}}e^{-\alpha_0|x-x_0|^{k_0}}\, d x}{\int_{\set{\alpha_0|x-x_0|^{k_0}+\alpha_{u}\leq c}}e^{-\alpha_0|x-x_0|^{k_0}}\, d x}\\
  &=e^{\alpha_{u}-\alpha_{b}}
  \frac{\int_{\set{|x|> \paren{\frac{c-\alpha_{u}}{\alpha_0}}^{1/k_0}}}e^{-\alpha_0|x|^{k_0}}\, d x}{\int_{\set{|x| \leq \paren{\frac{c-\alpha_{u}}{\alpha_0}}^{1/k_0}}}e^{-\alpha_0|x|^{k_0}}\, d x}\\
  &=\frac{e^{\alpha_u - \alpha_b} \Gamma(\frac{d}{k_0},c-\alpha_{u})}{\Gamma(\frac{d}{k_0})-\Gamma(\frac{d}{k_0},c-\alpha_{u})}\leq \frac{e^{\alpha_u - \alpha_b}}{\frac{\Gamma(\frac{d}{k_0})}{\Gamma(\frac{d}{k_0},c-\alpha_{u})}-1}. \label{eq:scubV0}
\end{align}
By the estimate of incomplete gamma function \cite[Satz 4.4.3]{Gabcke1979}, when $\frac{d}{k_0}\geq 1$ and $c-\alpha_{u}\geq \frac{d}{k_0}$, we have
\begin{equation*}
  \Gamma\paren[\Big]{\frac{d}{k_0},c-\alpha_{u}}\leq \frac{d}{k_0}\exp\paren[\big]{-(c-\alpha_{u})} (c-\alpha_u)^{\frac{d}{k_0}-1}.
\end{equation*}
On the other hand, Stirling's formula gives
\begin{equation}
  \Gamma\paren[\Big]{\frac{d}{k_0}}\geq C_{\Gamma}\paren[\Big]{\frac{d}{k_0}}^{\frac{d}{k_0}-\frac{1}{2}}\exp\paren[\Big]{-\frac{d}{k_0}}
  ,
\end{equation}
for a positive constant $C_{\Gamma}$.
We now choose $c-\alpha_{u}=\tilde{\nu} d/k_0$, where $\tilde{\nu}>1$ is such that
\begin{equation}\label{eq:tildenu}
    \tilde{\nu}-\log(\tilde{\nu})\geq \frac32+\log\paren[\Big]{\frac{2}{C_{\Gamma}}}
    .
\end{equation}
This gives
\begin{align}
    \frac{\Gamma(\frac{d}{k_0})}{\Gamma(\frac{d}{k_0},c-\alpha_{u})}&\geq \frac{C_{\Gamma}\paren[\big]{\frac{d}{k_0}}^{\frac{d}{k_0}-\frac{1}{2}}\exp\paren[\big]{-\frac{d}{k_0}}}{\frac{d}{k_0}\exp(-(c-\alpha_{u})) (c-\alpha_u)^{\frac{d}{k_0}-1}}\\
&\overset{\mathclap{z=d/k_0}}{=} \quad C_{\Gamma} \exp\paren[\Big]{(\tilde{\nu}-1)z-\frac12\log(z)-(z-1)\log(\tilde{\nu})}\\
    &\geq C_{\Gamma} \exp\paren[\Big]{\paren[\Big]{\tilde{\nu}-\log(\tilde{\nu})-\frac32}z}\overset{\eqref{eq:tildenu}}{\geq}2\label{e:GammaEst},
\end{align}
where the last inequality we use the fact that $z=\frac{d}{k_0}\geq 1$. Therefore, for this choice of $c$,
\begin{equation}
  s_c(V_0)\overset{\eqref{eq:scubV0},\eqref{e:GammaEst}}{\leq} e^{\alpha_u - \alpha_b},
\end{equation}
 which implies that
\begin{equation*}
  T_2\leq (1+s_{c}(V_0))\exp\paren[\Big]{ \frac{\tilde{\nu}}{k_0}+\frac{\alpha_{u}}{d} }\leq (1 + e^{\alpha_u - \alpha_b}) \exp\paren[\Big]{ \frac{\tilde{\nu}}{k_0}+\alpha_u },
\end{equation*}
which is independent of $d$.
\end{proof}

\section{Verification of properties~\ref{p:spectral},~\ref{p:var} and~\ref{property:eigenfunctions}}\label{sec:propertycheck}

In this section, we check Properties~\ref{p:spectral},~\ref{p:var} and ~\ref{property:eigenfunctions} under Assumptions~\ref{a:criticalpts},~\ref{assumption: nondegeneracy} and~\ref{a:massRatioBound}. 

\subsection{Lower bounds on eigenvalues (Property~\ref{p:spectral})} \label{sec:prop69}
In this section, we show that Property~\ref{p:spectral} holds for a double well potential.

\begin{proposition}\label{l: lower_bound_next_eigenvalue}
  If~$U$ satisfies Assumptions~\ref{a:criticalpts} and~\ref{assumption: nondegeneracy},
  then Property~\ref{p:spectral} holds.
\end{proposition}

We now prove Proposition~\ref{l: lower_bound_next_eigenvalue}.
Let~$\mathcal X$ be an Euclidean space, and~$\mu$ be a probability measure on~$\mathcal X$.
Recall~$\mu$ satisfies the \emph{Poincar\'e inequality} with constant~$\varrho$ (denoted by~$\PI(\varrho)$)
if for all test functions~$f \in H^1(\mu)$ we have
\begin{equation}\label{def: PI}\tag{$\PI(\varrho)$}
    \var_{\mu}(f)
    \leq \frac{1}{\varrho}\int_{\mathcal X} |\nabla f|^2\, d\mu
    .
\end{equation}
Here~$\var_{\mu}(f)$ is the variance of~$f$ with respect to the measure~$\mu$ and is defined by
\begin{equation}
  \var_{\mu}(f)\defeq \int_{\mathcal X}\left(f-\int_{\mathcal X} f\, d\mu\right)^2\, d\mu.
\end{equation}

The bound~\eqref{e:egvalLEpsilon} in Property~\ref{p:spectral} follows from a lower bound on the Poincar\'e constant proved in Menz and Schlichting~\cite{MenzSchlichting14}.
\begin{proposition}[Corollary 2.15 in \cite{MenzSchlichting14}]\label{prop: PI constant}
  If~$U$ satisfies Assumption~\ref{a:criticalpts} and~\ref{assumption: nondegeneracy} then $\pi_{\epsilon}$ satisfies~$\PI(\varrho_{\epsilon})$ with
    \begin{equation}\label{eq: PI_constant}
        \frac{1}{\varrho_{\epsilon}}\lesssim
	  \frac{\pi_{\epsilon}(\Omega_1)\pi_{\epsilon}(\Omega_2)}{(2\pi\epsilon)^{\frac{d}{2}-1}}
	  \frac{\sqrt{|\det(\nabla^2( U(s_{1,2})))|}}{|\lambda^{-}(s_{1,2})|}\exp\paren[\Big]{\frac{U(s_{1,2})}{\epsilon}}
	  \int_{\mathcal X} e^{-U/\epsilon} \, dx
	  ,
    \end{equation}
   where $\lambda^{-}(s_{1,2})$ denotes the negative eigenvalue of the Hessian $\nabla^2( U(s_{1,2}))$
  at the communicating saddle $s_{1,2}$.
\end{proposition}

We note that in~\cite{MenzSchlichting14} their domain is the whole space~$\R^d$.
The proof can easily be modified to work in the setting of the compact torus. We now prove Proposition~\ref{l: lower_bound_next_eigenvalue}.
\begin{proof}[Proof of Proposition~\ref{l: lower_bound_next_eigenvalue}]
The lower bound~\eqref{e:lambdaiLower} is well known and can, for instance, be found in~\cite[Proposition~2.1,Chapter~8]{Kolokoltsov00}.
  We now prove~\eqref{e:egvalLEpsilon} using Proposition~\ref{prop: PI constant}.
  Since
  \begin{equation}
    \int_{\mathcal X} \abs{\grad f}^2 \pi_\epsilon \, dx
      = \int_{\mathcal X} f L_\epsilon f \, \pi_\epsilon \, dx 
  \end{equation}
  we immediately see~$\lambda_{2, \epsilon} \geq \varrho_\epsilon$.
  Thus, Proposition \ref{prop: PI constant} implies
\begin{equation*}
    \limsup_{\epsilon\to 0}-(\epsilon\log(\lambda_{2,\epsilon}))\leq \limsup_{\epsilon\to 0}-(\epsilon\log(\varrho_{\epsilon}))\leq \hat{U}.
\end{equation*}

Thus, for every $H> \hat{U}\geq\limsup_{\epsilon\to 0}-(\epsilon\log(\lambda_{2,\epsilon}))$, there exists $\epsilon_{H}$ such that $\lambda_{2,\epsilon}\geq \exp(-\frac{H}{\epsilon})$ for every $\epsilon<\epsilon_{H}$.
Choosing
\begin{equation}
  C_{H}\defeq\min\set[\bigg]{\inf_{\epsilon_{H}\leq \epsilon \leq 1}\paren[\bigg]{\lambda_{2,\epsilon}\exp \paren[\Big]{\frac{H}{\epsilon}} },1}
\end{equation}
  immediately implies~\eqref{e:egvalLEpsilon} as desired.
\end{proof}

\subsection{Variation of eigenfunctions in temperature (Property~\ref{p:var})} \label{sec:prof6566}

We now prove Property~\ref{p:var} holds for potentials~$U$ that satisfy Assumptions~\ref{a:criticalpts} and~\ref{a:massRatioBound}.
\begin{proposition}\label{prop:pvar}
  If the potential~$U$ satisfies Assumptions~\ref{a:criticalpts} and~\ref{a:massRatioBound}, then Property~\ref{p:var} holds for the function~$g$ defined by
  \begin{equation}\label{e:gdef}
    g(\epsilon) \defeq \pi_\epsilon(\Omega_1)
    .
  \end{equation}
\end{proposition}

To prove~Proposition~\ref{prop:pvar}, we need a few preliminary lemmas.
Since it will be useful later, and it does not require too much additional effort,
we directly state and prove these lemmas allowing for the more general case where~$U$ has~$J$ wells, with~$J \geq 2$.

We begin by quantifying the fact that when $\epsilon$ is small, the eigenfunctions~$\psi_{2,\epsilon}$, \dots, $\psi_{J,\epsilon}$ are very close to linear combinations of the indicator functions~$\one_{\Omega_j}$, for $j=1,\dots,J$.
To state a precise bound, consider  the subspaces~$E_\epsilon, F_\epsilon \subseteq L^2(\pi_\epsilon)$ defined by
\begin{equation}\label{eq:defEepsFeps}
   E_{\epsilon}\defeq\operatorname{span}\{\one_{\Omega_1},\dots,\one_{\Omega_J}\},
   \quad F_{\epsilon}\defeq\operatorname{span}\{1,\psi_{2,\epsilon},\dots,\psi_{J,\epsilon}\}
    \subseteq L^2(\pi_\epsilon)
   ,
\end{equation}
We measure closeness of each~$\psi_{i, \epsilon}$ to a linear combination of $\set{\one_{\Omega_1}}$, by measuring the ``distance'' between the subspaces~$E_\epsilon$ and~$F_\epsilon$.
Explicitly, define
\begin{equation}
   d(E_{\epsilon}, F_{\epsilon})\defeq \norm{P_{E_{\epsilon}}-P_{E_{\epsilon}}P_{F_{\epsilon}}}=\norm{P_{E_{\epsilon}}-P_{F_{\epsilon}}P_{E_{\epsilon}}}.
\end{equation}
Here $P_{E_{\epsilon}}$, $P_{E_{\epsilon}}$ are the~$L^2(\pi_\epsilon)$ orthonormal projectors onto $E_{\epsilon}$ and $F_{\epsilon}$ respectively,
and~$\norm{\cdot}$ denotes the operator norm on~$L^2(\pi_\epsilon)$.
Proposition 2.2 in Chapter~8 of~\cite{Kolokoltsov00} gives an estimate on~$d(E_\epsilon, F_\epsilon)$, and we reproduce the result here for easy reference.

\begin{proposition}[Chapter 8,  Proposition 2.2 of \cite{Kolokoltsov00}]\label{prop_kolo_chap8_prop2.2}
Let~$\hat{\gamma}$ be the smallest well-escape energy defined by
\begin{equation}\label{e:gammaHatDefM}
  \hat{\gamma}\defeq \min_{1 \leq i \leq J} \Delta_i,
  \quad\text{where}\quad
  \Delta_i \defeq \min_{x\in\partial\Omega_i}U(x)- U(x_{\min,i}).
\end{equation}
  For any~$\gamma<\hat{\gamma}$,
  there exists a constant $ C_{\gamma}>0$ such that for all~$\epsilon\leq 1$, we have
\begin{equation}\label{e:dEF}
    d(E_{\epsilon}, F_{\epsilon})\leq  C_{\gamma}\exp
\Big(\frac{-\gamma}{\epsilon}\Big).
\end{equation}
\end{proposition}

Note that when~$J=2$, the quantity~$\hat \gamma$ defined in~\eqref{e:gammaHatDefM} is precisely the energy barrier defined in~\eqref{e:gammaHatDef}.
To prove  Proposition~\ref{prop:pvar} we need two more lemmas concerning the ratio of the densities~$\pi_\epsilon / \pi_{\epsilon'}$.
 
\begin{lemma}\label{lem: integral of two eigenfunction from different levels}
Let $\epsilon'<\epsilon$ and define~$r_\epsilon$ by
  \begin{equation}\label{e:rdefEps}
    r_{\epsilon}\defeq\frac{\pi_{\epsilon'}}{\pi_{\epsilon}}
    \,.
  \end{equation}
  Then, for each $2\leq i\leq J$,
\begin{equation}
    \norm{P_{F_{\epsilon}}(\psi_{i,\epsilon'}r_{\epsilon})}_{L^2(\pi_{\epsilon})}
     \leq 1+\frac12 C_m^2\sum_{j=1}^{J}\abs[\big]{\pi_{\epsilon'}(\Omega_j)-\pi_{\epsilon}(\Omega_j)}+\|r_{\epsilon}\|^{\frac12}_{L^{\infty}}d(E_{\epsilon},F_{\epsilon}).
  \label{eq:intPsiPsiPi}
\end{equation}
\end{lemma}

\begin{lemma}\label{lem: integral of eigenfunction on different measure}
    Let $\epsilon'<\epsilon$. Then
    \begin{equation}
      \norm{P_{F_{\epsilon}}(r_{\epsilon}-1)}_{L^2(\pi_{\epsilon})}\leq
    C_m\sum_{j=1}^{J}\abs[\big]{\pi_{\epsilon'}(\Omega_j)-\pi_{\epsilon}(\Omega_j)} + \|r_{\epsilon}\|^{\frac12}_{L^{\infty}}d(E_{\epsilon},F_{\epsilon}).  \label{eq:intPsiPi}
    \end{equation}
\end{lemma}

Momentarily postponing the proof of Lemmas~\ref{lem: integral of two eigenfunction from different levels} and~\ref{lem: integral of eigenfunction on different measure}, we now prove Proposition~\ref{prop:pvar}.

\begin{proof}[Proof of Proposition~\ref{prop:pvar}]
  Recall Lemma~\ref{lem: mass in well as BV} already shows that the function~$g$ defined in~\eqref{e:gdef} has bounded variation.
  Thus, to prove Proposition~\ref{prop:pvar} we only need to show that~\eqref{e:PsiEpPi} and~\eqref{e:PsiEpEpP} hold for the function~$g$ defined by~\eqref{e:gdef}.

We first prove~\eqref{e:PsiEpEpP}. In the double-well case, we have
\begin{equation}
    F_{\epsilon}=\operatorname{span}\{1,\psi_{2,\epsilon}\}.
\end{equation}
As a result, for any $0<\epsilon'<\epsilon\leq 1$,
\begin{align*}
    P_{F_{\epsilon}}(\psi_{2,\epsilon'}r_{\epsilon})&=\paren[\Big]{\int_{\T^d}\psi_{2,\epsilon'}r_{\epsilon}\pi_{\epsilon}\, d x}+\paren[\Big]{\int_{\T^d}\psi_{2,\epsilon'}r_{\epsilon}\psi_{2,\epsilon}\pi_{\epsilon}\, d x}\psi_{2,\epsilon}\\
    &\overset{\mathclap{\eqref{e:rdefEps}}}{=}\paren[\Big]{\int_{\T^d}\psi_{2,\epsilon'}\psi_{2,\epsilon}\pi_{\epsilon'}\, d x}\psi_{2,\epsilon}.
\end{align*}
The above equation immediately shows that
\begin{align}
    \abs[\Big]{\int_{\T^d}\psi_{2,\epsilon'}\psi_{2,\epsilon}\pi_{\epsilon'}\, d x}
    &=
   \norm{P_{F_{\epsilon}}(\psi_{2,\epsilon'}r_{\epsilon})}_{L^2(\pi_{\epsilon})}
   \\
   &\overset{\mathclap{\eqref{eq:intPsiPsiPi}}}{\leq} ~~ 1+\frac12 C_m^2\sum_{j=1}^{2}\abs[\big]{\pi_{\epsilon'}(\Omega_j)-\pi_{\epsilon}(\Omega_j)}+\|r_{\epsilon}\|^{\frac12}_{L^{\infty}}d(E_{\epsilon},F_{\epsilon})\\
    &\overset{\mathclap{\eqref{eq:defCr},\eqref{e:dEF}}}{\leq} ~~ 1+ C_m^2\abs[\big]{\pi_{\epsilon'}(\Omega_1)-\pi_{\epsilon}(\Omega_1)}+C_r^{\frac12}C_{\gamma}\exp
\Big(\frac{-\gamma}{\epsilon}\Big)
\,.
\end{align}
This implies that~\eqref{e:PsiEpEpP} in Property~\ref{p:var} holds for the function~$g(\epsilon)=\pi_{\epsilon}(\Omega_1)$.
Note that Lemma~\ref{lem: mass in well as BV} guarantees that the function~$g$ is of bounded variation.
\smallskip

We now prove~\eqref{e:PsiEpPi}.
  Using Lemma~\ref{lem: integral of eigenfunction on different measure}, and a similar argument shows
\begin{align*}
  \abs[\bigg]{\int_{\T^d}\psi_{2,\epsilon}\pi_{\epsilon'}\, d x}
    &=
   \norm{P_{F_{\epsilon}}(r_{\epsilon}-1)}_{L^2(\pi_{\epsilon})}
    \\
   &\overset{\mathclap{\eqref{eq:defCr},\eqref{e:dEF},\eqref{eq:intPsiPi}}}{\leq} \qquad 2C_m\abs[\big]{\pi_{\epsilon'}(\Omega_1)-\pi_{\epsilon}(\Omega_1)}+C_r^{\frac12}C_{\gamma}\exp
\Big(\frac{-\gamma}{\epsilon}\Big),
\end{align*}
  which implies~\eqref{e:PsiEpPi}, and concludes the proof.
\end{proof}

We now prove Lemmas~\ref{lem: integral of two eigenfunction from different levels} and~\ref{lem: integral of eigenfunction on different measure}.
\begin{proof}[Proof of Lemma \ref{lem: integral of two eigenfunction from different levels}]
Fix~$i$, where~$2\leq i\leq J$, and define
\begin{equation}
    \psi=\frac{P_{F_{\epsilon}}(\psi_{i,\epsilon'}r_{\epsilon})}{\norm{P_{F_{\epsilon}}(\psi_{i,\epsilon'}r_{\epsilon})}_{L^2(\pi_{\epsilon})}}
    \,.
\end{equation}
  Clearly~$\norm{\psi}_{L^2(\pi_{\epsilon})}=1$.
  Moreover, since~$P_{F_\epsilon}$ is an orthogonal projection, we see
\begin{equation}
    \abs[\bigg]{\int_{\T^d}(\psi_{i,\epsilon'}r_{\epsilon})\psi\pi_{\epsilon}\, d x}=\norm{P_{F_{\epsilon}}(\psi_{i,\epsilon'}r_{\epsilon})}_{L^2(\pi_{\epsilon})}
    .
\end{equation}
Using the Cauchy-Schwarz inequality and recalling~$r_\epsilon$ is defined by~\eqref{e:rdefEps}, we obtain
\begin{align}
    \bigg|\int_{\mathbb T^d}\psi_{i,\epsilon'}r_{\epsilon}\psi\pi_{\epsilon}\, d x\bigg|
    &= \bigg|\int_{\mathbb T^d}\psi_{i,\epsilon'}\psi\pi_{\epsilon'}\, d x\bigg|
    \\
    \label{eq:leqpsil2}
    &\leq
      \norm{\psi}_{L^2(\pi_{\epsilon'})}
      \norm{\psi_{i, \epsilon'}}_{L^2(\pi_{\epsilon'})}
    = \|\psi\|_{L^2(\pi_{\epsilon'})}.
\end{align}

It remains to compute $\|\psi\|_{L^2(\pi_{\epsilon'})}$.
For this, we decompose $\psi$ into the sum of the projection into $E_{\epsilon}$ and $E^{\perp}_{\epsilon}$, where $E_{\epsilon}$ is defined in \eqref{eq:defEepsFeps}. Explicitly,
\begin{equation}\label{eq:decompEepsilon}
    \psi=\sum_{j=1}^{J}a_{j,\epsilon}\one_{\Omega_j}+\upsilon_{\epsilon},
\end{equation}
where $\upsilon_{\epsilon}\in E^{\perp}_{\epsilon}$. Thus we can bound $\|\psi\|_{L^2(\pi_{\epsilon'})}$ by
\begin{equation}\label{eq:psidecompEproj}
    \|\psi\|_{L^2(\pi_{\epsilon'})}\leq \norm[\bigg]{\sum_{j=1}^{J}a_{j,\epsilon}\one_{\Omega_j}}_{L^2(\pi_{\epsilon'})}+\|\upsilon_{\epsilon}\|_{L^2(\pi_{\epsilon'})}.
\end{equation}

Since $\upsilon_{\epsilon}$ is orthogonal to each $\one_{\Omega_j}$ in $L^2(\pi_{\epsilon})$, we have
\begin{equation}\label{eq:akeq1}
   \sum_{j=1}^{J}\paren{a_{j,\epsilon}}^2\pi_{\epsilon}(\Omega_j) \leq \norm{\psi}^2_{L^2(\pi_{\epsilon})}=1.
\end{equation}
Notice that~\eqref{eq:akeq1}  implies that for every $a^{(i)}_{j,\epsilon}$, according to Assumption~\ref{a:massRatioBound},
\begin{equation}\label{eq:akub}
    \paren{a_{j,\epsilon}}^2\overset{\eqref{eq:akeq1}}{\leq} \frac{1}{\pi_{\epsilon}(\Omega_j)
    }\overset{\eqref{e:massRatioBound}}{\leq} C^2_m.
\end{equation}

Now we compute for $i=2,\dots,J$,
\begin{align}
        \norm[\bigg]{\sum_{j=1}^{J}a_{j,\epsilon}\one_{\Omega_j}}^2_{L^2(\pi_{\epsilon'})}
        &= \sum_{j=1}^{J}\paren{a_{j,\epsilon}}^2\pi_{\epsilon'}(\Omega_j)\\
        &=\sum_{j=1}^{J}\paren{a_{j,\epsilon}}^2\paren[\big]{\pi_{\epsilon'}(\Omega_j)-\pi_{\epsilon}(\Omega_j)}+ \sum_{j=1}^{J}\paren{a_{j,\epsilon}}^2\pi_{\epsilon}(\Omega_j)\\
        &\overset{\mathclap{\eqref{eq:akeq1},\eqref{eq:akub}}}{\leq }\qquad 1+C_m^2\sum_{j=1}^{J}\abs[\big]{\pi_{\epsilon'}(\Omega_j)-\pi_{\epsilon}(\Omega_j)}.
\end{align}
Using the inequality $(1+y)^{\frac12}\leq 1+\frac12 y$ for $y>0$ yields
\begin{equation}\label{eq:Eproj}
    \norm[\bigg]{\sum_{j=1}^{J}a_{j,\epsilon}\one_{\Omega_j}}_{L^2(\pi_{\epsilon'})}\leq 1+\frac{C_m^2}{2}\sum_{j=1}^{J}\abs[\big]{\pi_{\epsilon'}(\Omega_j)-\pi_{\epsilon}(\Omega_j)}.
\end{equation}
On the other hand, using the fact that~$\psi\in F_{\epsilon}$, we have that
\begin{align}
\|\upsilon_{\epsilon}\|_{L^2(\pi_{\epsilon})}&=\|(I-P_{E_{\epsilon}})\psi\|_{L^2(\pi_{\epsilon})}\\
       &=\|(P_{F_{\epsilon}}-P_{E_{\epsilon}}P_{F_{\epsilon}})\psi\|_{L^2(\pi_{\epsilon})}\leq \|P_{F_{\epsilon}}-P_{E_{\epsilon}}P_{F_{\epsilon}}\|
       = d(E_{\epsilon},F_{\epsilon}), \label{eqVkLessThanDEkFk} 
\end{align}
where the last equality is based on the fact that~$d(F_{\epsilon},E_{\epsilon})=d(E_{\epsilon},F_{\epsilon})$ since~$E_{\epsilon}$ and~$F_{\epsilon}$ are finite dimensional subspaces with equal dimensions~\cite[Page 240]{Kolokoltsov00}. Hence we obtain
\begin{align}
  \|\upsilon_{\epsilon}\|_{L^2(\pi_{\epsilon'})}
    &= \paren[\Big]{\int_{\mathbb T^d}(\upsilon_{\epsilon})^2r_{\epsilon}\pi_{\epsilon}\, d x}^{\frac12}
  \\
  \label{eq:upsEst}
  &\leq \|r_{\epsilon}\|^{\frac12}_{L^{\infty}}\|\upsilon_{\epsilon}\|_{L^2(\pi_{\epsilon})}\overset{\eqref{eqVkLessThanDEkFk}}{\leq} \|r_{\epsilon}\|^{\frac12}_{L^{\infty}}d(E_{\epsilon},F_{\epsilon}).
\end{align}

Therefore,
\begin{align}
  \MoveEqLeft   \norm{P_{F_{\epsilon}}(\psi_{i,\epsilon'}r_{\epsilon})}_{L^2(\pi_{\epsilon})}
\overset{\eqref{eq:leqpsil2},\eqref{eq:psidecompEproj}}{\leq} \norm[\bigg]{\sum_{j=1}^{J}a_{j,\epsilon}\one_{\Omega_j}}_{L^2(\pi_{\epsilon'})}+\|\upsilon_{\epsilon}\|_{L^2(\pi_{\epsilon'})}
  \\
&\overset{\eqref{eq:Eproj},\eqref{eq:upsEst}}{\leq}
  1+\frac12 C_m^2\sum_{j=1}^{J}\abs[\Big]{\pi_{\epsilon'}(\Omega_j)-\pi_{\epsilon}(\Omega_j)}+\|r_{\epsilon}\|^{\frac12}_{L^{\infty}}d(E_{\epsilon},F_{\epsilon})
  ,
\end{align}
which yields~\eqref{eq:intPsiPsiPi}, concluding the proof.
\end{proof}

\begin{proof}[Proof of Lemma \ref{lem: integral of eigenfunction on different measure}]
  The proof is similar to the proof of Lemma~\ref{lem: integral of two eigenfunction from different levels}.
  Define~$\tilde \psi$ by
  \begin{equation}
    \tilde{\psi}= \frac{P_{F_{\epsilon}}(r_{\epsilon}-1)}{\norm{P_{F_{\epsilon}}(r_{\epsilon}-1)}_{L^2(\pi_{\epsilon})}}
    .
  \end{equation}
  Clearly~$\norm{\tilde \psi}_{L^2(\pi_\epsilon)} = 1$ and
\begin{equation}\label{eq:psirk}
    \abs[\bigg]{\int_{\T^d}(r_{\epsilon}-1)\tilde{\psi}\pi_{\epsilon}\, d x}=\norm{P_{F_{\epsilon}}(r_{\epsilon}-1)}_{L^2(\pi_{\epsilon})}
\end{equation}
Moreover, since
  \begin{equation}
    \int_{\T^d}P_{F_\epsilon}(r_{\epsilon}-1) \, \pi_{\epsilon}\, d x
    = \int_{\T^d}(r_{\epsilon}-1)\pi_{\epsilon}\, d x
      =0
  \end{equation}
  we see that
\begin{equation}\label{eq:intpsi0}
    \int_{\T^d}\tilde{\psi} \pi_{\epsilon}=0.
\end{equation}
Using a decomposition similar to~\eqref{eq:decompEepsilon}, we see
\begin{equation}\label{eq:decomptildePsi}
    \tilde{\psi}=\sum_{j=1}^{J}\tilde{a}_{j,\epsilon}\one_{\Omega_j}+\tilde{\upsilon}_{\epsilon},\quad\text{where}\quad \tilde\upsilon_{\epsilon}\in E^{\perp}_{\epsilon}.
\end{equation}
Combined with~\eqref{eq:intpsi0} this implies
\begin{equation}\label{eq:intpsi0a}
    \sum_{j=1}^{J}\tilde{a}_{j,\epsilon}\pi_{\epsilon}(\Omega_j)=0.
\end{equation}
Thus we have
\begin{align}
    \abs[\bigg]{\int_{\T^d}(r_{\epsilon}-1)\tilde{\psi}\pi_{\epsilon}\, d x}~~&\overset{\mathclap{\eqref{eq:intpsi0}}}{=}\quad \abs[\bigg]{\int_{\T^d}r_{\epsilon}\tilde{\psi}\pi_{\epsilon}\, d x}\overset{\eqref{eq:decomptildePsi},\eqref{e:rdefEps}}{=}\abs[\bigg]{\int_{\T^d}\paren[\Big]{\sum_{j=1}^{J}\tilde{a}_{j,\epsilon}\one_{\Omega_j}+\tilde{\upsilon}_{\epsilon}}\pi_{\epsilon'}\, d x}\\
    &\leq \abs[\bigg]{\sum_{j=1}^{J}\tilde{a}_{j,\epsilon}\pi_{\epsilon'}(\Omega_j)}+\norm{\tilde{v}_{\epsilon}}_{L^2(\pi_{\epsilon'})}\\
    &\overset{\mathclap{\eqref{eq:intpsi0a}}}{=}\quad\abs[\bigg]{\sum_{j=1}^{J}\tilde{a}_{j,\epsilon}(\pi_{\epsilon'}(\Omega_j)-\pi_{\epsilon}(\Omega_j))}+\norm{\tilde{v}_{\epsilon}}_{L^2(\pi_{\epsilon'})}\\
    &\overset{\mathclap{\eqref{eq:akub},\eqref{eq:upsEst}}}{\leq}\qquad C_m\sum_{j=1}^{J}\abs[\big]{\pi_{\epsilon'}(\Omega_j)-\pi_{\epsilon}(\Omega_j)} + d(E_{\epsilon},F_{\epsilon})\|r_{\epsilon}\|^{\frac12}_{L^{\infty}(\pi_{\epsilon})}
    ,
\end{align}
  which yields~\eqref{eq:intPsiPi}, concluding the proof.
\end{proof}

\subsection{Uniform boundedness of eigenfunctions (Property~\ref{property:eigenfunctions})}\label{sec:Cpsi}

Finally we prove Property~\ref{p:var} holds for potentials~$U$ that satisfy Assumptions~\ref{a:criticalpts} and~\ref{a:massRatioBound}.
\begin{proposition}\label{lem: uniform boundedness of eigenfunction} 
    If the potential~$U$ satisfies Assumptions~\ref{a:criticalpts} and~\ref{a:massRatioBound}, then Property~\ref{property:eigenfunctions} holds.
\end{proposition}

We prove Proposition~\ref{lem: uniform boundedness of eigenfunction} in two steps.
First we find small neighborhoods of the local minima and use a local maximum principle to show that $\psi_{2,\epsilon}$ is uniformly bounded in $\epsilon$ in these neighborhoods.
Then we apply global maximum principle to show the uniform boundedness outside these regions.
This proof strategy also works when~$U$ has~$J\geq 2$ local minima. 

Define the $R_i>0$, $i\in\set{1,2}$ by
\begin{equation}\label{eq: def of R1 R2}
  R_i\defeq \sup\set[\Big]{r \st B(x_{\min,i},r)\subseteq\Omega_i, \sup_{x\in B(x_{\min,i},r)}U(x)-U(x_{\min,i})\leq \frac{\hat \gamma}{8}} 
\end{equation}
and then define
\begin{equation}
   \widetilde{B}_i=B(x_{\min,i},R_i),\quad 
    B_i=B\paren[\Big]{x_{\min,i},\frac{3 R_i}{4}}.
\end{equation}
We will show $\psi_{2,\epsilon}$ is uniformly bounded in $\epsilon$ both on $B_1\cup B_2$ and $\T^d\setminus (B_1\cup B_2)$.
We first bound~$\psi_{2, \epsilon}$ in the regions $B_1$ and $B_2$. 
\begin{lemma}\label{lem: boundedness of eigenfunction around local minima}
  There exists  a constant~$C_a= C_a(d,U,C_m)$ and $\tilde{\epsilon}=\tilde{\epsilon}(d,U)$ such that for every
  \begin{equation}\label{eq:epsilonCri}
      0 < \epsilon\leq \min\set[\Big]{\frac{1}{12}\min\{R_1,R_2\},\tilde{\epsilon}}
  \end{equation}
  we have
\begin{equation}\label{eq: eigenvector in Bi}
    \forall x\in B_i, \quad  |\psi_{2,\epsilon}(x)-a_{i,\epsilon}|\leq C_a
      \exp\paren[\Big]{-\frac{ \hat{\gamma}}{4 \epsilon}}, \quad i=1,2.
\end{equation}
 Here   $a_{1,\epsilon}$ and $a_{2,\epsilon}$ are defined by
\begin{equation}\label{eq:aksol}
     a_{i,\epsilon}\defeq \frac{1}{\pi_{\epsilon}(\Omega_i)}\int_{\T^d}\one_{\Omega_i}\psi_{2,\epsilon}\pi_{\epsilon}\, d x,\quad i=1,2.
 \end{equation}
\end{lemma}

Before we prove Lemma~\ref{lem: boundedness of eigenfunction around local minima}, we introduce two preliminary estimates.
The first one is an upper bound on the second eigenvalue~$\lambda_{2,\epsilon}$.
It is well known that the second eigenvalue~$\lambda_{2,\epsilon}$ is exponentially small in~$\epsilon$.
Explicitly, in the proof of Proposition~2.2 in Chapter~8 of~\cite{Kolokoltsov00} we see that there exists a constant~$\tilde{C}$, such that for every $\epsilon\in (0,1]$,
\begin{equation}\label{eq:egvalub}
    \lambda_{2,\epsilon}\leq \tilde{C}\exp\paren[\Big]{-\frac{\hat\gamma}{4\epsilon}}.
\end{equation} 

Next, we notice that~\eqref{eq:aksol} implies
\begin{equation}
    a_{1,\epsilon}\one_{\Omega_1}+a_{2,\epsilon}\one_{\Omega_2}=\P_{E_{\epsilon}}(\psi_{2,\epsilon}).
\end{equation}
Therefore,
\begin{align}
    \MoveEqLeft\int_{\T^d}\abs{\psi_{2,\epsilon}-a_{1,\epsilon}\one_{\Omega_1}-a_{2,\epsilon}\one_{\Omega_2}}^2\pi_{\epsilon}\, d x = \norm{\psi_{2,\epsilon}-\P_{E_{\epsilon}}(\psi_{2,\epsilon})}^2_{L^2(\pi_\epsilon)}\\
    &= \norm{\P_{F_{\epsilon}}(\psi_{2,\epsilon})-\P_{E_{\epsilon}}(\P_{F_{\epsilon}}(\psi_{2,\epsilon}))}^2_{L^2(\pi_\epsilon)}\\
    &\leq \norm{\P_{F_{\epsilon}}-\P_{E_{\epsilon}}\P_{F_{\epsilon}}}^2=d(E_{\epsilon},F_{\epsilon})^2. \label{eqpsiklemmaCpsi}
\end{align}

We are now equipped to prove Lemma~\ref{lem: boundedness of eigenfunction around local minima}. In the proof, the constant $C= C(U,d, C_m)$ may change from line to line.

\begin{proof}[Proof of Lemma \ref{lem: boundedness of eigenfunction around local minima}]

Fix $i=1$ or $2$,
for each $x\in B_i$, there exists $y\in B_i$ such that $x\in B(y,\epsilon)$.
  By the triangle inequality, we note~$B(y,2\epsilon)\subseteq \widetilde{B}_i$.
  Thus,
\begin{equation}\label{eq:concirc}
    B(y,2\epsilon)\subseteq \widetilde{B}_i\subset\Omega_i.
\end{equation}

First notice that the function $\psi_{2,\epsilon}-a_{i,\epsilon}$ satisfies
\[
(L_{\epsilon}-\lambda_{2,\epsilon})(\psi_{2,\epsilon}-a_{i,\epsilon})=\lambda_{2,\epsilon}a_{i,\epsilon}.
\]
Thus, for the range of~$\epsilon$ we consider, we may use~\cite[Corollary 9.21]{GT} to obtain the existence of a dimensional constant~$C$ such that for every $y\in B_i$ for which $B(y,2\epsilon)\subseteq \widetilde{B}_i$, we have
\begin{align}
    \sup_{x\in B(y,\epsilon)} \abs{\psi_{2,\epsilon}(x)-a_{i,\epsilon}}
	\leq C\biggl(&\paren[\Big]{\frac{1}{|B(y,2\epsilon)|}\int_{B(y,2\epsilon)}|\psi_{2,\epsilon}(x)-a_{i,\epsilon}|^2 \, d x}^{\frac12}
    \\\label{eq:GTlocal}
	& +|\lambda_{2,\epsilon}a_{i,\epsilon}|\biggr).
\end{align}

Now we bound $\int_{B(y,2\epsilon)}|\psi_{2,\epsilon}(x)-a_{i,\epsilon}|^2 \, d x$ that appears on the right hand side of~\eqref{eq:GTlocal}. Using the fact that when \eqref{eq:epsilonCri} holds, for $i=1,2$,
\begin{align}
     \int_{\mathbb{T}^d}e^{-\frac{U}{\epsilon}}\, d x&=\Big(\int_{\Omega_i}e^{-\frac{U}{\epsilon}}\, d x\Big)\cdot\Big( 1+\frac{\int_{\T^d\setminus\Omega_i}e^{-\frac{U}{\epsilon}}\, d x}{\int_{\Omega_i}e^{-\frac{U}{\epsilon}}\, d x}\Big)\\
     &=\Big(\int_{\Omega_i}e^{-\frac{U}{\epsilon}}\, d x\Big)\cdot\Big( 1+\frac{1-\pi_{\epsilon}(\Omega_i)}{\pi_{\epsilon}(\Omega_i)}\Big)\overset{\eqref{e:massRatioBound}}{\leq} \Big(\int_{\Omega_i}e^{-\frac{U}{\epsilon}}\, d x\Big)\cdot( 1+C_m^2)\\
    \label{eq:intoi} &\overset{\mathclap{\eqref{eq: est integral exp(-U/eta)}}}{\leq} C(2\pi\epsilon)^{\frac{d}{2}} e^{-\frac{U(x_{\min,i})}{\epsilon}},
\end{align}
we have that 
\begin{align}
   \MoveEqLeft \int_{B(y,2\epsilon)}|\psi_{2,\epsilon}-a_{i,\epsilon}|^2 \, d x= \int_{B(y,2\epsilon)}|\psi_{2,\epsilon}-a_{i,\epsilon}\one_{\Omega_i}|^2 \, d x \\
        &\leq  \Big(\sup_{z\in B(y,2\epsilon)}e^{\frac{U(z)}{\epsilon}}\Big)\int_{B(y,2\epsilon)}|\psi_{2,\epsilon}-a_{i,\epsilon}\one_{\Omega_i}|^2 e^{-\frac{U}{\epsilon}}\, d x\\
        &\overset{\mathclap{\eqref{eq:concirc}}}{\leq}\Big(\int_{\mathbb{T}^d}e^{-\frac{U}{\epsilon}}\, d x\Big)\Big(\sup_{z\in \widetilde{B}_i}e^{\frac{U(z)}{\epsilon}}\Big)\int_{\Omega_i}|\psi_{2,\epsilon}(x)-a_{i,\epsilon}\one_{\Omega_i}|^2 \, d \pi_{\epsilon}(x)\\
       &\overset{\mathclap{\eqref{eq:intoi}}}{\leq} C(2\pi\epsilon)^{\frac{d}{2}}\Big(\sup_{z\in \widetilde{B}_i}e^{\frac{U(z)-U(x_{\min,i})}{\epsilon}}\Big)\norm[\big]{\psi_{2,\epsilon}-a_{1,\epsilon}\one_{\Omega_1}-a_{2,\epsilon}\one_{\Omega_2}}_{L^2(\pi_{\epsilon})}^2\\
       &\overset{\mathclap{\eqref{eqpsiklemmaCpsi}}}{\leq}C(2\pi\epsilon)^{\frac{d}{2}}\Big(\sup_{z\in \widetilde{B}_i}e^{\frac{U(z)-U(x_{\min,i})}{\epsilon}}\Big)d(E_{\epsilon},F_{\epsilon})^2\\
       &\overset{\mathclap{\eqref{e:dEF}}}{\leq}C(2\pi\epsilon)^{\frac{d}{2}}\Big(\sup_{z\in \widetilde{B}_i}e^{\frac{U(z)-U(x_{\min,i})}{\epsilon}}\Big)\exp\paren[\Big]{-\frac{7\hat{\gamma}}{4\epsilon}}\\
       &\overset{\mathclap{\eqref{eq: def of R1 R2}}}{\leq}C(2\pi\epsilon)^{\frac{d}{2}}\exp\paren[\Big]{-\frac{13\hat{\gamma}}{8 \epsilon}},  \label{eq:13over8hatgamma}
\end{align}
where the second last inequality we use \eqref{e:dEF} with $\gamma =\frac{7}{8}\hat{\gamma}$.

Notice that there exists constant $\tilde{\epsilon}=\tilde{\epsilon}(d,\hat\gamma)$ that whenever $\epsilon<\tilde{\epsilon}$, 
\begin{equation}\label{eq:tildeEps}
    \exp\paren[\Big]{-\frac{\hat{\gamma}}{8 \epsilon}}< (2\pi\epsilon)^{\frac{d}{2}}.
\end{equation}
Thus, for $\epsilon<\tilde{\epsilon}$,
\begin{equation} \label{eq:nearBi}
    \int_{B(y,2\epsilon)}|\psi_{2,\epsilon}-a_{i,\epsilon}|^2 \, d x  \overset{\eqref{eq:13over8hatgamma},\eqref{eq:tildeEps}}{\leq} C(2\pi\epsilon)^{d}\exp\paren[\Big]{-\frac{3\hat{\gamma}}{2 \epsilon}}. 
\end{equation}

Therefore, plugging \eqref{eq:nearBi} into \eqref{eq:GTlocal} gives
\begin{align*}
      \sup_{x\in B(y,\epsilon)} |\psi_{2,\epsilon}(x)-a_{1,\epsilon}|
	~&\overset{\mathclap{\eqref{eq:nearBi}}}{\leq}~
	    \paren[\Big]{
	      \frac{C(2\pi\epsilon)^{d}}{|B(y,2\epsilon)|}\exp\paren[\Big]{-\frac{3\hat{\gamma}}{2\epsilon}}
	    }^{\frac12}
	  +C|\lambda_{2,\epsilon}a_{1,\epsilon}|
      \\
       &\overset{\mathclap{\eqref{eq:egvalub},\eqref{eq:nearBi}}}{\leq}\quad
	  \paren[\Big]{\frac{C(2\pi\epsilon)^{d}}{( 2\epsilon)^d}\exp\paren[\Big]{-\frac{3\hat{\gamma}}{2\epsilon}}}^{\frac12}
        +C \tilde{C}C_m\exp\paren[\Big]{-\frac{\hat{\gamma}}{4\epsilon}}
	\\
        &\leq  C_a\exp\paren[\Big]{-\frac{\hat{\gamma}}{4\epsilon}},
\end{align*}
which implies \eqref{eq: eigenvector in Bi}.
\end{proof}

The above gives a bound for~$\psi_{2,\epsilon}$ in~$B_1 \cup B_2$.
We will now bound~$\psi_{2, \epsilon}$ on~$\T^d\setminus (B_1\cup B_2)$ by bounding the~$L_\epsilon$-harmonic extension of~$\psi_{2, \epsilon}$ of the boundary data, and then bounding solutions to the inhomogeneous problem.
For simplicity of notation, define
\begin{equation}
    \tilde{\Omega}\defeq \T^d\setminus (B_1\cup B_2).
\end{equation}
Lemma~\ref{lem: boundedness of eigenfunction around local minima} can be used to immediately bound the~$L_\epsilon$-harmonic extensions of~$\psi_{2, \epsilon}$.

\begin{lemma}\label{corBoundednessF0k}
    For $i=1,2$, let $f^{(i)}_{0,\epsilon}$ be the solution to
   \begin{align}
     L_{\epsilon}f^{(i)}_{0,\epsilon}(y)&=0,  && y\in \tilde \Omega \\
      f^{(i)}_{0,\epsilon}(y) &= \psi_{2,\epsilon}(y),
	&& y\in\partial B_i\\
        f^{(i)}_{0,\epsilon}(y) &=0,
	&& y\in (\partial B_1\cup \partial B_2)\setminus \partial B_i.
\end{align}
  For every~$\epsilon$ satisfying~\eqref{eq:epsilonCri} and every~$y\in \tilde \Omega$ we have
\begin{equation}\label{eq:OutBi}
  \abs{f^{(i)}_{0,\epsilon}(y)}\leq \abs{a_{i,\epsilon}}+C_a\exp\paren[\Big]{-\frac{\hat{\gamma}}{4\epsilon}}.
\end{equation}
\end{lemma}
\begin{proof}
    Observe that $f^{(i)}_{0,\epsilon}$ satisfies $L_{\epsilon}f^{(i)}_{0,\epsilon}(y)=0$ on $\tilde \Omega$.
    Thus, by weak maximum principle \cite[Section 6.4.1, Theorem 1]{Evans}, 
        \begin{equation}
        \sup_{\tilde \Omega}\abs{f^{(i)}_{0,\epsilon}}
	  =\sup_{\partial B_1\cup\partial B_2}\abs{f^{(i)}_{0,\epsilon}}=\sup_{\partial B_1\cup\partial B_2}\abs{\psi_{2,\epsilon}}\\
	  \overset{\eqref{eq: eigenvector in Bi}}{\leq}\abs{a_{i,\epsilon}}+C_a\exp\paren[\Big]{-\frac{\hat{\gamma}}{4 \epsilon}},
    \end{equation}
    which implies the inequality \eqref{eq:OutBi}.
\end{proof}

We now bound the eigenfunction of~$L_\epsilon$ in~$\tilde \Omega$ with inhomogeneous Dirichlet boundary conditions specified by~$\psi_{2, \epsilon}$.
\begin{lemma}\label{lem: boundedness of eigenfunction in each region}
For $i=1,2$, let $f_{\lambda}$ solve
\begin{align}
  (L_{\epsilon}-\lambda_{2,\epsilon}) f_{\lambda}(y)&=0,  && y\in \tilde{\Omega}\\
     f_{\lambda}(y) &= \psi_{2,\epsilon},&& y\in\partial B_i\\
	f_{\lambda}(y) &=0, && y\in (\partial B_1\cup \partial B_2)\setminus \partial B_i
      .
\end{align}
There exist $\epsilon_0>0$, $T_0'>0$ such that for  
\begin{equation}\label{eq:epsilonthreshold}
  \epsilon\leq \min\set[\Big]{
    \epsilon_0,
    \frac{ \hat{\gamma}}{4 \log(2\tilde{C} T_0')},
    \frac{\min\{R_1,R_2\}}{12}, \tilde{\epsilon}}\defeq \epsilon_1,
\end{equation}
we have
    \begin{equation}\label{eq:eigenOutBi}
        |f_{\lambda}(y)|\leq 2\paren[\Big]{\abs{a_{i,\epsilon}}+C_a\exp\paren[\Big]{-\frac{ \hat{\gamma}}{4 \epsilon}}},\quad \forall y\in \tilde{\Omega}.
    \end{equation}
  Here $C_a$ is the constant in \eqref{eq: eigenvector in Bi} and $\tilde{C}$ is the constant in~\eqref{eq:egvalub}.
\end{lemma}
\begin{proof}[Proof of Lemma \ref{lem: boundedness of eigenfunction in each region}] 
We only prove when $i=1$. The proof in the case $i=2$ is identical. Let $f_{0}$ solve
    \begin{alignat}{2}
      L_{\epsilon} f_{0}(y)&=0,  &\qquad& y\in \tilde{\Omega}\\
      f_{0}(y) &= \psi_{2,\epsilon},&& y\in\partial B_1\\
	f_{0}(y) &=0, && y\in \partial B_2.
\end{alignat}
Let $\delta f_{\lambda}=f_{\lambda}-f_{0}$. Then $\delta f_{\lambda}$ satisfies
\begin{alignat}{2}
  L_{\epsilon}\delta f_{\lambda}(y)&=\lambda_{2,\epsilon} f_{\lambda}(y),
    &\qquad& y\in \tilde{\Omega}\\
  \delta f_{\lambda}(y) &= 0,
    && y\in\partial B_1\cup \partial B_2
    .
\end{alignat}

Let~$\tau$ be the first exit time of~$X^\epsilon$ to~$B_1 \cup B_2$.
  We know~$g(y) \defeq \E^y \tau$ solves the Poisson equation
  \begin{alignat}{2}\label{eq: exit problem of B1 and B2}
    L_{\epsilon}g_\epsilon&=1,  &\qquad& y\in \tilde{\Omega}\\
      g_\epsilon &= 0,&& y\in\partial B_1\cup \partial B_2.
\end{alignat}
  Thus if $M'\defeq \sup_{z\in\tilde{\Omega}}|f_{\lambda}(z)|<\infty$,
  the comparison principle immediately implies
  \begin{equation}
    \sup_{\tilde{\Omega}} \abs{\delta f_\lambda}
      \leq \lambda_{2, \epsilon} M' \sup_{\tilde{\Omega}} g_\epsilon
    .
  \end{equation}

Since~$f_\lambda = \delta f_\lambda + f_0$, we see
\begin{equation}\label{eq: M}
  M'\leq \norm{f_{0}}_{L^\infty(\tilde \Omega)}
   +\lambda_{2,\epsilon} M' \norm{g_\epsilon}_{L^\infty(\tilde \Omega)}
    .
\end{equation}

According to \cite[Corollary of Lemma 1.9, Chapter 6]{FreidlinWentzell}, for $\epsilon$ smaller than some $\epsilon_0$, there exist constants $T_0$ and $c$ such that
    \begin{equation}\label{eq:defT0prime}
        \sup_{y\in\tilde{\Omega}}\E^y\tau
	  \leq T_0+\frac{\epsilon^2}{c}
	  <T_0+\frac{\epsilon_0^2}{c}
	  \defeq T_0'.
    \end{equation}
Notice that the choice \eqref{eq:epsilonthreshold} ensures that
  $\epsilon\leq\epsilon_0$,
  the condition~\eqref{eq:epsilonCri} holds, and
\begin{equation}\label{eq:epsilonEtau}
  \lambda_{2,\epsilon} T_0'\overset{\eqref{eq:egvalub}}{\leq} \tilde{C}\exp\paren[\Big]{-\frac{\hat{\gamma}}{4\epsilon}}T_0'\overset{\eqref{eq:epsilonthreshold}}{\leq}\frac12
  .
\end{equation}
This implies
\begin{align}
    M &\overset{\mathclap{\eqref{eq: M}}}{\leq}
      \frac{ \norm{f_{0}}_{L^\infty(\tilde \Omega)}  }{1-\lambda_{2,\epsilon} \norm{g_\epsilon}_{L^\infty(\tilde \Omega)}}
      \overset{\eqref{eq:defT0prime}, \eqref{eq:epsilonEtau}}{\leq}
      \frac{ \norm{f_{0}}_{L^\infty(\tilde \Omega)} }{1-\lambda_{2,\epsilon} T_0'}\\
    &\overset{\mathclap{\eqref{eq:OutBi},\eqref{eq:epsilonEtau}}}{\leq}\qquad
      2\paren[\Big]{\abs{a_{i,\epsilon}}+C_a\exp\paren[\Big]{-\frac{\hat{\gamma}}{4\epsilon}}}. \qedhere
\end{align}
\end{proof}

\begin{proof}[Proof of Proposition \ref{lem: uniform boundedness of eigenfunction}]
We will show that there exists constant~$C_{\psi} = C_{\psi}(U,d,C_m)$ independent of $\epsilon$ 
    such that~\eqref{eq:defCpsi} holds. We discuss two cases, $\epsilon\leq \epsilon_1$ and $\epsilon> \epsilon_1$, where $\epsilon_1$ is defined in \eqref{eq:epsilonCri}. 

\restartcases

\case[$\epsilon\leq \epsilon_1$]
For $y\in B_1\cup B_2$, we apply Lemma \ref{lem: boundedness of eigenfunction around local minima}, to obtain
\begin{equation}\label{ESinBi}
    \sup_{y\in B_1\cup B_2}|\psi_{2,\epsilon}(y)|\overset{\eqref{eq: eigenvector in Bi}}{\leq} \max\{|a_{1,\epsilon}|,|a_{2,\epsilon}|\}+C_a\exp\paren[\Big]{-\frac{3\hat{\gamma}}{4\epsilon}}\leq C_{m}+C_a.
\end{equation}
To obtain the last inequality above we used the fact that by Cauchy-Schwarz inequality
\begin{equation}\label{eq:boundA}
|a_{i,\epsilon}|\overset{\eqref{eq:aksol}}{=} \abs[\Big]{ \int_{\T^d}\frac{1}{\pi_{\epsilon}(\Omega_i)}\one_{\Omega_i}\psi_{2,\epsilon}\pi_{\epsilon}\, d x}\leq \frac{1}{\pi_{\epsilon}(\Omega_i)}\norm{\psi_{2,\epsilon}}_{L^2{\pi_\epsilon}}(\pi_{\epsilon}(\Omega_i))^{\frac12}\overset{\eqref{e:massRatioBound}}{\leq}C_{m}.
\end{equation}

For $y\in \tilde{\Omega}$, we write 
\begin{equation*}
    \psi_{2,\epsilon}=\psi^{(1)}_{2,\epsilon}+\psi^{(2)}_{2,\epsilon},
\end{equation*}
where $\psi^{(i)}_{2,\epsilon}$ solves
\begin{align}
  (L_{\epsilon}-\lambda_{2,\epsilon})\psi^{(i)}_{2,\epsilon}(y)&=0,  && y\in \tilde{\Omega}\\
      \psi^{(i)}_{2,\epsilon}(y) &= \psi_{2,\epsilon},&& y\in\partial B_i\\
	 \psi^{(i)}_{2,\epsilon}(y)&=0, && y\in (\partial B_1\cup \partial B_2)\setminus \partial B_i
      .
\end{align}
Applying Lemma~\ref{lem: boundedness of eigenfunction in each region} to $\psi^{(i)}_{2,\epsilon}$ gives
\begin{equation}
    \sup_{y\in \tilde{\Omega}}|\psi_{2,\epsilon}(y)|\overset{\eqref{eq:eigenOutBi}}{\leq}4\paren[\Big]{\max\{\abs{a_{1,\epsilon}},\abs{a_{2,\epsilon}}\}+C_a\exp\paren[\Big]{-\frac{\hat{\gamma}}{4\epsilon}}}
    \overset{\eqref{eq:boundA}}{\leq} 4(C_m+C_a).  \label{ESoutBi}
    \end{equation}

Combining \eqref{ESoutBi} and \eqref{ESinBi}, we obtain that for $0<\epsilon\leq \epsilon_1$, there exists $C= C(U,d,C_m)$ independent of $\epsilon$ such that $\norm{\psi_{2,\epsilon}}_{L^{\infty}(\mathbb{T}^d)}\leq C.$

\case[$\epsilon> \epsilon_1$]
According to \cite[Corollary 9.21]{GT}, for $y\in \mathbb{T}^d$,
\begin{align*}
    \sup_{x\in B(y,\epsilon)} |\psi_{2,\epsilon}(x)|
        &\leq \paren[\Big]{\frac{C}{|B(y,2\epsilon)|}\int_{B(y,2\epsilon)}|\psi_{2,\epsilon}(x)|^2 \, d x}^{\frac12}
       \\
        &\leq \paren[\Big]{\frac{C}{|B(y,2\epsilon)|}\Big(\sup_{z\in \mathbb{T}^d}e^{\frac{U(z)-U_{\min}}{\epsilon}}\Big)\int_{\mathbb{T}^d}|\psi_{2,\epsilon}(x)|^2 \, d \pi_{\epsilon}(x)}^{\frac12}\\
        &=C(\epsilon_1)^{-\frac{d}{2}}\exp\paren[\Big]{\frac{\norm{U}_{\osc}}{2\epsilon_1}}= C(U,d, C_m).
\end{align*}

We conclude from the above two cases that  \eqref{eq:defCpsi} holds.
\end{proof}

\section{The multi-well case}\label{sec:multi}

In this section, sketch how Theorem~\ref{thm: main} can be adapted to the setting where~$U$ has more than two local minima. Assumption~\ref{a:criticalpts} and~\ref{a:massRatioBound} remain the same, except that now ~$U$ has $J\geq 2$ local minima located at~$x_{\min, 1},\dots, x_{\min, J}$.
The nondegeneracy assumption Assumption~\ref{assumption: nondegeneracy} needs to be modified as we now explain.
Given~$i, j \in \set{1, \dots, J}$ with~$i \neq j$, define the saddle height between~$x_{\min,i}$ and $x_{\min,j}$ by
\begin{equation}\label{e:UHatDef1}
    \hat{U} = \hat{U}(x_{\min,i},x_{\min,j}) \defeq \inf_\omega \sup\limits_{t\in [0,1]}U(\omega(t))
    ,
\end{equation}
where the infimum is taken over all continuous paths~$\omega \in C( [0, 1]; \mathbb T^d)$ such that~$\omega(0) =x_{\min,1}$, $\omega(1) = x_{\min,2}$.
(Note that this generalizes the definition of the saddle height in~\eqref{e:UHatDef} for two well case.)
The generalization of Assumption~\ref{assumption: nondegeneracy} can now be stated as follows.

\begin{assumption}\label{assumption: nondegeneracyM}
   There exists $\delta>0$ such that the following hold.
   \begin{enumerate}[(1)]
     \item 
       For every~$i, j \in \set{1, \dots, J}$ with~$i \neq j$, the saddle height between the two local minima $x_{\min,i}$ and $x_{\min,j}$ is attained at a unique critical point $s_{i,j}$ of index one.
       That is, 
       the first eigenvalue of $\Hess U(s_{i,j})$ is negative and the others are positive. The point $s_{i,j}$ is called communicating saddle between the minima $x_{\min,i}$ and $x_{\min,j}$.

 \item  The set of local minima~$\set{x_{\min,1},\dots, x_{\min,J}}$ is ordered such that $x_{\min,1}$
is a global minimum and  the energy drop between the saddle and the local minimum is the largest for the second minimum in the following quantitative sense: for all $i\in\set{3,\dots,J}$
\begin{equation*}
    U(s_{1,2})- U(x_{\min,2}) \geq U(s_{1,i})- U(x_{\min,i}) +\delta.
\end{equation*}
   \end{enumerate}
\end{assumption}

This assumption is the same as Assumption~1.7 in~\cite{MenzSchlichting14}, which guarantees that the estimate \eqref{e:egvalLEpsilon} still holds in the multi-modal case.

Next, we state the version of Theorem~\ref{thm: main} that applies to $U$ with multiple wells.
\begin{theorem}[Multiple modes]\label{thm: mainM}
  Suppose for some~$0 \leq \eta_{\min} < \eta_{\max} \leq \infty$, the function~$U$ satisfies 
Assumptions~\ref{a:criticalpts},~\ref{a:massRatioBound} and~\ref{assumption: nondegeneracyM}.
  Let~$\hat \gamma_r \geq 1$ be defined as in~\eqref{e:gammaHatRDef} with~$\hat\gamma$ defined by~\eqref{e:gammaHatDefM} and~$\hat{U}$ defined by~\eqref{e:UHatDef}. In the same setting as in Theorem~\ref{thm: main}, for every bounded test function~$h$ and arbitrary initial points~$\set{y^i_1}$, the points~$(x^1, \dots, x^N)$ returned by Algorithm~\ref{a:ASMC}  satisfy that~\eqref{e:MCerror}.
\end{theorem}

We now sketch the proof of Theorem~\ref{thm: mainM}.
We begin with stating the multi-well version of the required Properties~\ref{p:spectral},~\ref{p:var} and~\ref{property:eigenfunctions}. In particular, we have the following proposition including all the required properties, whose proofs are analogous to those of Proposition~\ref{l: lower_bound_next_eigenvalue},~\ref{prop:pvar},~\ref{lem: uniform boundedness of eigenfunction} and Lemma~\ref{lem: mass in well as BV}, respectively.
\begin{proposition}[Properties]
    Suppose that the potential~$U$ satisfies
Assumptions~\ref{a:criticalpts},~\ref{a:massRatioBound} and~\ref{assumption: nondegeneracyM}, then the following properties hold:
\begin{enumerate}
    \item There exists~$\Lambda$ such that for all~$\epsilon \in (0, 1]$ such that
  \begin{equation}
    \lambda_{i,\epsilon}\geq \Lambda,\quad 
    \text{for all } i\geq J+1
    .
  \end{equation}
  Moreover, for $H > \hat U$, there exist a constant~$C_H$ such that for every~$\epsilon \in (0, 1]$, the bound~\eqref{e:egvalLEpsilon} holds.

\item For every~$i \in \set{2, \dots, J}$ and every~$\gamma < \hat \gamma$, for every~$0 < \epsilon' < \epsilon \leq 1$ we have
  \begin{gather}
    \label{e:PsiEpPiM}
     \norm{P_{F_{\epsilon}}(\psi_{i,\epsilon'}r_{\epsilon})}_{L^2(\pi_{\epsilon})}
      \leq 1+\frac{1}{2}C_m^2\sum_{j=1}^{J}\abs[\big]{\pi_{\epsilon'}(\Omega_j)-\pi_{\epsilon}(\Omega_j)}+
    C_r^{\frac12}C_{\gamma}\exp\paren[\Big]{-\frac{\gamma}{\epsilon}},
  \\
    \label{e:PsiEpEpPM}
    \norm{P_{F_{\epsilon}}(r_{\epsilon}-1)}_{L^2(\pi_{\epsilon})}
    \leq C_m\sum_{j=1}^{J}\abs[\big]{\pi_{\epsilon'}(\Omega_j)-\pi_{\epsilon}(\Omega_j)}+C_r^{\frac12}C_{\gamma}\exp\paren[\Big]{-\frac{\gamma}{\epsilon}}
    .
  \end{gather}
Here~$C_{r}$ be the constants defined in \eqref{eq:defCr},~$C_m$ is the constant in~\eqref{e:massRatioBound} and ~$C_{\gamma}$ is the constant in~\eqref{e:dEF}.

    \item  There exists constant~$C_{\psi}$, independent of~$\epsilon$ such that for all~$2\leq i\leq J$,
  \begin{equation}\label{eq:defCpsiM}
    \sup\limits_{0< \epsilon\leq 1}\|\psi_{i,\epsilon}\|_{L^{\infty}(\mathbb{T}^d)}\leq C_{\psi}.
  \end{equation}

  \item There exists a constant $C_{\BV}$ such that such that for every $\eta\in (0,1)$ and every $i\in \set{1,\dots, J}$, the bound~\eqref{eq:massBV} holds.
\end{enumerate}
\end{proposition}

In the double-well case, we estimated the error using the second eigenfunctions~$\psi_{2, \epsilon}$.
In the~$J$-well case, we will do the same, but use all the eigenfunctions~$\set{ \psi_{j, \epsilon} \st 2 \leq j \leq J}$.
More precisely, we have the following multi-modal version of Lemma~\ref{l:langevinError}, whose proof is very similar to that of Lemma~\ref{l:langevinError}.

\begin{lemma}\label{l:langevinErrorM}
   Assume that for each $i \in \set{1, \dots, N}$, $\PDF(Y_{\epsilon,0}^{i}) = q^{i}_{\epsilon,0}$.  
Then  for any bounded test function $h$,
 \begin{equation}\label{e:langevinErrorM}
    \Err_{\epsilon,T}(h)
    \leq e^{-\lambda_{2, \epsilon} T}  \norm{P_{F_{\epsilon}}(h)}_{L^2(\pi_{\epsilon})}\max_{2\leq j\leq J}\Err_{\epsilon,0}(\psi_{j,\epsilon})
    +\frac{1}{2\sqrt{N}}\|h\|_{\osc}
    + \mathcal E_{\epsilon,T}(h)
\end{equation}
where $\mathcal E_{\epsilon,T}(h)$ is defined as in~\eqref{eq: def of tilde epsilon}.
\end{lemma}

Next, we introduce the multi-modal versions of the remaining preparatory lemmas in Section~\ref{s:mainProof}.
Of these, Lemma~\ref{lem: probability density L infty norm between levels} remains unchanged.
Lemma~\ref{lem: first level} is almost identical except that we now require~\eqref{e:ErrPsi2} now holds for~$\Err_{2,0}(\psi_{j,2})$ for all~$2\leq j\leq J$. The multi-modal version of Lemma~\ref{l:iteration} is stated below, and its proof is similar to that of Lemma~\ref{l:iteration}, modulo a more tedious computation.

\begin{lemma}\label{l:iterationM}
Fix $\delta>0$, let~$C_q = C_q(U, 1)$ be the constant from Lemma~\ref{lem: probability density L infty norm between levels} with~$T_0 = 1$. Define $\tilde{C}_{N}$ by~\eqref{eq:deftildeCN}. For any~$\alpha> 0$, let~$H$ and~$\gamma$ defined as in~\eqref{eq:Hgammachoice}. Choose~$C_{\alpha}$ as in~\eqref{eq:deftildeCalpha} and
\begin{equation}\label{eq:CbetaM}
    C_{\beta}\defeq \exp\paren[\bigg]{\paren[\Big]{\frac12C_m^2+C_rC_m}J C_{\BV}+C_r^{\frac12}C_{\gamma}(1+C_r)}.
\end{equation}
Here~$C_r$,~$C_m$,~$C_{\BV}$ and~$C_{\gamma}$ are the constants in~\eqref{eq:defCr},~\eqref{e:massRatioBound},~\eqref{eq:massBV} and~\eqref{e:dEF}, respectively.
If~$N,T$ satisfy~\eqref{eq:TNpropN}--\eqref{eq:TNpropT},
  then for every~$k \in \set{2, \dots, M-1}$, we have
  \begin{equation}\label{e:iterationM}
   \max_{2\leq j\leq J} \Err_{k+1,0}(\psi_{j,k+1})
       \leq
	\beta_k \max_{2\leq j\leq J}\Err_{k,0}(\psi_{j,k})
	+ c_k
  \end{equation}
  with constants~$\beta_k, c_k$ satisfying~\eqref{e:CBeta}.
\end{lemma}

The proof of Lemma~\ref{l:iterationM} is very similar to that of Lemma~\ref{l:iteration}.
In particular, instead of~\eqref{eq: defBetak}, the constants~$\beta_k$ will now be chosen according to
\begin{equation}
    \beta_k\defeq e^{-\lambda_{2, k} T} \max_{2\leq j\leq J}\Big(\norm{P_{F_{\eta_k}}(r_k-1)}_{L^2(\pi_k)}\cdot\|\psi_{j,k+1}\|_{L^{\infty}}
    \mathbin{+}\norm{P_{F_{\eta_k}}(\psi_{j,k}r_k)}_{L^2(\pi_k)}\Big).
\end{equation}
The upper bound of $c_k$ is almost identical to its bound in the double-well case (Lemma~\ref{lem: cj}). 

The constants~$C_N$ and~$C_T$ in Theorem~\ref{thm: mainM} are defined the same as in~\eqref{eq:defCN} and~\eqref{eq:defCT}.
The proof of Theorem~\ref{thm: mainM} is now similar to the proof of Theorem~\ref{thm: main}.

\bibliographystyle{halpha-abbrv}
\bibliography{reference,gautam-refs1,gautam-refs2,Ruiyu-refs1,preprints}
\end{document}

%% file: preamble.tex
\usepackage{marktext}
\usepackage{gimac}

\usepackage[english]{babel}
\usepackage[utf8]{inputenc}
\usepackage{algpseudocode, algorithm}

\newcommand{\tmix}{t_\mathrm{mix}}
\newcommand{\TV}{\mathrm{TV}}
\newcommand{\BV}{\mathrm{BV}}
\newcommand{\LBV}{\mathrm{LBV}}
\newcommand{\osc}{\mathrm{osc}}

\DeclareMathOperator{\PI}{PI}

\DeclarePairedDelimiter{\ip}{\langle}{\rangle}

\usepackage{colortbl}
\colorlet{tableHeader}{cyan!5}

\newtheorem{assumption}[theorem]{Assumption}
\newtheorem{property}[theorem]{Property}

\newcounter{broj}

\DeclareMathOperator{\Err}{Err}
\DeclareMathOperator{\PDF}{PDF}
\DeclareMathOperator{\Hess}{Hess}

\definecolor{mygreen}{rgb}{0.25,0.85,0.25}
\definecolor{dgrn}{rgb}{0.05,0.5,0.05}